\documentclass[onefignum,onetabnum]{siamart171218}

\makeatletter 
\@ifundefined{l@section}{\let\l@section\@gobble}{}
\@ifundefined{l@subsection}{\let\l@subsection\@gobble}{}
\@ifundefined{l@subsubsection}{\let\l@subsubsection\@gobble}{}
\@ifundefined{l@paragraph}{\let\l@paragraph\@gobble}{}
\@ifundefined{l@subparagraph}{\let\l@subparagraph\@gobble}{}
\@ifundefined{l@figure}{\let\l@figure\@gobble}{}
\@ifundefined{l@table}{\let\l@table\@gobble}{}
\makeatother

\usepackage{amsmath, array, mathtools, amssymb, pgfplots, cancel, listings, xspace, algorithm, algpseudocode, overpic, soul, listings, amsfonts, graphicx, color, threeparttable, tikz, tocloft}
\usepackage[margin=1in]{geometry}

\hypersetup{
  pdftitle={A Taxonomy of Numerical Differentiation Methods},
  pdfauthor={P. Komarov, F. van Breugel, J. N. Kutz}
}

\usetikzlibrary{positioning, shapes.geometric, arrows.meta, fit}

\tikzset{startstop/.style={rectangle, rounded corners, draw=black},
    process/.style={rectangle, aspect=2, draw=black},
    decision/.style={diamond, aspect=2, minimum width=6em, draw=black, inner sep=0pt, align=center},
    arrow/.style={thick, -{>[length=0.25em, width=0.4em]}}
} 

\pgfplotsset{compat=1.18} 

\newcommand{\textbox}[2][0em]{ 
  \parbox[c][1em][c]{5.5em}{
  \vspace{#1}
  \centering\linespread{0.85}\selectfont#2
  }
}
\newcommand{\tight}[2][\centering]{ 
  \parbox{\hsize}{\linespread{0.85}\selectfont#1#2}
}
\newcommand{\blhref}[2]{
  {\hypersetup{linkcolor=black}\hyperref[#1]{#2}}
}
\newcommand{\by}{\mathbf{y}}
\newcommand{\bu}{\mathbf{u}}
\newcommand{\bx}{\mathbf{x}}
\newcommand{\bth}{\boldsymbol\theta}
\newcommand{\hbx}{\mathbf{\hat{x}}}
\newcommand{\hby}{\mathbf{\hat{y}}}
\newcommand{\dbx}{\mathbf{\dot{x}}}
\newcommand{\hdbx}{\mathbf{\hat{\dot{x}}}}
\newcommand{\bw}{\mathbf{w}}
\newcommand{\bv}{\mathbf{v}}
\newcommand{\bc}{\mathbf{c}}
\newcommand{\be}{\mathbf{e}}

\newcommand{\tbw}{\mathbf{\tilde w}}
\newcommand{\tbz}{\mathbf{\tilde z}}
\newcommand{\bff}{\boldsymbol{f}}
\newcommand{\bh}{\boldsymbol{h}}
\newcommand{\bA}{\mathbf{A}}
\newcommand{\bB}{\mathbf{B}}
\newcommand{\bC}{\mathbf{C}}
\newcommand{\tbA}{\mathbf{\tilde A}}
\newcommand{\tbB}{\mathbf{\tilde B}}
\newcommand{\tbC}{\mathbf{\tilde C}}
\newcommand{\bQ}{\mathbf{Q}}
\newcommand{\bR}{\mathbf{R}}
\newcommand{\bP}{\mathbf{P}}
\newcommand{\bK}{\mathbf{K}}
\newcommand{\bM}{\mathbf{M}}
\newcommand{\bW}{\mathbf{W}}
\newcommand{\bL}{\mathbf{L}}
\newcommand{\bV}{\mathbf{V}}
\newcommand{\bSigma}{\boldsymbol{\Sigma}}
\newcommand{\bLambda}{\boldsymbol{\Lambda}}
\newcommand{\bmu}{\boldsymbol{\mu}}
\newcommand{\bbeta}{\boldsymbol{\beta}}
\newcommand{\rbchi}{\raisebox{0.95\depth}{\small$\boldsymbol\chi$}}

\newcommand{\rotvert}{\rotatebox[origin=c]{90}{$\vert$}}

\DeclareMathOperator*{\argmin}{argmin}
\DeclareMathOperator*{\argmax}{argmax}

\headers{Taxonomy of differentiation}{P. Komarov, F. van Breugel, J. N. Kutz}

\title{A taxonomy of numerical differentiation methods}

\author{Pavel Komarov\thanks{Department of Electrical and Computer Engineering, University of Washington, Seattle, WA
  (\email{pvlkmrv@uw.edu})}
\and Floris van Breugel\thanks{Department of Mechanical Engineering, University of Nevada, Reno, NV
  (\email{floris@unr.edu})}
\and J. Nathan Kutz$^{\ddagger,*}$}

\begin{document}

\maketitle
\begingroup\makeatletter 
\renewcommand{\thefootnote}{\fnsymbol{footnote}}
\long\def\@makefntext#1{\parindent 1em\noindent\hbox to1.8em{\hss$\m@th^{\@thefnmark}$}#1}
\footnotetext[3]{Autodesk Research, 6 Agar Street, London UK (\email{kutz@uw.edu})}
\endgroup

\begin{abstract}
Differentiation is a cornerstone of computing and data analysis in every discipline of science and engineering. Indeed, most fundamental physics laws are expressed as relationships between derivatives in space and time. However, derivatives are rarely directly measurable and must instead be computed, often from noisy, potentially corrupt data streams. There is a rich and broad literature of computational differentiation algorithms, but many impose extra constraints to work correctly, e.g.~periodic boundary conditions, or are compromised in the presence of noise or anomalies. It can therefore be challenging to select the method best-suited to any particular problem. Here we review a broad range of numerical methods for calculating derivatives, present important contextual considerations and choice points, compare relative advantages, and provide basic theory for each algorithm in order to assist users with the mathematical underpinnings. This serves as a practical guide to help scientists and engineers match methods to application domains. We also provide an open-source Python package, \texttt{PyNumDiff}, which contains a broad suite of methods for differentiating noisy data.
\end{abstract}


\begin{AMS} 
  65D25 (Numerical analysis: Numerical differentiation), 97M01 (Math Ed: Mathematical modeling, applications of mathematics: Instructional exposition / tutorial papers), 37M01 (Dynamical Systems: Approximation methods and numerical treatment of dynamical systems: Instructional exposition / tutorial papers)
\end{AMS}

{\footnotesize
\tableofcontents}

\section{Introduction}

Computing derivatives is a critically enabling task for almost any application in the engineering, physical, mathematical, and biological sciences. Mathematical models and governing equations often take the form of relationships between changing quantities, so whether the computation is of temporal or spatial data, simulated or measured, be it in the field of signal processing, system identification, or Machine Learning (ML), accurate derivative estimates are crucial. The great need in this area has fueled great innovation, resulting in a proliferation of differentiation methods specialized to take advantage of certain scenarios or system properties~\cite{savitzky1964smoothing, trefethen2000spectral, kutzbook, aravkin}. However, confronting so many options complicates algorithm selection.

Numerical algorithms including Finite Differences, Spectral Methods, and Finite Elements or Volumes are now mature technologies, with the latter powering commercialized (e.g.~COMSOL~\cite{comsol}, Ansys and Abaqus) and open-source (e.g.~FEniCS, Firedrake, OpenFOAM) PDE solvers. Likewise, Automatic Differentiation (\texttt{AutoDiff}) dominates in Machine Learning~\cite{jax2018github}. As practitioners in the rapidly-growing field of data-driven science, we also take particular interest in derivative estimation for \textit{noisy} data streams. Motivated by the definition of the derivative from calculus, beginners in numerical computation often use finite differences by default, but this turns out to be one of the worst choices for real data. We attempt to give comprehensive coverage to all these schemes, from the narrowly focused to the most general-purpose. The overall trend is that more constrictive situations allow for better performance, so we suggest the use of specialized tools if possible. Or, if faced with a dearth of simplifying assumptions, we have made it easy for the reader to skip ahead to the pertinent category.

This paper is organized as follows: In \autoref{sec:overview} we establish a common map for navigating key determining factors, so practitioners can quickly get a sense of which archetypal situation is most representative of their differentiation problem and see at a glance which solution strategies might work best. Then we follow that map, beginning with the easiest situations and most bulletproof methods, working toward more equivocal, ambiguous scenarios. In \autoref{sec:autodiff} we cover Automatic Differentiation tools like \texttt{JAX}~\cite{jax2018github}, which have been incredibly successful and revolutionized deep learning, but are by nature only narrowly applicable. We also touch on their use in Differentiable Physics simulators. In \autoref{sec:noiseless-simulation} we address noiseless simulation, including Finite Differencing, Spectral Methods, and Finite Elements. We then discuss noise as a concept in \autoref{sec:noise} and follow up in \autoref{sec:noisy-with-knowledge} with methods that can address uncertainty by taking advantage of additional system knowledge in the form of dynamics models and synchronized data streams. The more difficult case of differentiating in the presence of noise \textit{without} additional information is covered in \autoref{sec:noisy-without-knowledge}, where we describe an approach for measuring performance in the absence of ground truth by balancing fidelity and smoothness~\cite{floris2020}, use this metric to choose near-optimal hyperparameters, and present experimental results to compare accuracy and bias across methods. In doing so, we make heavy use of tools from our companion, open-source Python package \texttt{PyNumDiff}~\cite{pynumdiff}, which implements many of the methods discussed as well as hyperparameter optimization to tune them. In \autoref{sec:practical} we cover the practical consideration of data with irregular sampling rate. Finally, \autoref{sec:conclusion} concludes with general retrospective and recommendations. Throughout, we draw connections to Scientific Machine Learning methods that enforce or ingest derivatives: Physics-Informed Neural Networks (Sections \ref{sec:autodiff} and \ref{sec:noisy-with-knowledge}), Operator Learning (Section \ref{sec:noiseless-simulation}), and Sparse Model Discovery (Section \ref{sec:noisy-without-knowledge}).

\section{Choosing a Method}
\label{sec:overview}

\begin{figure}\label{fig:winners}
\centering
\vspace{2mm}
\begin{tikzpicture}
\def\h{2.6cm}
\foreach \i in {0,...,5} { 
    \draw (\i*\h,0) -- (\i*\h,\h);
}
\foreach \j in {0,...,1} { 
    \draw (0,\j*\h) -- (5*\h,\j*\h);
}
\def\colheaders{
    Analytic\\[-3pt]Relationships,
    Simple\\[-3pt]Simulations,
    Complicated\\[-3pt]Simulations,
    Noisy Data\\[-3pt]with Model,
    Model-Free\\[-3pt]Noisy Data}
\foreach \head [count=\i from 0] in \colheaders {
    \node[anchor=south, align=center, text height=1ex, text depth=0ex, font=\bfseries\small] at ({(\i+0.5)*\h},\h) {\head};
}
\node[align=left] at (0.5*\h, 0.5*\h) {\blhref{sec:autodiff}{\texttt{AutoDiff}}\\[-3pt]\footnotesize to numerical\\[-4pt]\footnotesize precision};
\node[align=left] at (1.5*\h, 0.72*\h) {\blhref{sec:spectral}{Spectral}\\[-3pt]\footnotesize with assumptions,\\[-3pt]\footnotesize $O(\Delta x^\infty)$ error};
\node[align=left] at (1.45*\h, 0.25*\h) {\blhref{sec:finite-difference}{Finite Diff.}\\[-3.5pt]\footnotesize no assumptions,\\[-3pt]\footnotesize $O(\Delta x^m)$ error};
\node[align=left] at (2.51*\h, 0.5*\h) {\blhref{sec:finite-elements}{Finite Element}\\[-2pt]\blhref{sec:finite-elements}{Method}\\[-2.5pt]\footnotesize $O(h^\alpha)$ error};
\node[align=left] at (3.5*\h, 0.5*\h) {\blhref{sec:noisy-with-knowledge}{Kalman}\\[-2pt]\blhref{sec:noisy-with-knowledge}{Smoothers}\\[-3.5pt]\footnotesize maximum\\[-4pt]\footnotesize \textit{a posteriori}};
\node[align=left] at (4.5*\h, 0.8*\h) {\blhref{sec:kalman-constant-deriv}{\texttt{RTSDiff}}\\[-4pt]\footnotesize for irregular $\Delta t$};
\node[align=left] at (4.47*\h, 0.49*\h) {\blhref{sec:kalman-constant-deriv}{\texttt{RobustDiff}}\\[-3.5pt]\footnotesize for outliers};
\node[align=left] at (4.48*\h, 0.17*\h) {\blhref{sec:local-linear}{\texttt{LinearDiff}}\\[-3pt]\footnotesize for large $\Delta t$};
\end{tikzpicture}
\vspace{-2mm}
\caption{Preferred differentiation algorithms at a glance, for the five major situations identified in this review. Names are shorthand, with details given in later sections, hyperlinked for convenience. Further details and rationale behind these distinctions and selections are mapped in \autoref{fig:flowchart}. $\Delta x$ and $\Delta t$ are spacing between samples, $m$ is FD scheme order, $h$ is element side-length, and $\alpha$ is a constant, making both the FEM and FD error bounds similarly algebraic, while error of Spectral Methods decreases ``super-algebraically" as the number of samples and basis functions increases.}
\end{figure}

\begin{figure}[t]\label{fig:flowchart}
\centering
\begin{tikzpicture}[scale=0.85, transform shape, node distance=2.2em, every node/.style={font=\small}]
    \node (start) [startstop] {Start};
    \node (static-analytic?) [decision, right=of start] {\textbox[-0.2em]{static, analytic?}};
    \node (autodiff) [startstop, right=of static-analytic?, text width=6em, draw=blue, thick] {\tight{\blhref{sec:autodiff}{Automatic Differentiation}}};
    \node (simulation?) [decision, below=of static-analytic?] {\textbox[-0.5em]{noiseless simulation?}};
    \node (sim-smooth?) [decision, right=of simulation?, xshift=2em] {\textbox[-0.4em]{smooth derivative?}};
    \node (regular?) [decision, right=of sim-smooth?] {\textbox[-0.2em]{regular domain?}};
    \node (noisy-smooth?) [decision, below=of simulation?] {\textbox[-0.4em]{smooth derivative?}};
    \node (sim-periodic?) [decision, right=of regular?] {\textbox[0.2em]{periodic BC?}};
    \node (bad) [startstop, below=of noisy-smooth?, yshift=1em, text width=8.5em, align=center, draw=red, thick] {\tight{Can't distinguish noise from derivative}};
    \node (finite-elements) [startstop, below=of regular?, text width=4.3em, align=center, draw=violet, thick] {\tight{\blhref{sec:finite-elements}{Finite Elements}}};
    \node (sim-fourier) [startstop, right=of sim-periodic?, text width=3.3em, draw=blue, thick] {\tight{\blhref{sec:fourier}{Fourier\\Spectral}}};
    \node (cosine-sample?) [decision, below=of sim-periodic?] {\textbox[-0.1em]{can sample\\cos-spaced?}};
    \node (chebyshev) [startstop, right=of cosine-sample?, text width=4.3em, draw=blue, thick] {\tight{\blhref{sec:chebyshev}{Chebyshev\\Spectral}}};
    \node (finite-difference) [startstop, below=of cosine-sample?, draw=violet, thick] {\blhref{sec:finite-difference}{Finite Difference}};
    \node (noisy-periodic?) [decision, below=of noisy-smooth?, xshift=8em, yshift=-2em] {\textbox[0.2em]{periodic BC?}};
    \node (noisy-fourier) [startstop, right=of noisy-periodic?, draw=violet, thick] {\blhref{sec:noise-recs}{Fourier Spectral}};
    \node (known-dynamics?) [decision, below=of noisy-periodic?] {\textbox[-0.2em]{known dynamics?}};
    \node (linear?) [decision, right=of known-dynamics?, xshift=2em] {\textbox{linear?}};
    \node (nice-noise?) [decision, right=of linear?] {\textbox[-0.4em]{Gaussian white noise?}};
    \node (kalman) [startstop, right=of nice-noise?, text width=4em, align=center, draw=violet, thick] {\tight{\blhref{sec:kalman-filter}{Kalman} \blhref{sec:rauch-tung-striebel}{Smoother}}};
    \node (robust) [startstop, below=of nice-noise?, text width=4em, align=center, draw=violet, thick] {\tight{\blhref{sec:robust-estimation}{Robust Kalman Smoother}}};
    \node (nonlinear-KF) [startstop, below=of linear?, xshift=2em, text width=4.1em, align=center, draw=purple, thick] {\tight{\blhref{sec:UKF}{Nonlinear Kalman Smoother}}};
    \node (irregular-steps?) [decision, below=of known-dynamics?] {\textbox[0.2em]{\blhref{sec:practical}{irregular steps?}}};
    \node (smooth-fd) [process, below=of irregular-steps?, xshift=-3em, align=center, text width=4.5em] {\tight{\blhref{sec:prefilt}{Kernel Smoothing\\$\rightarrow$ Finite Difference}}};
    \node (iterated-fd) [process, right=of smooth-fd, xshift=-2em, align=center, text width=4.1em] {\tight{\blhref{sec:iteratedfinitedifference}{Iterated Finite Difference}}};
    \node (polynomial-fits) [process, right=of iterated-fd, xshift=-2em, align=center, text width=4.5em] {\tight{\blhref{sec:polyfit}{Polynomial Fits}}};
    \node (basis-fits) [process, right=of polynomial-fits, xshift=-2em, align=center, text width=2.5em] {\tight{\blhref{sec:basis-fit}{Basis Fits}}};
    \node (tvr) [process, right=of basis-fits, xshift=-2em, align=center, text width=4.8em] {\tight{\blhref{sec:tvr}{Total Variation Regularized Derivative}}};
    \node (rts) [process, right=of tvr, xshift=-2em, align=center, text width=6.5em] {\tight{\blhref{sec:kalman-constant-deriv}{Kalman Smoothing with Naive Model}}};
    \node (ll) [process, right=of rts, xshift=-2em, align=center, text width=5.2em] {\tight{\blhref{sec:local-linear}{Learned Local Linear Dynamics}}};
    \node (smoothing-methods) [draw, fit=(smooth-fd)(polynomial-fits)(basis-fits)(iterated-fd)(tvr)(rts)(ll), rounded corners, draw=purple, thick] {};

    \draw[arrow] (start) -- (static-analytic?);
    \draw[arrow] (static-analytic?) -- node[anchor=south] {yes} (autodiff);
    \draw[arrow] (static-analytic?) -- node[anchor=east, text width=9em] {\tight[\raggedleft]{no, dynamic or non-analytic solutions}} (simulation?);
    \draw[arrow] (simulation?) -- node[anchor=south, text width=6em] {\tight{yes, governing equations}} (sim-smooth?);
    \draw[arrow] (simulation?) -- node[anchor=east] {no, measurements} (noisy-smooth?);
    \draw[arrow] (noisy-smooth?) -- node[anchor=east] {no} (bad);
    \draw[arrow] (sim-smooth?) |- node[anchor=south west] {no} (finite-elements);
    \draw[arrow] (sim-smooth?) -- node[anchor=south] {yes} (regular?);
    \draw[arrow] (regular?) -- node[anchor=south] {yes} (sim-periodic?);
    \draw[arrow] (regular?) -- node[anchor=west] {no} (finite-elements);
    \draw[arrow] (noisy-smooth?) -| node[anchor=south] {yes} node[anchor=west, yshift=-3em, text width=8em] {\tight[\raggedright]{\blhref{sec:noise}{noisy signal,\\smooth derivative}}} (noisy-periodic?);
    \draw[arrow] (sim-periodic?) -- node[anchor=south] {yes} (sim-fourier);
    \draw[arrow] (sim-periodic?) -- node[anchor=west] {no} (cosine-sample?);
    \draw[arrow] (cosine-sample?) -- node[anchor=south] {yes} (chebyshev);
    \draw[arrow] (cosine-sample?) -- node[anchor=west] {no, uniform samples} (finite-difference);
    \draw[arrow] (noisy-periodic?) -- node[anchor=south] {yes} (noisy-fourier);
    \draw[arrow] (noisy-periodic?) -- node[anchor=east] {no} (known-dynamics?);
    \draw[arrow] (known-dynamics?) -- node (known-dynamics-yes) [anchor=south] {yes} (linear?);
    \draw[arrow] (linear?) -- node[anchor=south] {yes} (nice-noise?);
    \draw[arrow] (nice-noise?) -- node[anchor=south] {yes} (kalman);
    \draw[arrow] (nice-noise?) -- node (robust-no) [anchor=east] {no} (robust);
    \draw[arrow] (robust-no) -| node[anchor=north, text width=5em] {\tight{classic\\often works}} (kalman);
    \draw[arrow] (linear?) |- node[anchor=east] {no} ++(2em,-3em) -- (nonlinear-KF);
    \draw[arrow] (known-dynamics?) -- node (known-dynamics-no) [anchor=east] {no} (irregular-steps?);
    \draw[arrow] (known-dynamics-yes) |- node[anchor=north, yshift=0.1em, text width=7em] {\tight{modeling may not be worth it}} (known-dynamics-no);
    \draw[arrow] (irregular-steps?) -- node[anchor=east] {no} (irregular-steps? |- smoothing-methods.north);
    \draw[arrow] (irregular-steps?) -| ++(5em,-3em) -| node[anchor=south, xshift=-1.3em] {yes} (polynomial-fits);
    \draw[arrow] (irregular-steps?) -| ++(5em,-3em) -| (rts);
    \draw[arrow] (irregular-steps?) -| ++(5em,-3em) -| (ll);
\end{tikzpicture}
\caption{Suggested procedure for choosing a differentiation method. Methods are boxed with blue, violet, purple, or red to give an overall impression of how computationally efficient, accurate, robust, and user-friendly they are, with blue meaning a method excels in all respects for its specialization and red meaning a solution is hopeless. Groupings in the lowermost node represent many methods, expounded in \autoref{sec:noisy-without-knowledge}. Methods and a couple of other labels are hyperlinked for quick navigation to corresponding sections.}
\vspace{-2mm}
\end{figure}

In the broadest possible strokes, there are three scenarios for computing derivatives, two of which can be further refined:

\vspace{2mm}\begin{enumerate}
  \item Analytic relationships with unchanging structure
  \item Simulations: (a) less complicated, (b) more complicated
  \item Noisy data: (a) with a model, (b) without a model
\end{enumerate}\vspace{2mm}

\noindent A digest of these with some preferred algorithms is given in \autoref{fig:winners}. \autoref{fig:flowchart} shows a more complete and detailed overview of key considerations and methods, organized as a flowchart. The decisions therein are important because different techniques make different assumptions. For example, smoothing methods assume the signal is continuously differentiable up to some order (possibly $C^\infty$), and the Kalman filter requires a system model. Under changing assumptions, the most suitable approach changes, so it is important to understand whether data is periodic, noisy, etc. Notably missing among these considerations is dimensionality: Because differentiation is a linear operator, all methods can be extended to the multidimensional case, either by multivariable chain rule, by using higher-dimensional meshes and basis functions, or by sequential application along grid dimensions (in any order).\footnote{Not all methods are extensible with equal ease. On grids the curse of dimensionality dominates computational cost, while irregular domains require extra machinery like boundary geometry and mesh generation.} We therefore primarily treat the one-dimensional case in this review, for clarity of exposition.

For most differentiation scenarios (Sections \ref{sec:autodiff}, \ref{sec:noiseless-simulation}, \ref{sec:noisy-with-knowledge}), there are only a handful of clear winners that dominate other methods due to computational efficiency, accuracy, or the unique ability to handle relevant structure. The exceptional case is differentiation of noisy data without a model (\autoref{sec:noisy-without-knowledge}), because there are many possible underlying generating functions, making the problem ill-posed. Choosing a method becomes significantly more complicated in this situation, because so many methods have been invented, none foolproof. Methods are compared on particular properties throughout the review: error bounds at a glance in \autoref{fig:winners}, the decision path encompassing many factors in \autoref{fig:flowchart}, boundary handling for noiseless cases in Sections \ref{sec:finite-difference}, \ref{sec:fourier}, and \ref{sec:weak-form}, measured accuracy and bias in \autoref{sec:accuracy-and-bias}, computational complexity in \autoref{ta:complexity}, and tolerance of irregular sampling in \autoref{ta:irregular-dt-workability}.

Describing the mathematical underpinnings of all these methods, as well as the problems they fit, requires an extensive vocabulary of symbols, so a partial lexicon is given in \autoref{ta:common-symbols} to aid with orientation.

\begin{table}[!t]
\caption{\label{ta:common-symbols} Common symbols used throughout this review. For even more, used in \autoref{sec:noisy-with-knowledge} and beyond, see \autoref{ta:kalman-symbols}.}\vspace{-2mm}
\centering
\begin{tabular}{lp{10cm}}
  \hline $y$, $u$, $v$ & a function, although $y$ can also be a spatial independent variable\\
  $f$ & a function, also used to denote frequency\\
  $x$, $t$, $\Delta x$, $\Delta t$ & independent variable and its discrete increment\\
  $c$, $a$, $b$, $\alpha$, $\beta$ & constants or coefficients\\
  $\epsilon$ & small value\\
  $h$ & side-length of a finite element\\
  $n$ & data sample iterator\\
  $k$ & basis function or term iterator\\
  $j$ & generic iterator\\
  $N$ & number of total data samples\\
  $\xi$, $\varphi$ & basis functions\\
  $\nu$ & derivative order, i.e.~1$^\text{st}$, 2$^\text{nd}$, or 3$^\text{rd}$\\
  $\omega$ & angular frequency\\
  $\bth$ & parameters, stacked together as a vector\\
  $X$ & input data examples stacked together in a matrix or tensor\\
  $Y$ & target examples stacked together, or DFT coef.~with subscript $k$\\
  $L$ & a loss function, or occasionally a length\\ \hline
\end{tabular}
\vspace{-4mm}
\end{table}

\subsection{Recommendations}

When suitable, Automatic Differentiation (\autoref{sec:autodiff}) and Spectral Differentiation with the Fourier (\autoref{sec:fourier}) or Chebyshev basis (\autoref{sec:chebyshev}) are excellent, both accurate and efficient, but use of these is quite restrictive. The Finite Element Method (\autoref{sec:finite-elements}), by contrast, is exceptionally versatile, capable of differentiating even non-smooth functions (in the weak sense developed in \autoref{sec:weak-form}) over irregularly-shaped domains. But it is also rather mathematically involved, involves more setup, and is relatively computationally demanding---although much of its administrative overhead can be offloaded to commercial or open-source packages. Finite Difference (\autoref{sec:finite-difference}) is less accurate and robust but is very straightforward, works well in the absence of noise, and has well-defined error bounds~\cite{kutzbook}.

Down the ``noisy signal, smooth derivative" path, the Fourier spectral method (i.e.~transform, zero out higher modes, and inverse transform) is ideal for isolating high-frequency noise (\autoref{sec:noise}), but it should only be used when signals are periodic (can be joined smoothly end-to-end), lest Gibbs phenomenon (\autoref{fig:gibbs-phenomenon}) corrupt the approximation. Otherwise, one should exploit a model of system dynamics and noise when possible (\autoref{sec:noisy-with-knowledge}), although system models can be challenging to derive beyond simple systems and can be of limited utility unless fairly accurate. In the naive case (\autoref{sec:noisy-without-knowledge}), sophisticated methods tend to perform similarly (\autoref{sec:performance-comparison}), but additional considerations---such as robustness to outliers, fast hyperparameter optimization, and ability to handle large or variable step size---distinguish a few methods for specialized use (\autoref{fig:flexibility}). In some cases, a method may even be preferred for narrative reasons, having theoretical structure that matches assumptions about an underlying process or context.

\section{Differentiating Analytic Functions with Static Structure}
\label{sec:autodiff} 

Given the dramatic success of deep learning, optimized by hardware-accelerated gradient descent implemented in libraries like \texttt{JAX}~\cite{jax2018github} and \texttt{PyTorch}~\cite{pytorch}, it is natural to wonder how \texttt{AutoDiff} fits into the larger numerical differentiation ecosystem. In a word: \texttt{AutoDiff} is not for differentiating data samples, but rather for sampling fixed derivative relationships. Automatic Differentiation is most heavily used in the context of Machine Learning, where a model, typically named $f$ for ``function", is parameterized by a big vector $\bth$, takes a batch of inputs $X$, and produces a batch of outputs.

To illustrate the machinery, the paradigmatic supervised ML problem is:
$$\bth^* = \argmin_{\bth \in \mathbb{R}^m}\ L(f(X; \bth), Y)$$

\noindent where $L$ is a loss function, minimized when model outputs equal targets; $f$ is a model with possibly multidimensional input and output; $X$ is a set of inputs to the model; $Y$ is a set of corresponding target outputs; and $\bth$ affects how the model produces its answers, sometimes shown after a semicolon for visual separation from model inputs. Several other major classes of problem can be addressed with a similar setup, including unsupervised, semisupervised, and self-supervised learning~\cite{barlowtwins}, as well as fitting key evaluation functions in reinforcement learning, because ultimately all boil down to minimizing \textit{loss} based on known correct answers\footnote{A computer scientist we know jokes, ``Supervised learning is the only thing that really works."~\cite{abhishek}} or clever assumptions that enforce self-consistency and regularization~\cite{sindy-shred, parsimony}. In scientific computing, for example, Physics-Informed Neural Networks (PINNs) generalize with fewer training examples by exploiting a loss supplemented with residuals of known physical equations~\cite{pinn}, i.e.~$\|f(X;\bth)-Y\|_2^2 + \lambda\|\mathcal{R}[f(\cdot;\bth)](\tilde X)\|_2^2$, where $\mathcal{R}[\cdot]$ denotes an (analytic, possibly nonlinear) operator applied to the predicted physical quantity, written in residual form (so as to equal zero when the governing equation is satisfied),\footnote{As a concrete example, consider a mass on a damped spring under known forcing, $A\sin(\omega t)$: The model $f(t;\bth)$ predicts position at time $t$, $X$ holds the times at which position was measured, and $Y$ holds measurements. $\mathcal{R}[f](t) = m\ddot{f} + c\dot{f} + kf - A\sin(\omega t)$ for mass $m$, damping $c$, and spring constant $k$, which vanishes exactly when the equation of motion is satisfied, with $\dot{f}$ and $\ddot{f}$ supplied by \texttt{AutoDiff} through the network. Because the forcing is known in closed form, the residual can be imposed at as many additional times as one likes, so a handful of noisy position samples suffices to pin down a trajectory.} $\tilde X$ are freely chosen points where the residual is evaluated, and $\lambda$ balances the physical prior against observed data.

The key realization of ML is that we can consider $X$ and $Y$ to be static and given, and treat the model's parameterization, $\bth$, as an input, thus obtaining a relationship between overall error and the parameters. If we then differentiate this relationship with respect to the parameters and evaluate at the given input data, we obtain a gradient, the tilt of the loss function at the point encoded by current $\bth$, which can be used to navigate ``downhill" toward the loss function's minimum by updating model parameters with:

$$\bth_{j+1} = \bth_j - \epsilon \nabla_{\bth} L \Big|_{X,Y;\bth_j}$$

\noindent where $\epsilon$ is a learning rate and $\nabla_{\bth} L \big|_{X,Y;\bth_j}$ indicates the gradient of $L$, using the current parameterization of $f$, averaged over evaluations with input-output data pairs~\cite{goodfellow}.

The total size of the data may be massive, so we often choose to evaluate over small sub-samples of the data for computational efficiency, a practice known as \textit{stochastic} gradient descent. This leads to ``sampling noise", because the mean gradient over the sub-population is not exactly the same as the mean gradient evaluated over the full population. To handle this noise, we may use a more sophisticated update rule, like Adam~\cite{adam-optimizer}, which effectively serves as a low-pass filter to smooth out variations in the descent direction. For more on noise, see \autoref{sec:noise}.

To get a better sense of what taking the gradient of such a complicated construction means, we can expand the calculation using the multivariable chain rule:
$$\nabla_{\bth} L = \frac{d}{d\bth} L(f(\bx; \bth), \by) = \frac{\partial L}{\partial f} \frac{d f(\bx; \bth)}{d\bth} + \overset{\text{\normalsize 0}}{\cancel{\frac{\partial L}{\partial \by}\frac{d\by}{d\bth}}} = \frac{\partial L}{\partial f} \Big(\overset{\text{\normalsize 0}}{\cancel{\frac{\partial f}{\partial \bx} \frac{d\bx}{d\bth}}} + \frac{\partial f}{\partial\bth} \overset{\text{\normalsize 1}}{\cancel{\frac{d\bth}{d\bth}}} \Big) = \frac{\partial L}{\partial f}\frac{\partial f}{\partial\bth}$$

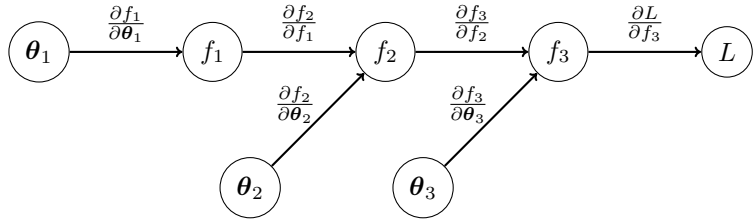
\begin{figure}[!t]\label{fig:computational-graph}
\centering
\begin{tikzpicture}[node distance=1.5cm]

\node[circle, draw] (f1) {$f_1$};
\node[circle, draw, right=of f1] (f2) {$f_2$};
\node[circle, draw, right=of f2] (f3) {$f_3$};
\node[circle, draw, right=of f3] (L) {$L$};
\node[circle, draw, left=of f1] (theta1) {$\bth_1$};
\node[circle, draw, below=of f1, xshift=5mm, yshift=5mm] (theta2) {$\bth_2$};
\node[circle, draw, below=of f2, xshift=5mm, yshift=5mm] (theta3) {$\bth_3$};

\draw[arrow] (theta1) -- (f1) node[midway, above] {$\frac{\partial f_1}{\partial \bth_1}$};
\draw[arrow] (f1) -- (f2) node[midway, above] {$\frac{\partial f_2}{\partial f_1}$};
\draw[arrow] (f2) -- (f3) node[midway, above] {$\frac{\partial f_3}{\partial f_2}$};
\draw[arrow] (f3) -- (L) node[midway, above] {$\frac{\partial L}{\partial f_3}$};
\draw[arrow] (theta2) -- (f2) node[midway, left, xshift=1mm, yshift=2mm] {$\frac{\partial f_2}{\partial \bth_2}$};
\draw[arrow] (theta3) -- (f3) node[midway, left, xshift=1mm, yshift=2mm] {$\frac{\partial f_3}{\partial \bth_3}$};

\end{tikzpicture}
\caption{Example of a simple computational graph for $\nabla_{\bth} L$ with compositional $f(\bx; \bth) = f_3(f_2(f_1(\bx; \bth_1);\bth_2);\bth_3)$, which results in $\frac{\partial f}{\partial \bth} = \frac{\partial f_3}{\partial f_2} \frac{\partial f_2}{\partial f_1} \frac{\partial f_1}{\partial \bth_1} + \frac{\partial f_3}{\partial f_2} \frac{\partial f_2}{\partial \bth_2} + \frac{\partial f_3}{\partial \bth_3}$ by chain rule.}
\end{figure}

If the model's evaluation function is compositional, as in the layers of a deep net, we perform the chain rule through all the layers. This involves differentiating vectors, matrices, or even tensors w.r.t.~other multidimensional objects, which can be somewhat confusing but is ultimately a pile of scalar derivatives stacked together. The details are well explained by~\cite{erik} and~\cite{clark}. The final product of this process is a relationship in terms of the loss function and model evaluation function, both of which are \textit{fixed}. \texttt{AutoDiff} exploits these unchanging relationships to find the analytic derivative (in many variables) of an analytic function in a very mechanistic way (computational graph, \autoref{fig:computational-graph}), then evaluates it at many data points in parallel to arbitrarily high accuracy.

\subsection{When Automatic Differentiation is Inappropriate}
\label{sec:wave-example}

\begin{figure}[t]\label{fig:wavesim}
  \centering
  \includegraphics[width=0.99\textwidth]{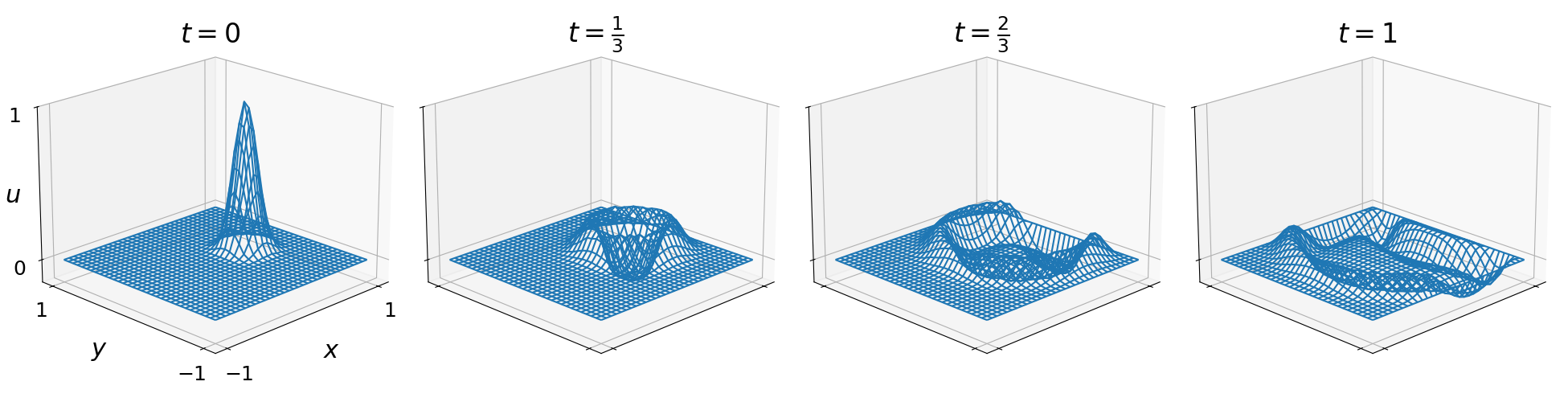}
  \vspace{-4mm}
  \caption{Simulation of the wave equation for initial condition $u(x, y, 0) = e^{-40((x - 0.4)^2 + y^2)}$, $u_t(x, y, 0) = 0$ and boundary condition $u = 0$.}
\end{figure}

To get a sense of why such a dominant tool as \texttt{AutoDiff} does not comprehensively handle our differentiation needs, consider simulating the 2D wave equation~\cite{trefethen2000spectral}, which obeys the partial differential equation (PDE):
$$\frac{\partial^2 u}{\partial t^2} = \frac{\partial^2 u}{\partial x^2} + \frac{\partial^2 u}{\partial y^2}$$

\noindent where $x$ and $y$ are dimensions of the input, and $u(x,y)$ is the wave amplitude. We choose the domain $-1 \leq x, y \leq 1$, boundary condition $u = 0$, and initial condition $u(x, y, 0) = e^{-40((x - 0.4)^2 + y^2)}, \ \frac{\partial}{\partial t}u(x, y, 0) = 0$, which evolves as illustrated in \autoref{fig:wavesim} through time.

To efficiently simulate the evolution of this system, we can discretize and use \autoref{algo:wavesim}, which employs the ``leapfrog" scheme, known to be stable for the second-order wave equation with small enough $\Delta t$~\cite{strang-notes, kutzbook}. At each time step, $u$ is updated by accumulating small numerical integrals of its second-order spatial derivatives. This means that although the initial condition of $u$ is simply analytic, simulated $u$ rapidly acquires more terms across future time steps. Furthermore, for arbitrary initial and boundary conditions, it is possible an analytic expression for $u$ does not exist.

\begin{algorithm}[t]
\caption{\textbf{Wave Simulation}}
\label{algo:wavesim}
\begin{algorithmic}[1] 
\State choose final time $T$, discrete time increment $\Delta t$, and spatial gaps $\Delta x$ and $\Delta y$
\State $u =$ sample the initial condition on a grid, where samples are $\Delta x$ and $\Delta y$ apart in respective dimensions
\State $u_\text{prev} = u$
\For{$t_n \in \{0, \Delta t, 2\Delta t,...,T\}$}
    \State $\text{right-hand side} = \Big[ \frac{\partial^2 u}{\partial x^2} + \frac{\partial^2 u}{\partial y^2} \Big] \Big|_{\text{grid points}}$ \Comment{evaluate at grid points}
    \State $u_\text{next} = 2u - u_\text{prev} + \Delta t^2\cdot\text{(right-hand side)}$ \Comment{leapfrog update}
    \State $u_\text{next} [\text{boundary}] = 0$\Comment{pin the boundary to 0}
    \State $u_\text{prev} = u$
    \State $u = u_\text{next}$
\EndFor
\end{algorithmic}
\end{algorithm}

Notice that the derivatives here are taken on a \textit{changing} function at \textit{fixed} sample points rather than on a \textit{fixed} relationship at \textit{changing} parameterization points. This distinction is subtle, because a parameterized relationship with sufficient expressive power can model a function of any shape~\cite{hornik1989}, so in principle we could define $u(x, y, t; \bth)$ to follow a changing function and use \texttt{AutoDiff} to evaluate its derivatives w.r.t.~$x$ and $y$ at space-time grid points. But this would either require knowledge of the true shape of $u(x_n, y_n)\ \forall\ t_n$ to use as a training target or a penalty on the PDE residual as in PINNs, and the whole point of simulation is to compute these answers using only the PDE relationship, without having to solve analytically or encode $u$ as a functional relationship when numerical values will do.

\subsection{Automatically Differentiable Physics}
\label{sec:autodiff-physics}

There is a case where \texttt{AutoDiff} can play an important role in simulations: If we assume a particular, unchanging form between a function's values and those of its derivative, like a Finite Difference equation (see \autoref{sec:finite-difference}) or FFT-based Spectral Method (see \autoref{algo:deriv-via-fft} and \autoref{algo:deriv-via-cheb}), then Automatic Differentiation can be used to take gradients of this relationship, a kind of meta-derivative, and thence by chain rule obtain gradients between initial and final states across many steps.

Several software libraries implement this architecture, including \texttt{JAX-Fluids}~\cite{jaxfluids-article, jaxfluids-library}, \texttt{PhiFlow}~\cite{phiflow}, \texttt{JAX-FEM}~\cite{jax-fem}, and \texttt{Exponax}~\cite{exponax}. ``Differentiable Physics"~\cite{differentiable-physics} refers to the passage of gradients all the way through a physics simulation, enabling seamless assimilation in ML training and inference pipelines. Note, however, that the physical derivatives \textit{per se} of such simulations are still found via whichever proximal differentiation algorithm the simulator implements (Finite Difference in the case of~\cite{jaxfluids-library, phiflow}), and so the simulation itself is only as accurate as that method.\footnote{\texttt{AutoDiff} through a simulator returns the exact gradient of the discrete system model, i.e.\ the \textit{discrete adjoint}, not of the continuous problem it approximates. A complementary route derives adjoint equations from the governing PDE, i.e.\ the \textit{continuous adjoint}, and discretizes and solves those instead. The two approaches generally disagree at finite resolution, converging to the same true gradient under refinement.}

\begin{figure}[!t]\label{fig:diff-phys}
  \centering
  \includegraphics[width=0.99\textwidth]{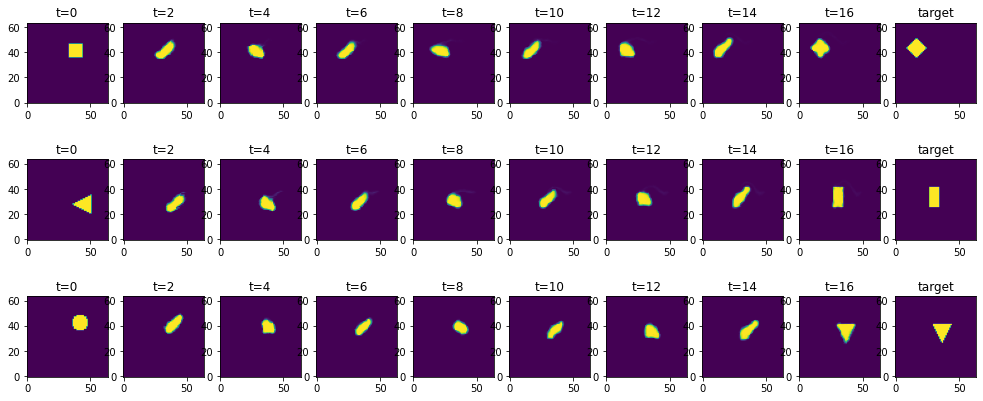}
  \vspace{-3mm}
  \caption{From~\cite{differentiable-physics}, used with permission, a neural network's solution to an advection flow inverse problem, demonstrating the strength of end-to-end training with Differentiable Physics. Fluid with marker dye is simulated at a particular location, and the goal is to control high-dimensional fluid flows to force the marker patch into a target shape.}
  \vspace{-3mm}
\end{figure}

\subsection{Recommendations} In some cases it may be possible to cast a problem into a static, analytic form. If so, use Automatic Differentiation without reservation, knowing it has been stress tested and scaled up by the prolific Machine Learning community. But recognize these are constrictive limitations, that for many problems analytic forms are unknown or inconvenient, necessitating a purely numeric representation of functions and other differentiation methods.

\section{Differentiating Noiseless Simulation Data}
\label{sec:noiseless-simulation}

When a function $y(x)$ ceases to be representable with few analytic terms, it is more conveniently represented numerically, in a vector. However, continuous functions are infinitely dense, so it is impossible to completely represent an arbitrary function in a finite vector. Instead, we must either sample $y_n = y(x_n),\ n \in \{0, ..., N-1\}$, over some domain $x \in [a, b]$, or keep a list of coefficients to direct reconstruction of $y$ from a set of basis functions. Because isolated values have undefined slope, we implicitly assume some interpolation to fill the gaps before differentiation~\cite{johnson-fft-notes}. Fortunately, in the absence of noise, function approximations of this kind can be made highly accurate.\footnote{These approximations also appear in Operator Learning, where a model consumes and produces \textit{functions} (e.g.~mapping between continuous system states through time), represented either by samples on a grid or coefficients in a basis~\cite{fno}. Such maps can be approximated linearly from a series of discretized snapshots, as in Dynamic Mode Decomposition (DMD)~\cite{dmd, optimized-dmd}, or, with enough data, built by optimization through a computational graph (\autoref{sec:autodiff}), i.e.~Neural Operators. The latter implementation enables incorporation of physical regularization terms like PINNs, except there is no coordinate input to differentiate against, so derivatives must be encoded \textit{structurally} by a stencil or transform inside the graph, as in Differentiable Physics simulators (\autoref{sec:autodiff-physics})~\cite{pino}.}

\subsection{Finite Difference}
\label{sec:finite-difference}

Consider the most basic definition of a derivative, which takes the difference of two nearby function values, divides by the distance between them, and shrinks that distance to an infinitesimal value:

$$y'(x) = \lim_{h \to 0} \frac{y(x+h) - y(x)}{h}; \quad y'_n = \lim_{\Delta x \to 0} \frac{y_{n+1}-y_n}{\Delta x}$$

\noindent On the left is the classical continuous definition, based on a linear interpolation, and on the right is the discrete analog with $y_n = y(x_n)$. Clearly, if we could sample a function infinitely densely, this would evaluate the derivative correctly.

But in computation, the limit can never be perfectly evaluated numerically; rather $\Delta x$ is taken to be sufficiently small so the approximation is accurate to within some acceptable error. For smooth functions, we can bound this error by substituting a truncated Taylor series, which effectively interpolates the function by substituting a degree-$d$ polynomial (with independent variable $\Delta x$) for $y_{n+1}$, e.g.~for $d=2$:

$$\frac{y_{n+1}-y_n}{\Delta x} = \frac{y_n + \Delta x \cdot y'_n + \frac{\Delta x^2}{2!}y''(c) - y_n}{\Delta x} = y'_n + \frac{\Delta x}{2!} y''(c)$$

\noindent for some $c\in [x,x+\Delta x]$. Because $\frac{y''(c)}{2!}$ is simply a constant, we can say this differencing equation is within $O(\Delta x)$ of the true sampled derivative.\footnote{The use of big-O notation is common with Taylor series, but it is technically only allowed in cases where the remainder term can be given an explicit form based on Intermediate Value Theorem. If the function under consideration, $y$, is \textit{not} continuously differentiable $d$ times on the closed interval $[x, x+\Delta x]$, then the error of the numerator Taylor expansion (before division by $\Delta x$) is $o(\Delta x^{d-1})$ rather than $O(\Delta x^d)$~\cite{little-o-big-o}. Little-o, ``$f$ is $o(g)$", means $\lim_{N \to \infty} \frac{f(N)}{g(N)} = \lim_{\Delta x \to 0} \frac{f(\Delta x)}{g(\Delta x)} = 0$, while big-O, ``$f$ is $O(g)$", means $|f(N)| \leq c |g(N)|\ \forall\ N \geq N_0 \Leftrightarrow |f(\Delta x)| \leq c |g(\Delta x)|\ \forall\ \Delta x \leq \Delta x_0$ for some $c > 0$ and $N_0$ or $\Delta x_0\in \mathbb{R}$. For $g$ of the same degree, little-o is (perhaps counterintuitively) a stronger statement than big-O, but little-o of one degree down is a weaker statement. At this place in the flowchart (\autoref{fig:flowchart}), we are considering smooth derivatives, so we assume infinite continuous differentiability, $y \in C^\infty$, and adopt big-O.} This bound can be shrunk by using more points from a ``stencil", like $[-1, 0, 1]$ (i.e.~the previous point, current point, and following point, ``centered" on the current point), and setting up a system of equations to cleverly cancel terms:
\begin{equation}\label{eqn:fd-taylor-expansions}
\begin{aligned}
y'_n \approx&\ c_0 y_{n-1} + c_1 y_n + c_2 y_{n+1}\quad\quad\text{(ansatz)}\\
=&\ c_0 [y_n - \Delta x\cdot y'_n + \frac{(-\Delta x)^2}{2!} y''_n + O(\Delta x^3)] + c_1 y_n \\
&+ c_2 [y_n + \Delta x\cdot y'_n + \frac{\Delta x^2}{2!} y''_n + O(\Delta x^3)]\quad\parbox{4.25cm}{(substitute $d\!=\!3$ expansions\\[-1.5pt]to match 3 unknowns)}
\end{aligned}
\end{equation}

\noindent By grouping terms according to derivative order, we can obtain as many equations as unknown coefficients:
\begin{align*}
y_n:\quad 0 &= c_0 + c_1 + c_2\\
y'_n:\quad 1 &= -\Delta x \cdot c_0 + \Delta x \cdot c_2\\
y''_n:\quad 0 &= \frac{(-\Delta x)^2}{2!}c_0 + \frac{\Delta x^2}{2!}c_2 \rightarrow 0 = c_0 + c_2
\end{align*}

\noindent These can be set up as a linear inverse problem:

$$\begin{bmatrix}
1 & 1 & 1\\
-1 & 0 & 1 \\
1 & 0 & 1
\end{bmatrix}
\begin{bmatrix} c_0 \\ c_1 \\ c_2 \end{bmatrix}
= \begin{bmatrix} 0 \\ \frac{1}{\Delta x} \\ 0 \end{bmatrix}$$

\noindent and solved to yield $c_0 = \frac{-1}{2\Delta x},\ c_1 = 0,\ c_2 = \frac{1}{2\Delta x}$. Plugging these back into \autoref{eqn:fd-taylor-expansions}, the $\Delta x$ in the denominators cancels one of those in the error terms, making the final error $O(\Delta x^2)$.

This methodology can be generalized~\cite{taylor-fd} to higher, $\nu^\text{th}$-order, derivatives and arbitrary stencils of length $S$, $[s_0, ..., s_{S-1}]$, including lopsided, one-sided (``forward" or ``backward" difference), and non-contiguous ones; the linear inverse problem just becomes:
\begin{equation}\label{eqn:fd-generalized}
\begin{bmatrix}
s_0^0 & \cdots & s_{S-1}^0\\
\vdots & \ddots & \vdots \\
s_0^{S-1} & \cdots & s_{S-1}^{S-1}
\end{bmatrix}
\begin{bmatrix} c_0 \\ \vdots \\ c_{S-1} \end{bmatrix}
= \begin{bmatrix} 0 \\ \vdots \\ \frac{\nu!}{\Delta x^\nu} \\ \vdots \\ 0 \end{bmatrix} \leftarrow \parbox{2cm}{nonzero only\\[-1.5pt]at $\nu^\text{th}$ index}
\end{equation}

\noindent The error of the resulting finite difference scheme, $y^\nu_n \approx c_0y_{n+s_0} + ... + c_{S-1}y_{n+s_{S-1}}$, is at most $O(\Delta x^{S-\nu})$, because plugging coefficients back into the scheme's Taylor-expansion causes all terms up to and including the $(S-1)^\text{th}$ order---except the $\nu^\text{th}$---to cancel by construction, leaving only $O(\Delta x^S)$ beyond, which is then divided by ${\Delta x^\nu}$. In the special case of even $\nu$ and an odd-length, centered stencil, the solution causes cancellation of all odd-order terms, including the most dominant error term, $O(\Delta x^S)$, allowing us to achieve higher accuracy with the same number of points. \autoref{ta:fd-schemes} gives formulas with second-order accuracy. For greater accuracy, one can solve \autoref{eqn:fd-generalized} with longer stencils, i.e.~more neighboring points.

\begin{table}[t]
\caption{\label{ta:fd-schemes} Second-order-accurate center-difference formulas along with forward- and backward-difference schemes for endpoints.}
\centering
\begin{tabular}{l}
  \hline $\mathbf{O(\Delta x^2)}$ \textbf{center-difference schemes} \\ \hline $y'_n=[y_{n+1} - y_{n-1}]/2\Delta x$ \\
  $y''_n=[y_{n+1} -2y_n + y_{n-1}]/\Delta x^2$ \\
  $y'''_n=[y_{n+2} -2y_{n+1} +2y_{n-1} - y_{n-2} ]/2\Delta x^3$ \\
  $y^{(4)}_n=[y_{n+2} -4y_{n+1} +6 y_n -4 y_{n-1} +y_{n-2} ]/\Delta x^4$ \\
  \hline $\mathbf{O(\Delta x^2)}$ \textbf{forward- and backward-difference schemes} \\
  \hline $y'_n=[-3y_n + 4y_{n+1} - y_{n+2}]/2\Delta x$ \\
  $y'_n=[3y_n - 4y_{n-1} +  y_{n-2}]/2\Delta x$ \\
  $y''_n=[2y_n-5y_{n+1} +4y_{n+2} -y_{n+3}]/\Delta x^2$ \\
  $y''_n=[2y_n-5y_{n-1} +4y_{n-2} -y_{n-3}]/\Delta x^2$ \\
  \hline
\end{tabular}
\vspace{-3mm}
\end{table}

Finite difference formulas can be compactly represented as sparse matrices acting on a vector of data by simple multiplication. For example, the following two matrices compute the first and second derivatives of data vector $\by$ using the formulas from \autoref{ta:fd-schemes}. Center differencing is used for interior points, but the first and last rows reflect forward and backward difference formulas, respectively.

$$\by' = \frac{1}{2\Delta x} \begin{bmatrix}
  -3 & 4  & -1 & 0 & \cdots & 0\\
  -1 & 0 & 1 & 0 & \cdots & 0\\
  0 & -1 & 0 & 1 & \cdots & 0\\
  \ddots & \ddots & \ddots & \ddots & \ddots & \ddots \\
  0 & \cdots & 0 & -1 & 0 & 1 \\
  0 & \cdots & 0 & 1 & -4 & 3 
\end{bmatrix}
\begin{bmatrix}
  y_1 \\ y_2 \\ \vdots \\ y_n
\end{bmatrix} = {\bf D} \cdot \by$$

$$\by'' = \frac{1}{\Delta x^2} \begin{bmatrix}
  2 & -5 & 4 & -1 & 0 & \cdots\\
  1 & -2 & 1 & 0 & 0 & \cdots\\
  0 & 1 & -2 & 1 & 0 & \cdots\\
  \ddots & \ddots & \ddots & \ddots & \ddots & \ddots\\
  \cdots & 0 & 0 &  1 & -2 & 1 \\
  \cdots & 0 &-1 & 4 & -5 & 2 
\end{bmatrix}
\begin{bmatrix}
  y_1 \\ y_2 \\ \vdots \\ y_n
\end{bmatrix} = {\bf D_2} \cdot \by$$

Finite Differences are \textit{local} approximations, using only a small set of neighboring data points to estimate derivatives. The method is typically constructed to use the minimal number of neighbors to achieve $O(\Delta x^m) = O(N^{-m})$ error, where $m$ is called the ``order"~\cite{kutzbook} and $N$ is the number of samples, because $\Delta x \propto 1/N$ on a fixed domain.

\subsubsection{Recommendations} Finite Difference is simple and flexible. It realizes a relationship between function and derivative values that is at once easy to implement and understand, which has made it the method of choice for several meta-differentiation libraries (\autoref{sec:autodiff-physics}). It works with arbitrary boundary conditions and is fast to run, $O(N m)$, requiring only one major iteration along the data. However, Finite Difference's error bounds only hold for \textit{smooth} functions (capable of being approximated by a Taylor series) and in the absence of noise, which differencing amplifies by $\frac{1}{\Delta x}$ (see \autoref{sec:noise}). This is the same class of functions where Spectral Methods, which achieve better accuracy, can excel, so Finite Difference should only be used if Fourier (\autoref{sec:fourier}) or Chebyshev (\autoref{sec:chebyshev}) do not apply.

\subsection{Spectral Methods}
\label{sec:spectral}

This approach forms a continuous interpolant, $y(x)$, from a set of ``basis functions", $\{\xi_k(x)\},\ k \in \mathbb{N}$, defined over the domain, $x \in [a, b]$. The basis is chosen to span a function space containing the function of interest, meaning $y$ can be represented as a linear combination of the basis functions:
\begin{equation}
\label{eqn:reconstruction-sum}
y(x) = \sum_k c_k \xi_k(x), \quad x \in [a, b]
\end{equation}

\noindent where each $c_k$ is a scalar coefficient. A basis set (e.g.~the Fourier functions in the next section) can have infinitely many members, and we may need them all to fully represent an arbitrary function. But in computation, we always truncate to consider only a finite number of terms, often $N = |\{y_n\}|$, so $y$ has just enough expressive power to uniquely interpolate $N$ samples.

Members of a basis set are linearly independent, i.e.~not mere multiples of one another, but often basis functions are chosen to meet the stronger condition of orthogonality, meaning the inner product between them is defined:

$$\langle \xi_i,\xi_j \rangle = \int\limits_{a}^{b} \xi_j(x) \overline{\xi_i(x)} dx = \begin{cases}0 & i \neq j\\\text{constant} \in \mathbb{R} & i = j\end{cases}$$

\noindent where the overbar denotes a complex conjugate (for complex-valued basis functions).

Taking the inner product of both sides of \autoref{eqn:reconstruction-sum} against the $k^\text{th}$ orthogonal basis function, all terms of the sum but the $k^\text{th}$ disappear, and we can rearrange to find the coefficient, $c_k$:
\begin{equation}
\label{eqn:coefficient-finding-integral}
\frac{\langle \xi_k,y \rangle}{\langle \xi_k,\xi_k \rangle} = \frac{\int\limits_{a}^{b} y(x) \overline{\xi_k(x)} dx}{\|\xi_k(x)\|_2^2} = c_k
\end{equation}

\begin{figure}[!t]\label{fig:vector-projection}
\centering
\begin{tikzpicture}[scale=3, >=stealth]
  \filldraw (0,0) circle (0.5pt);
  \draw[->, thick] (0,0) -- (1,0) node[anchor=west] {$\vec{e}_0$};
  \draw[->, thick] (0,0) -- (0,1) node[anchor=east] {$\vec{e}_1$};
  \draw[->, thick] (0,0) -- (0.9,1.2) node[anchor=west] {$\vec{y}$};
  \draw[->, thick, green!60!black] (0,0) -- (1.35, 0.45) node[anchor=south] {$\vec{\xi}_0$};
  \draw[->, thick, red] (0,0) -- (-0.2213594362, 0.66407830863) node[anchor=north east] {$\vec{\xi}_1$};
  \draw[dashed, gray] (0.9,1.2) -- (-0.27,0.81);
  \draw[dashed, gray] (0.9,1.2) -- (1.17, 0.39);
  \draw[->, dashed] (0,0) -- (1.17,0.39) node[below] {$\mathrm{proj}_{\vec{\xi}_0} \vec{y}$};
  \draw[->, dashed] (0,0) -- (-0.27,0.81) node[left] {$\mathrm{proj}_{\vec{\xi}_1} \vec{y}$};
\end{tikzpicture}
\vspace{-2mm}
\caption{A vector change of coordinates, expressing vector $\vec{y}$ in terms of the spanning basis vectors $(\vec{\xi}_0, \vec{\xi}_1)$ instead of axis-aligned unit vectors $(\vec{e}_0,\vec{e}_1)$, can be accomplished by $\vec{y} = \alpha \vec{\xi}_0 + \beta \vec{\xi}_1$. If $\vec{\xi}_0$ and $\vec{\xi}_1$ are orthogonal, then $\alpha = \frac{\langle \vec{\xi_0},\vec{y}\rangle}{\|\vec{\xi}_0\|_2^2},\ \beta = \frac{\langle \vec{\xi}_1,\vec{y}\rangle}{\|\vec{\xi}_1\|_2^2}$.}
\vspace{-3mm}
\end{figure}
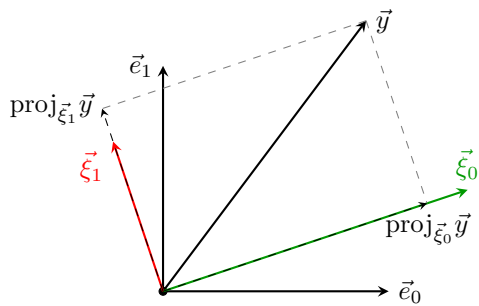

\noindent Note that none of the inner product integrals here are allowed to diverge, which limits $\{\xi_k\}$ and $y$ to the function space $L_2$, the square integrable functions, by definition. This is the only Lebesgue space which is also a Hilbert space (stronger condition), due to the presence of the inner product. This means functions in $L_2$ can be treated completely analogously to vectors, so finding $\{c_k\}$ can be thought of as similar to a vector change of basis, as shown in \autoref{fig:vector-projection}.

If we have only discrete samples, then we need a discrete analog to \autoref{eqn:coefficient-finding-integral}; we can set up a square system to uniquely determine $N$ coefficients:
\begin{equation}
\label{eqn:spectral-linear-inverse}
\begin{bmatrix}
\xi_{0,0} & \xi_{1,0} & \cdots & \xi_{N-1,0}\\
\xi_{0,1} & \xi_{1,1} & \cdots &\xi_{N-1,1} \\
\vdots & \vdots & \ddots & \vdots\\
\xi_{0,N-1} & \xi_{1,N-1} & \cdots & \xi_{N-1,N-1}
\end{bmatrix}
\begin{bmatrix} c_0 \\ c_1 \\ \vdots \\ c_{N-1} \end{bmatrix}
= \begin{bmatrix} y_0 \\ y_1 \\ \vdots \\ y_{N-1} \end{bmatrix}
\end{equation}

\noindent where $\xi_{k,n}$ is the value of the $k^\text{th}$ basis function at the $n^\text{th}$ sample point, $x_n$. This equation specifies that each point $y_n$ should be approximated by a weighted sum of the basis functions at the corresponding $x_n$, so the reconstruction recipe, \autoref{eqn:reconstruction-sum}, still applies. This allows us to represent the function's value over the \textit{entire domain} despite starting with only discrete samples.

The linear inverse problem in \autoref{eqn:spectral-linear-inverse} works for orthogonal and merely linearly independent basis sets alike, because both make the square matrix full rank. This gives the system a unique solution, although its condition number and consequent sensitivity to inexactness in the target, $\{y_n\}$, may vary. In the general case, solving the system requires $O(N^3)$ operations, a relatively expensive computation. However, the Fourier (\autoref{sec:fourier}) and Chebyshev (\autoref{sec:chebyshev}) bases enable much faster solutions on the order of $O(N\log N)$.

Because this method exploits all available information (all samples) to fit all basis functions, it is called \textit{global} and is able to converge to a true underlying smooth function super-algebraically, i.e.~faster than any $O(N^{-m})$ for $m \in \mathbb{R}$~\cite{kutzbook, trefethen2000spectral}. This is sometimes called ``infinite order accuracy" and written $O(N^{-\infty})$\footnote{If $y$ is analytic, convergence is even faster, exponential, i.e.~$O(c^N)$ for $c \in (0,1)$~\cite{trefethen2000spectral}, but in this case one should probably just use \texttt{AutoDiff} (\autoref{sec:autodiff}).}~\cite{trefethen2000spectral}. This implies phenomenal accuracy, but for nonsmooth functions there is no such guarantee, limiting the applicability of spectral interpolants.

\subsubsection{Differentiation with the Fourier Basis}
\label{sec:fourier}

The most ubiquitous~\cite{oppenheim, hamming} basis is the complex exponentials $e^{i\omega x} = \cos(\omega x) + i\sin(\omega x)$, where $i$ is the imaginary unit and $\omega$ is angular frequency. For real-valued functions and signals, each complex exponential is multiplied by a complex coefficient that perfectly compensates for the imaginary component to produce a real, phase-shifted sinusoid. Fitting against and reconstructing with this orthogonal basis in continuous variables gives rise to the Fourier transform, often denoted with a hat $\hat{\ }$, and its inverse:
\begin{equation}
\label{eqn:fourier-transform-pair}
\begin{aligned}
\hat{y}(\omega) &= \int\limits_{-\infty}^{\infty} y(x) e^{-i \omega x} dx = \mathcal{F}\{y(x)\} \\
y(x) &= \frac{1}{2\pi} \int\limits_{-\infty}^{\infty} \hat{y}(\omega) e^{i \omega x} d \omega = \mathcal{F}^{-1}\{\hat{y}(\omega)\}
\end{aligned}
\end{equation}

\noindent If $y$ is absolutely integrable, as it must be unless it contains infinite energy, we can apply integration by parts and find a simple expression for the derivative in the transformed domain:
\begin{equation}
\label{eqn:derivation-of-fourier-deriv}
\begin{aligned}
\mathcal{F}\{\frac{d}{dx} y(x)\} &=\! \int\limits_{-\infty}^{\infty}\! \underbrace{e^{-i \omega x}}_{u} \underbrace{\frac{dy}{dx} dx}_{dv} = \underbrace{y(x) e^{-i \omega x} \Big|_{-\infty}^{\infty}}_{\mathclap{\text{\shortstack{0 for Lebesgue-\\integrable functions}}}} - \int\limits_{-\infty}^{\infty}\! \underbrace{y(x)}_{v} \underbrace{(-i \omega) e^{-i \omega x} dx}_{du}\\[-5pt]
&= i\omega \cdot \hat{y}(\omega)
\end{aligned}
\end{equation}

\noindent That is, we can accomplish differentiation in the Fourier domain with mere multiplication, and we can even apply the above recursively to get a formula for the $\nu^\text{th}$ derivative:
\begin{equation}
\label{eqn:fourier-deriv-relationship}
y^{(\nu)} ={\cal F}^{-1} \{  (i\omega)^\nu {\cal F}\{y\} \}
\end{equation}

\noindent Since we are performing numerical computation on sampled functions, we must use the Discrete Fourier Transform pair:
\begin{equation}
\label{eqn:dft}
\begin{aligned}
\text{DFT: \ \ } Y_k &= \sum_{n=0}^{N-1} y_n e^{-i \frac{2\pi}{N} n k} \\
\text{DFT} ^{-1} \text{: \ \ } y_n &= \frac{1}{N} \sum_{k=0}^{N-1} Y_k e^{i \frac{2\pi}{N} n k}
\end{aligned}
\end{equation}

\noindent where $k$ in this context is often called ``wavenumber" and the coefficient $c_k$ from last section is renamed $Y_k$. The DFT and its inverse can be computed in $O(N \log N)$ by the Fast Fourier Transform (FFT) algorithm~\cite{cooley1965algorithm}, allowing us to avoid inverting the matrix in \autoref{eqn:spectral-linear-inverse}.

There is a discrete analog of Equations \ref{eqn:derivation-of-fourier-deriv} and \ref{eqn:fourier-deriv-relationship} as well~\cite{johnson-fft-notes, spectral-derivatives}, which essentially band-limits the spectral reconstruction to frequencies of lower absolute value:
\begin{equation}\label{eqn:fourier-times-ik}
Y^{(\nu)}_k = \begin{cases} (ik)^\nu \cdot Y_k & k < \frac{N}{2} \\ (i \frac{N}{2})^\nu \cdot Y_k & k = \frac{N}{2} \text{ and } \nu \text{ even} \\ 0 & k = \frac{N}{2} \text{ and } \nu \text{ odd} \\ (i(k - N))^\nu \cdot Y_k & k > \frac{N}{2} \end{cases}
\end{equation}

\begin{figure}[t]\label{fig:gibbs-phenomenon}
  \centering
  \includegraphics[width=0.8\textwidth]{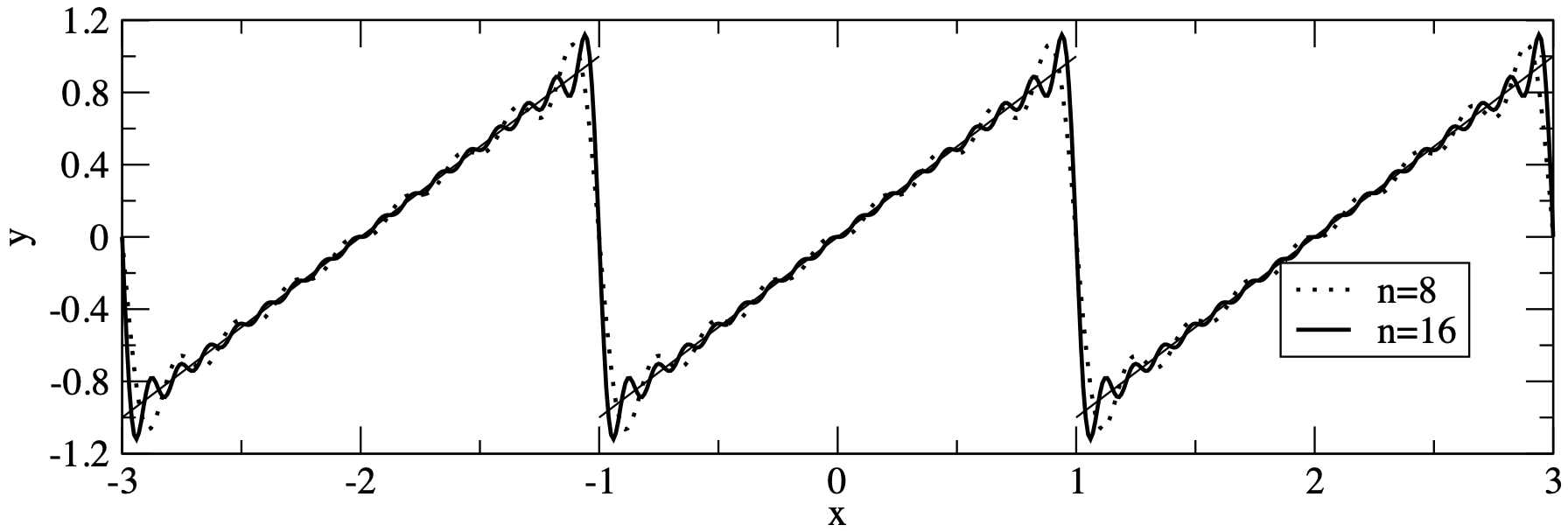}
  \vspace{-3mm}
  \caption{From~\cite{kutzbook}, periodic extension of a ramp function results in discontinuities at the domain endpoints, which cannot be perfectly reconstructed with 8, 16, or any number of smooth sines and cosines, leading to overshoots and wobbles known as Gibbs phenomenon.}
\end{figure}

\noindent This gives rise to the highly accurate and fast solution procedure in \autoref{algo:deriv-via-fft}. But there is a catch: The FFT assumes its input signal is defined on the domain $x\in[0,2\pi)$ and furthermore requires that the signal in question is \textit{periodic} on this interval, because $e^{-i\frac{2\pi}{N}n(k+N)} = e^{-i\frac{2\pi}{N}nk}$ and $e^{i\frac{2\pi}{N}(n+N)k} = e^{i\frac{2\pi}{N}nk}$ in \autoref{eqn:dft}, and therefore $Y_{k+N} = Y_k$ and $y_{n+N} = y_n$. If $y$ is periodic on some other interval of length $T$, we can map $x \in [x_0, x_0 + T)$ to $\theta \in [0, 2\pi)$, find the scaled $\nu^\text{th}$ derivative, and rescale by $(\frac{2\pi}{T})^\nu$ without a problem. But if $y$ is aperiodic, then when the DFT implicitly periodically extends $y$, there will be discontinuities or corners in the extended function, and these are impossible to fit with a finite number of smooth, sinusoidal basis functions, resulting in artifacts known as Gibbs phenomenon, illustrated in \autoref{fig:gibbs-phenomenon}.

\begin{algorithm}[t]
\caption{\textbf{Differentiation via FFT}}
\label{algo:deriv-via-fft}
\begin{algorithmic}[1] 
\State Find the Fourier coefficients $Y_k = \text{FFT}(y_n)$.
\State Multiply by appropriate $(ik)^\nu$, 0, or $(i(k-N))^\nu$ from \autoref{eqn:fourier-times-ik} to make Fourier coefficients of the derivative, $Y^{(\nu)}_k$.
\State Inverse transform to obtain samples of the derivative function $y^{(\nu)}_n = \text{FFT}^{-1}(Y^{(\nu)}_k)$
\end{algorithmic}
\end{algorithm}

\subsubsection{Differentiation with a Polynomial Basis: Chebyshev}
\label{sec:chebyshev}

\begin{figure}[!t]\label{fig:collocation-points}
  \centering
  \includegraphics[width=0.8\textwidth]{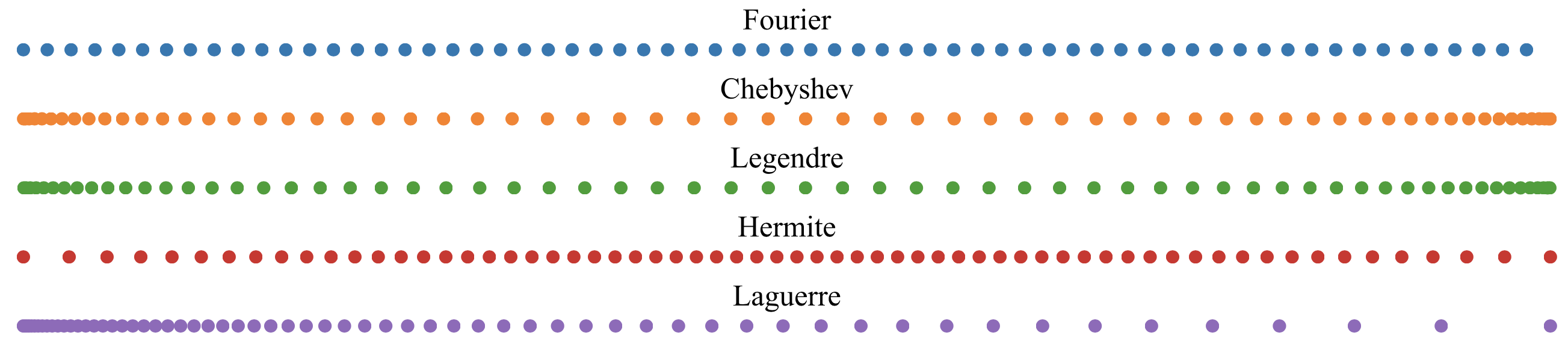}
    \vspace{-2mm}
  \caption{From~\cite{dedalus}, optimal sampling grids across interval domains, also called ``collocation points", for a few common choices of basis. Each basis is defined on its own canonical interval, e.g.~$[0,2\pi)$ for Fourier and $[-1,1]$ for Chebyshev, so other intervals are affine-transformed before basis function fit, which is easily reversible to return to the original domain.}
\end{figure}

\begin{figure}[!t]\label{fig:runge-phenomenon}
  \centering
  \includegraphics[width=0.8\textwidth]{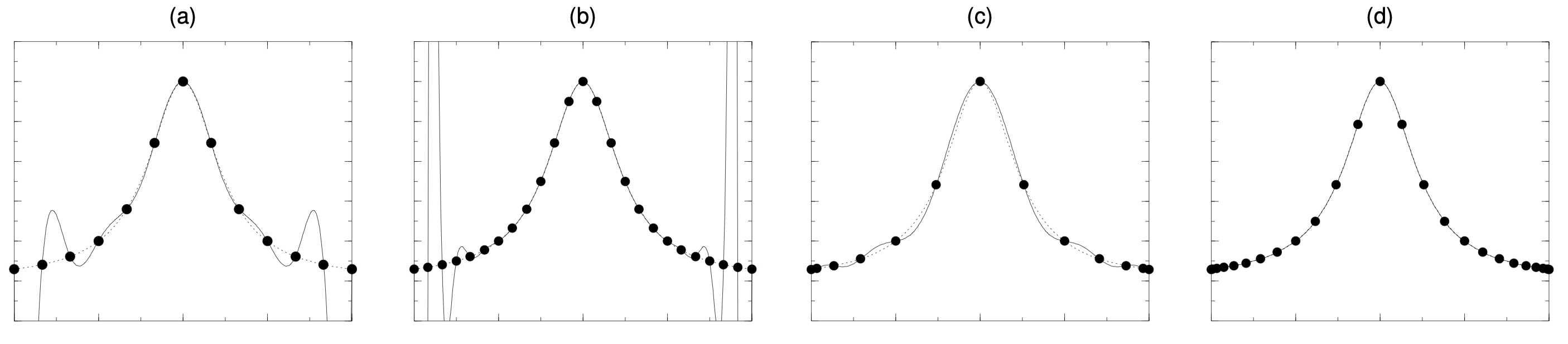}
    \vspace{-2mm}
  \caption{From~\cite{kutzbook}. With equispaced sampling grids, as in (a) and (b), high-order polynomial fits exhibit wobbles at the edges of the domain, known as Runge's phenomenon. This can be mitigated by sampling at nodes better suited to the polynomial basis, such as the Chebyshev-Lobatto nodes in the case of Chebyshev polynomials, shown in (c) and (d).}
\end{figure}

Polynomial bases have the requisite expressive power to represent aperiodic functions and come in many flavors, each with its own properties: Legendre polynomials are natively orthogonal on $[-1, 1]$, while Chebyshev polynomials are orthogonal only with the inclusion of the weight function $1/\sqrt{1-x^2}$ in the inner product integral over $[-1, 1]$. Hermite polynomials are orthogonal on $[-\infty, \infty]$ with weight function $e^{-x^2}$~\cite{dedalus}, but Bernstein polynomials, defined on $[0, 1]$, are not orthogonal~\cite{spectral-derivatives}, just linearly independent. Each has its own optimal sampling grid~\cite{dedalus} (\autoref{fig:collocation-points}) to best avert Runge's phenomenon (\autoref{fig:runge-phenomenon}) by lowering the condition number of the matrix in \autoref{eqn:spectral-linear-inverse}.

\begin{figure}[!t]\label{fig:cylinder}
  \centering
  \includegraphics[width=0.5\textwidth]{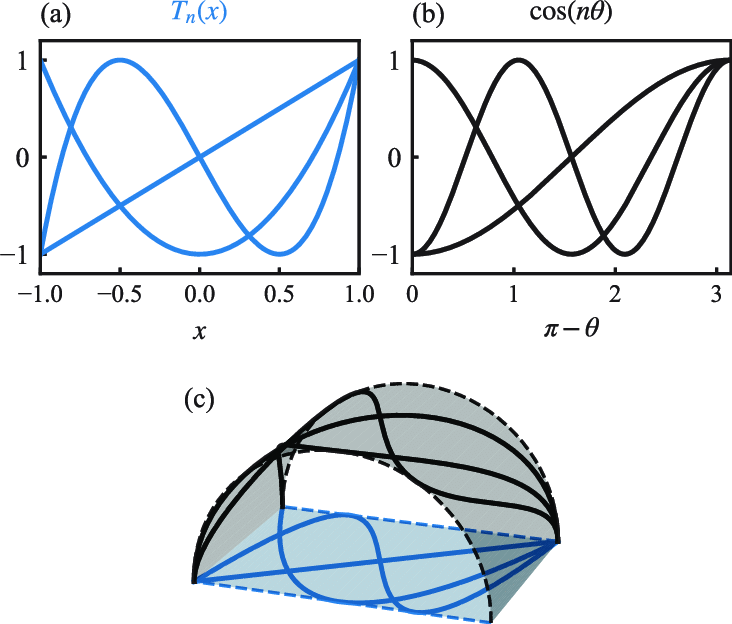}
    \vspace{-2mm}
  \caption{From~\cite{dedalus}, relationship of Chebyshev domain and Fourier domain. Chebyshev polynomials look like the $xy$-plane shadows of cosines wrapped around a cylinder. Notice the cosines are horizontally flipped.}
\end{figure}

Of the polynomial bases, the Chebyshev basis is the most practically useful, because it is identical to a basis of cosines under the change of variables $x = \cos(\theta),\ x \in [-1, 1],\ \theta \in [0, \pi]$ per \autoref{eqn:equivalent}, depicted in \autoref{fig:cylinder}~\cite{trefethen2000spectral, dedalus}.
\begin{equation}\label{eqn:equivalent}
\begin{aligned}
T_k (\cos \theta) = \cos(k\theta) \quad\quad\quad\quad\quad \\
y(x) = \sum_{k=0}^{N-1} a_k T_k(x)\ ;\quad y(\theta) = \sum_{k=0}^{N-1} a_k \cos(k \theta)
\end{aligned}
\end{equation}

\noindent Coefficients of a cosine basis can be found with the Discrete Cosine Transform (DCT),\footnote{Be cautious that, using discrete $\theta_n = \frac{\pi n}{N-1}$, the DCT-I$^{-1}$ says $y_n = \frac{1}{2(N-1)} \big(Y_0 + Y_{N-1}(-1)^n + 2 \sum_{k=1}^{N-2} Y_k\cos(\frac{\pi k n}{N-1})\big)$, but in the Chebyshev expansion $y_n = \sum_{k=0}^{N-1} a_k\cos(\frac{\pi k n}{N-1})$. $Y_0$ is multiplied by $1 = \cos(\frac{\pi \cdot 0 \cdot n}{N-1})$ and $Y_{N-1}$ by $(-1)^n = \cos(\frac{\pi\cdot (N-1)\cdot n}{N-1})$, so the first and last terms can be put inside the sum but are twice as large as the others, and all $Y_k$ have to be scaled to produce $a_k$~\cite{spectral-derivatives}.} based on a modification of the FFT, in $O(N\log N)$, so by reduction to this case, Chebyshev coefficients can be found equally fast. Recall the FFT requires equispaced points in variable $\theta \in [0, 2\pi)$, as does the DCT with $\theta \in [0, \pi]$, where half the interval is dropped because the function is assumed to be even, making information from the latter half redundant. Thus, to take advantage of these algorithms, samples in $x$ need to come from $x_n = \cos(\theta_n) = \cos(\frac{\pi n}{N-1}) \in [-1, 1]$~\cite{trefethen2000spectral}, which are the Chebyshev-Lobatto nodes~\cite{brown-cheb}, corresponding to the $N$ extrema of $T_{N-1}$~\cite{dedalus}.

Once the Chebyshev transform is achieved, then, analogous to the power rule for ordinary representations of polynomials as power series, we can find Chebyshev series coefficients of the derivative in $O(N)$ using a recurrence relation~\cite{brown-cheb}, which can be manipulated to produce \autoref{eqn:dcheb}~\cite{spectral-derivatives}, implemented in \texttt{numpy} as \texttt{chebder}~\cite{numpy-chebder}:
\begin{equation}\label{eqn:dcheb}
\frac{d}{dx} T_k(x) = k \cdot \begin{cases}
2\sum_{\text{odd } j>0}^{k-1} T_j(x) & \text{for even } k \\
-1 + 2\sum_{\text{even } j \geq 0}^{k-1} T_j(x) & \text{for odd } k
\end{cases}
\end{equation}

\begin{algorithm}[!t]
\caption{\textbf{Differentiation with the Chebyshev Basis}}
\label{algo:deriv-via-cheb}
\begin{algorithmic}[1] 
\State Find the Chebyshev coefficients with $a_k \propto Y_k = \text{DCT}(y_n)$
\State Use the series recurrence, \autoref{eqn:dcheb}, to make Chebyshev coefficients of the derivative, $a^{(\nu)}_k$.
\State Scale to make DCT coefficients, $Y^{(\nu)}_k$, and inverse transform to obtain samples of the derivative function $y^{(\nu)}_n = \text{DCT}^{-1}(Y^{(\nu)}_k)$
\end{algorithmic}
\end{algorithm}

\subsubsection{Recommendations}

\begin{figure}[!t]\label{fig:periodic-extensions}
  \centering
  \includegraphics[width=0.8\textwidth]{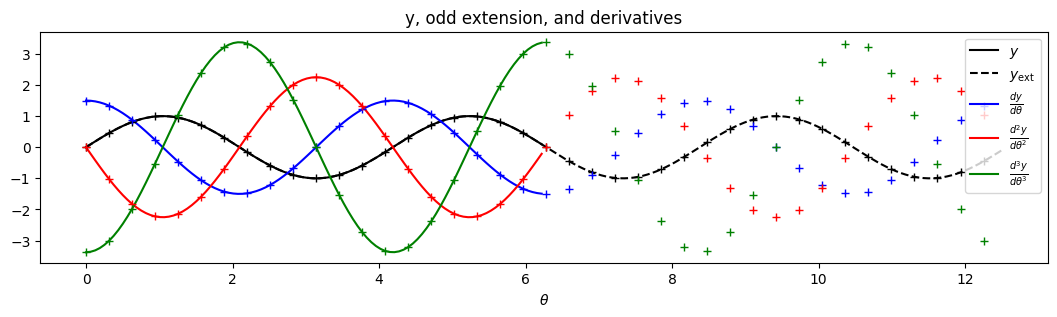}
  \includegraphics[width=0.8\textwidth]{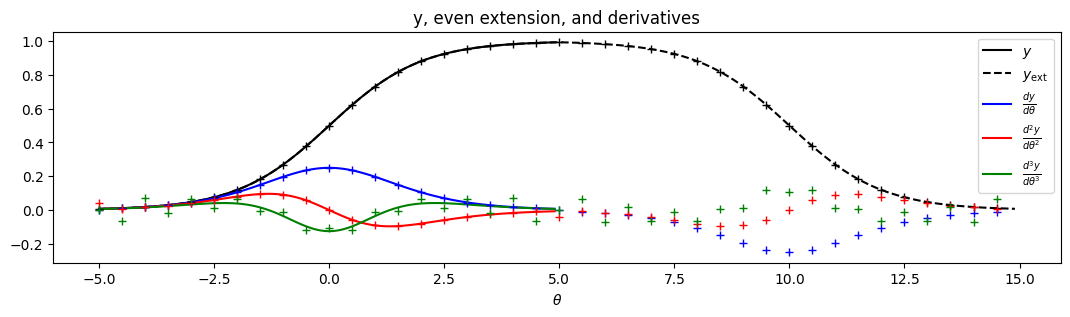}
  \vspace{-3mm}
  \caption{Examples of aperiodic functions that can be periodically extended, $\sin(\frac{3\theta}{2})$ on $\theta \in [0, 2\pi)$ and $\frac{1}{1 + e^{-\theta}}$ on $\theta \in [-5, 5)$. Black + marks indicate function samples, and colorful + marks show the Fourier-based differentiation algorithm's answers for subsequent derivatives. Answers remain essentially perfect for the truly periodic extension but start to drift for higher derivatives of the only-nearly-periodic extension.}
\end{figure}

Fourier spectral differentiation is best used with periodic signals. By including higher frequencies, the fit can get closer at discontinuities and corners, but Gibbs phenomenon will always remain~\cite{oppenheim, trefethen2000spectral}. In some cases this problem can be alleviated or mitigated by concatenating an aperiodic signal with its mirror image (even extension) or negated mirror image (odd extension) to produce a truly periodic signal or something very close to periodic, as in \autoref{fig:periodic-extensions}. The FFT is so fast that transforming a double-length signal is only marginally slower. But this trick does not produce a perfectly periodic extension for general signals.

\autoref{algo:deriv-via-cheb} offers a method with exceptional computational efficiency and accuracy while also being unconstrained by boundary conditions. However, function samples need to come from special locations or an affine transform thereof, and the method should only be used in the absence of noise (see \autoref{sec:noise}, particularly \autoref{fig:cheb-noise}). The Matlab package \texttt{chebfun}~\cite{driscoll2014chebfun} and Python packages \texttt{spectral-derivatives}~\cite{spectral-derivatives} and \texttt{dedalus}~\cite{dedalus} exploit Chebyshev polynomials to allow for accurate and efficient derivative estimates of arbitrary smooth functions.

\subsection{Finite Elements}
\begin{figure}[!t]\label{fig:fe-basis-1}
\centering
\begin{minipage}{0.32\textwidth}
    \centering
  \includegraphics[width=\textwidth]{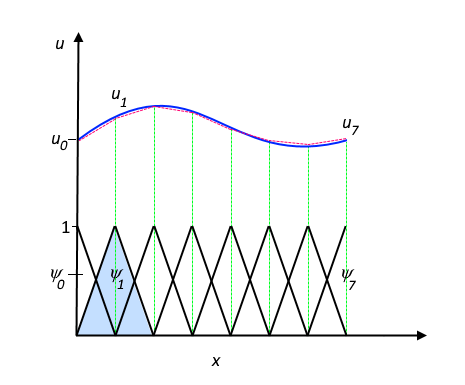}
    \footnotesize (a)
\end{minipage}
\begin{minipage}{0.33\textwidth}
  \centering
  \includegraphics[width=\textwidth]{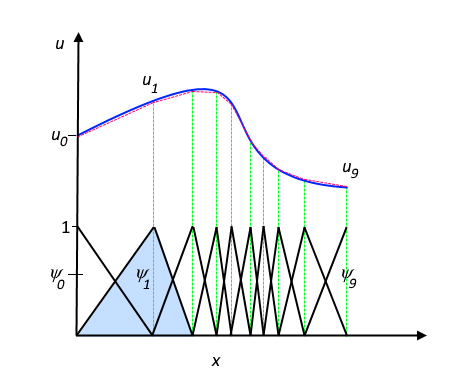}
    \footnotesize (b)
\end{minipage}
\begin{minipage}{0.33\textwidth}
  \centering
  \includegraphics[width=\textwidth]{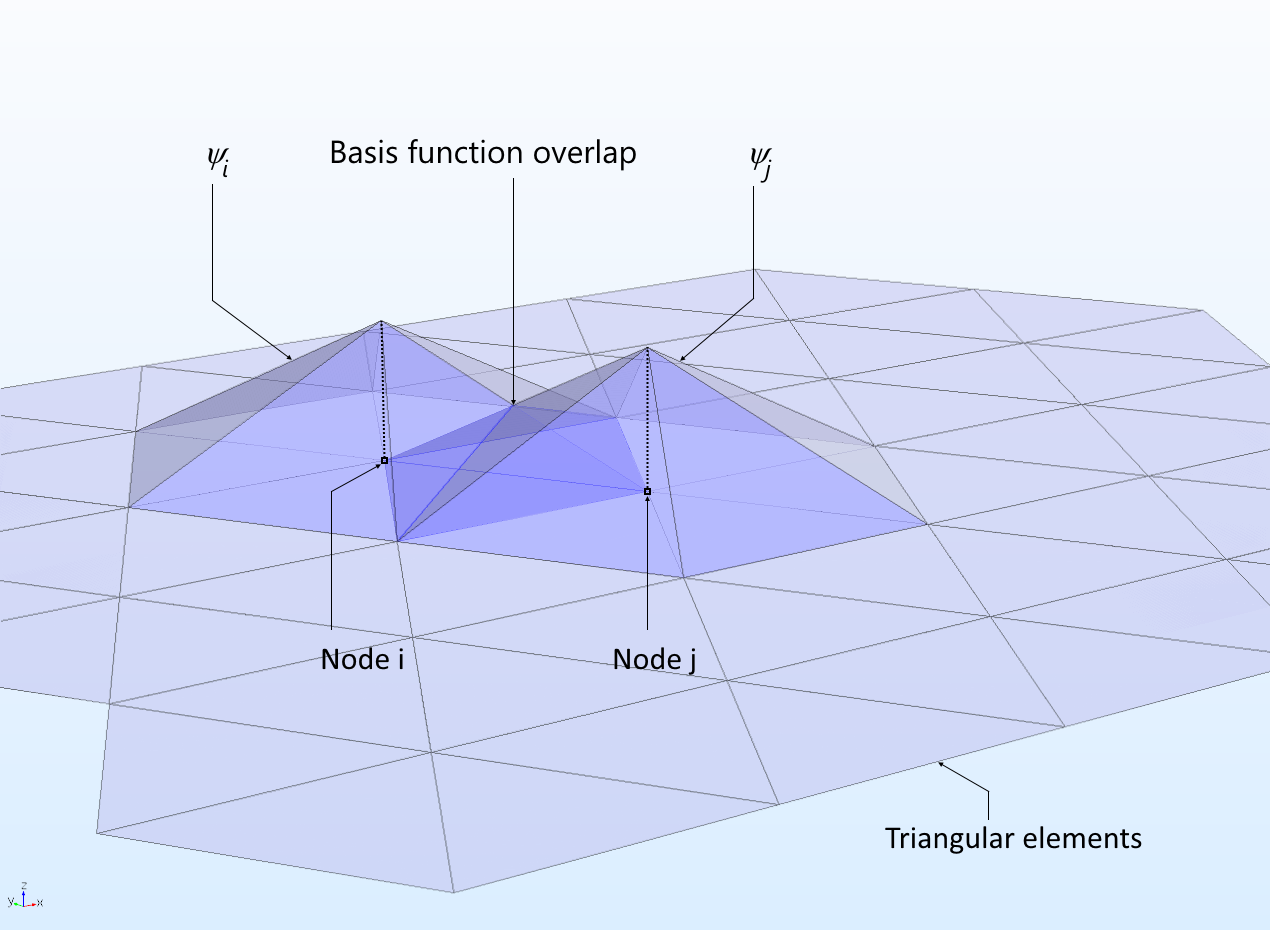}
    \footnotesize (c)
\end{minipage}
\vspace{-2mm}
\caption{From~\cite{comsol}, Finite Element basis functions, $\psi$, with only local support, of type Lagrange polynomial with order 1. (a) Equispaced, (b) on a mesh that captures a steeper section with higher fidelity, (c) in 2D. In (a) and (b) a target function is shown in blue, with the best-fit sum of basis functions shown as a dotted red line.}
\end{figure}

\begin{figure}[!t]\label{fig:fe-basis-2}
  \centering
    \begin{overpic}[width=0.65\textwidth]{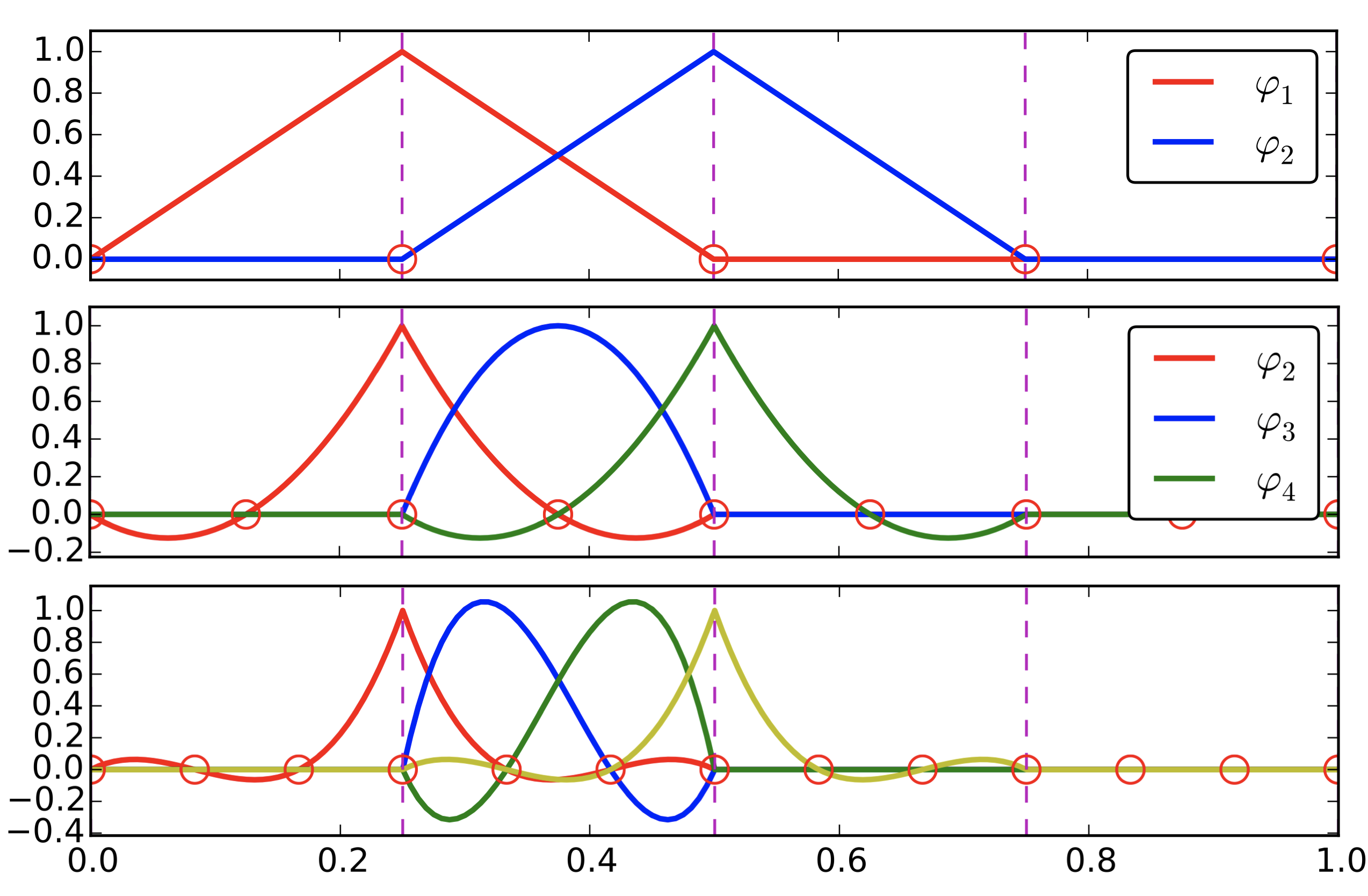}
    \put(100,53){\footnotesize (a)}
    \put(100,33){\footnotesize (b)}
    \put(100,13){\footnotesize (c)}
  \end{overpic}
    \vspace{-2mm}
  \caption{From~\cite{langtangen}, Finite Element basis functions, this time called $\varphi$, of type Lagrange polynomial, (a) linear, (b) quadratic, (c) cubic. Boundaries between elements are shown as pink dotted lines. As the order increases, the number of basis functions associated with each element increases, as does the number of ``nodes" circled in red where all-but-one of these basis functions is zero. Within an element, the basis forms a partition of unity, i.e.~the basis functions sum to 1.}
\end{figure}

\!\!\!\footnote{There is also Finite Volumes, which emphasizes conservation of fluxes across boundaries rather than fit to a basis, introducing integrals that resemble the weak form of Finite Elements (\autoref{sec:weak-form}). In fact, for conservation-law PDEs in divergence form, FVM can be viewed as a Galerkin method (see below) with piecewise constant test functions, though it is conventionally regarded as a distinct scheme.}\label{sec:finite-elements} Instead of projecting a function onto basis functions that span the whole domain, like Spectral Methods (\autoref{sec:spectral}), consider basis functions with local, overlapping domains, as shown in \autoref{fig:fe-basis-1} and \autoref{fig:fe-basis-2}. Basis functions of this kind are capable of representing a certain family of functions, together forming the ``trial space", often symbolized in the literature~\cite{langtangen, larson, dokken} as $V_h$:
$$V_h = \{v : v \in C^d(\Omega),\ v|_{\Omega_e} \in \text{span}(\{\varphi_k\})= \text{the basis space}\}$$

\noindent That is, the space is the set of functions $v$ that satisfy: (1) $d$ times continuous differentiability over domain $\Omega$\footnote{We use $\Omega$ to suggest that the domain and function can be multidimensional, but we primarily stay in 1D for illustration.} and (2) can be built piecewise by stringing together elements $v|_{\Omega_e}$, which are composed from the basis functions and only nonzero over subdomains $\Omega_e$ of size related to a scalar parameter $h$ (hence the subscript). 

Often, the basis space is Lagrange polynomials of low degree, because these are cheap to compute. If we demand only that the resulting $v$ be continuous, then $d=0$, and we build from piecewise linear components, resulting in the faceted appearance of many Finite Elements solutions. Using high-order basis functions is sometimes referred to as ``Spectral Elements", because each element's fit is comprised of many basis functions, reminiscent of Spectral Methods.

\begin{figure}[!t]\label{fig:fe-tessellation}
  \centering
  \includegraphics[width=0.7\textwidth]{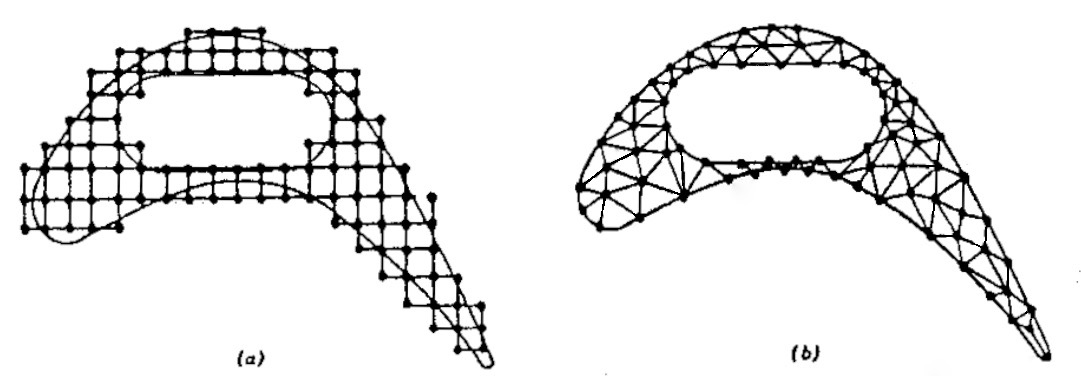}
  \vspace{-2mm}
    \caption{From~\cite{huebner}, a comparison of domain discretization for (a) Finite Difference and (b) Finite Elements.}
\end{figure}

Building from piecewise elements has great benefits for flexible representation of complicated functions on irregular domains, such as shown in \autoref{fig:fe-tessellation}. Whereas ordinary Finite Difference requires a regular grid,\footnote{Technically, the Finite Difference equations in \autoref{sec:finite-difference} can be extended to the irregular sampling case, as described in \autoref{sec:fd-irregular-dt}, but this is not the common case, especially in two or more dimensions where grids could be misaligned.} and Spectral Methods work most naturally on rectangular domains without holes, elements are suited to work in a variety of shapes, including line segments in 1D, triangles and squares in 2D, and tetrahedra and prisms in 3D.

But a piecewise function may not have the expressive power to fully represent the target function, which typically lives in a more general function space, often called $V$, for instance:
$$V = \{v : \|v\|_{L_2(\Omega)} < \infty, \|\nabla v\|_{L_2(\Omega)} < \infty, v(\partial\Omega) = \text{boundary conditions}\}$$

\noindent Here the $L_2$ norm of $v$, i.e.~the square root of the inner product of $v$ with itself, which is an integral over the domain, has to converge, as does the norm of $v$'s first derivative\footnote{Differentiation in the \textit{weak} sense, smuggled in by notation; see \autoref{sec:weak-form}. Spaces of finite-norm functions with finite-norm weak derivatives are called Sobolev spaces, of which $V$ is the subset also meeting the boundary conditions, and $V_h$ is the finite-dimensional subset spanned by the chosen basis functions. This allows corners but excludes jump discontinuities, which have Dirac delta weak derivatives (see \autoref{fig:weak-derivatives}), because although $\int_{-\infty}^\infty \delta(x) dx = 1$, $\|\delta\|_{L_2}^2 = \int_{-\infty}^\infty (\delta(x))^2 dx = \infty$.} ($\nabla$ if there are multiple dimensions), and $\partial\Omega$ is the domain boundary.

\begin{figure}[!t]\label{fig:fe-fn-projection}
\centering
\begin{minipage}{0.49\textwidth}
    \centering
  \includegraphics[width=\textwidth]{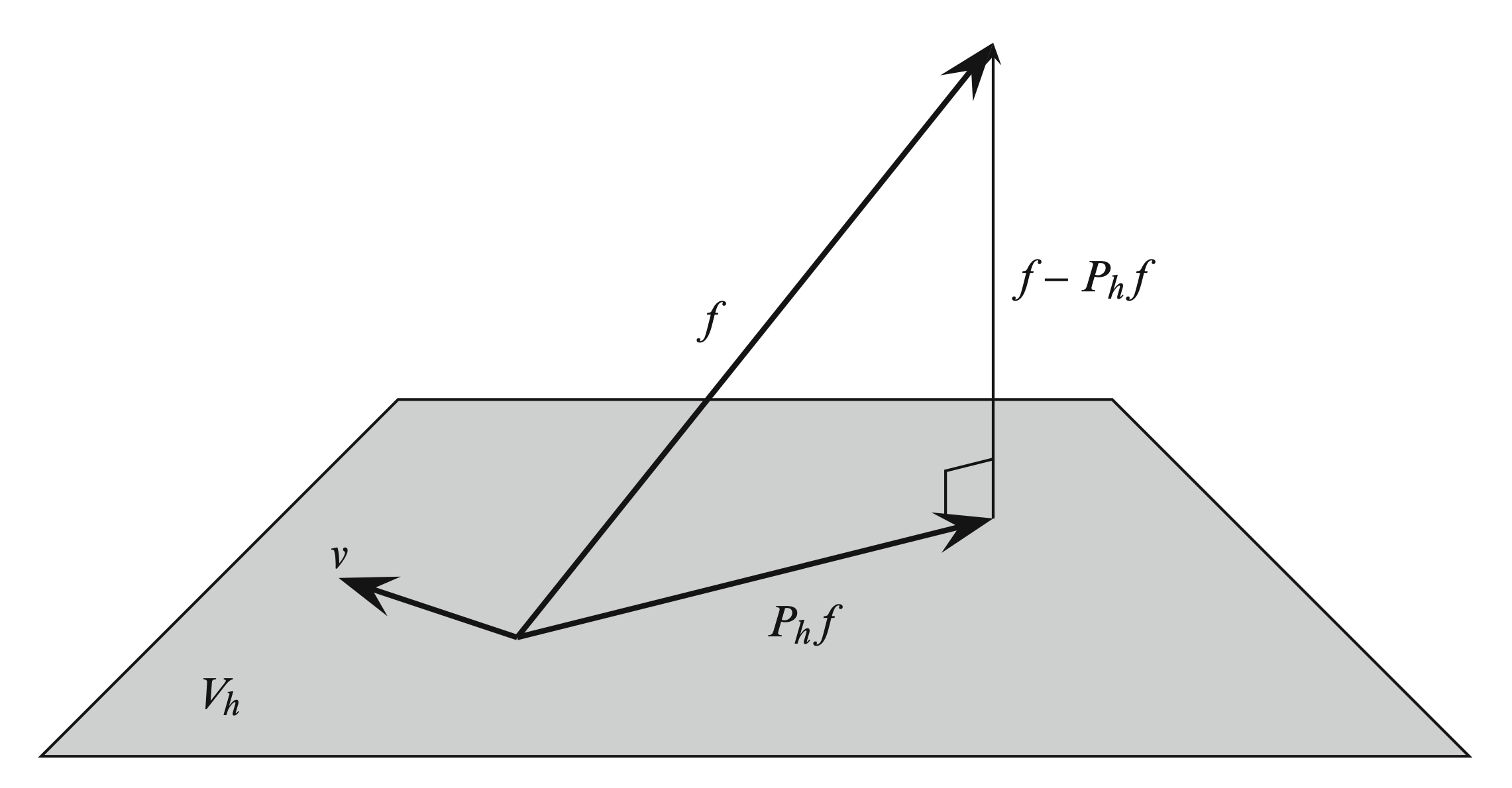}
    \footnotesize (a)
\end{minipage}
\begin{minipage}{0.49\textwidth}
  \centering
  \includegraphics[width=\textwidth]{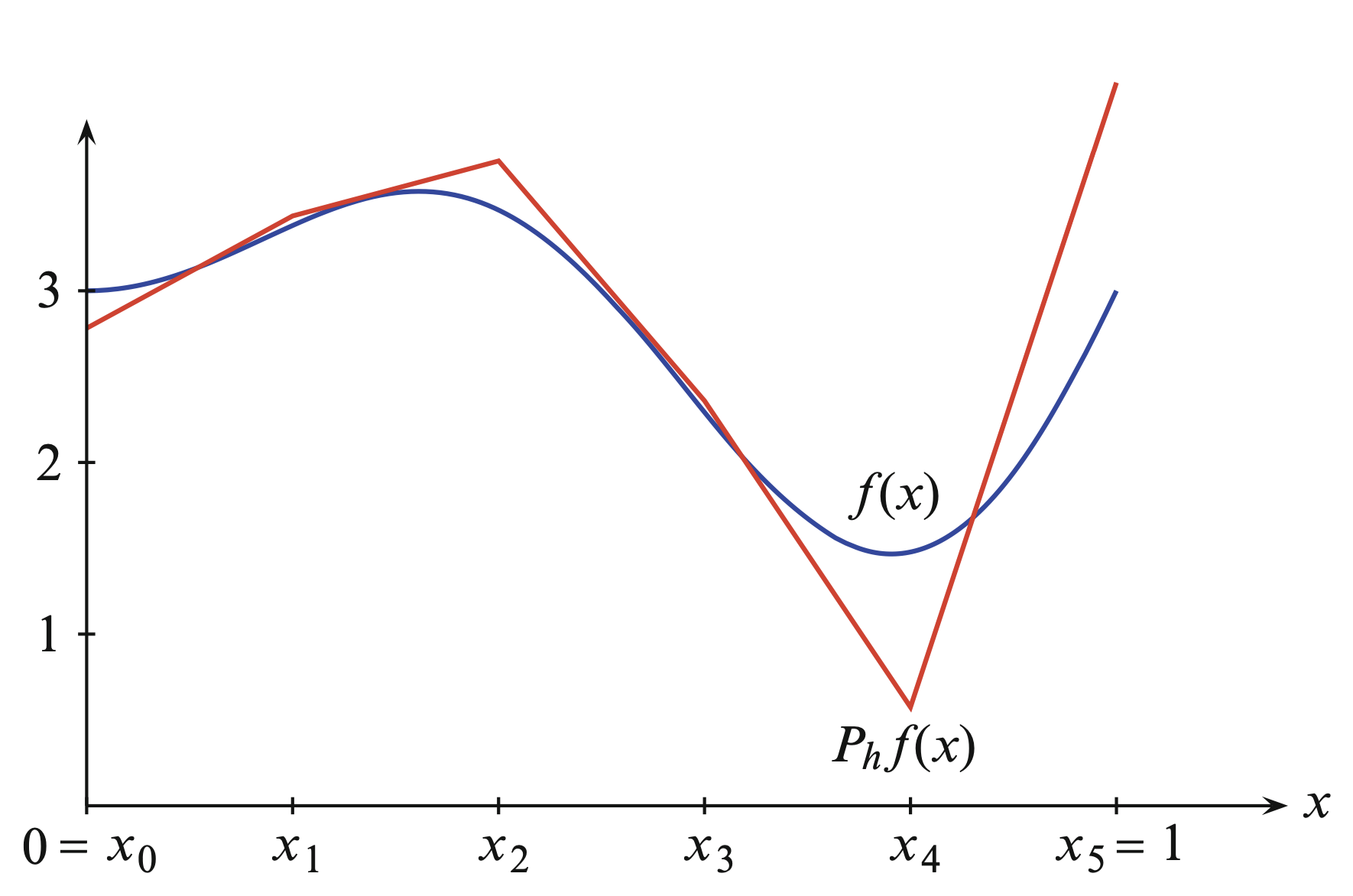}
    \footnotesize (b)
\end{minipage}
\vspace{-1mm}
\caption{From~\cite{larson}. (a) Notional illustration of a function, here called $f$, being $L_2$-projected onto the space $V_h$. (b) $f(x) = 2x\sin(2\pi x) + 3$ and its $L_2$-projection $P_h f(x)$ using the mesh $[0, 0.2, 0.4, 0.6, 0.8, 1]$. The fit does not perfectly fall on the function at the mesh points, the boundaries of elements, because the fit is not linear interpolation; rather it is best fit in an $L_2$-distance, inner-product induced-norm sense. See \autoref{eqn:L2-induced-norm}.}
\end{figure}

Our first task is to represent $y \in V$ with $y_h \in V_h \subset V$, as illustrated in \autoref{fig:fe-fn-projection}. We can naturally measure distance between these objects with the induced $L_2$ norm:
\begin{equation}\label{eqn:L2-induced-norm}\|y - y_h\|_{L_2(\Omega)}^2 = \int\limits_\Omega (y - y_h) \overline{(y-y_h)} d\Omega
\end{equation}

\noindent where the overbar denotes complex conjugate. The minimizer of this equation is provably\footnote{See the short proof of theorem 1.1 in~\cite{larson}.} $y_h = P_h y$, the projection of $y$ onto the subspace $V_h$. But we cannot minimize this equation directly, because continuous functions are infinite-dimensional: Fitting at all continuum points is impossible to compute.

Instead we can use what is called a Galerkin method, one among a group of methods known as the \textit{weighted residuals approach}~\cite{huebner}, to turn the problem into a discrete set of constraints and solve these to obtain $y_h$. The trick is to notice that $y-y_h$ (called $f-P_hf$ in \autoref{fig:fe-fn-projection}(a)) is orthogonal to all $v \in V_h$, which in this context are called ``test functions":
$$\langle y - y_h, v\rangle = 0
\Leftrightarrow \int\limits_\Omega (y - y_h)\overline{v}\ d\Omega = 0$$

Since $y_h$ is composed of basis functions, we can write it as $y_h = \sum_{j} c_j \varphi_j(x)$, where the $\{c_j\}$ are coefficients. The total number of basis functions is related to the mesh size, $N$, and the order of the basis, which can drive up the number of nodes per element. For linear basis functions as in \autoref{fig:fe-basis-1}(a), elements only have nodes at their boundaries, so $j \in \{0,..., N-1\}$.

A Galerkin projection uses the same family of functions for the basis and test sets. We have actually already encountered something which can be thought of as Galerkin projection, the Fourier transform, or more specifically the Fourier series~\cite{oppenheim, fourier-series}, which maps from a continuous function on a periodic interval, say $[-\pi, \pi]$, to discrete coefficients, using complex exponentials:
\begin{equation}\label{eqn:Fourier-is-Galerkin}
\int\limits_{-\pi}^\pi (y(x) - y_h(x))e^{-ikx} dx = 0
\end{equation}

The reconstruction sum $y_h(x) = \sum_{j=0}^{N-1} c_j e^{ijx}$ leaves us with $N$ unknown coefficients, so we need $N$ equations to lock them down. These we can obtain through \autoref{eqn:Fourier-is-Galerkin} by projecting against test functions with different $k$. Plugging in the ansatz for $y_h$, we find:
\begin{equation}\label{eqn:Fourier-is-Galerkin-algebraic}
\begin{aligned}
&\rightarrow \int\limits_{-\pi}^\pi y(x)e^{-ikx} dx = \int\limits_{-\pi}^\pi \Big(\sum_{j=0}^{N-1} c_j e^{ijx} \Big)e^{-ikx} dx\quad\text{for } k = 0, ..., N-1\\
&\rightarrow \underbrace{\langle y(x), e^{ikx} \rangle}_{\text{\small$b_k$}} = \sum_{j=0}^{N-1} c_j \underbrace{\int\limits_{-\pi}^\pi e^{ijx}e^{-ikx} dx}_{\text{\footnotesize$\begin{cases} 2\pi & j = k \\ 0 & j \neq k \end{cases}$}} = 2\pi c_k\quad\text{for } k = 0, ..., N-1\\[-15pt]
&\rightarrow \mathbb{I}\bc = \frac{1}{2\pi}\mathbf{b}
\end{aligned}
\end{equation}

\noindent The coefficients in \textbf{c} are of course the Fourier coefficients.

A linear inverse problem of this kind is typical in Finite Elements. In general, the matrix multiplying the coefficients, often called the ``mass matrix", $\bM$, is not identity, but thankfully using local, piecewise basis functions causes the entries of $\bM$ to be sparse, because most pairs of basis functions do not overlap, so their products integrate to 0. This allows the system to be solved more efficiently than by naive $O(N^3)$ matrix inversion.

The approximation error of fitting the function $y$ with piecewise $y_h$ is $O(h^\alpha)$, where $h$ is related to the side-length of a typical element~\cite{comsol, larson}. This is an algebraic bound, similar to that of Finite Difference, and indeed Finite Elements can be thought of as an alternative form of Finite Difference that interpolates over patches of the domain~\cite{kutzbook, langtangen}. The parameter $\alpha$ can be difficult or impossible to estimate for all but simple problems \textit{a priori}, so software solvers often use their solutions to calculate \textit{a posteriori} error estimates and choose how to make refinements~\cite{comsol, langtangen}.

\subsubsection{Differentiation with Finite Elements}
\label{sec:weak-form}

We can now compute a piecewise basis-function representation of a function we wish to differentiate, but the approximation involves inner product integrals. These are computed via quadrature rules~\cite{larson}, which plug particular ``quadrature points" into the integrand. To achieve an accurate projection, $y_h$, we need freedom to evaluate initial condition $y_0$ and boundary terms wherever the quadrature rule asks, as well as to place elements where the function needs them. Handed a fixed array of measurements instead, any interpolation presupposes a (probably) suboptimal choice of $y_h$, and the projection degenerates to fitting a piecewise polynomial through raw data, which is spline interpolation (\autoref{sec:splinesmoothing}) in the absence of noise. Thus Finite Element methods are best for working with analytic governing equations, not data samples.

\begin{figure}[!t]\label{fig:weak-derivatives}
\centering
\begin{tikzpicture}
\begin{axis}[at={(0,0)}, anchor=origin,
    width=5cm, height=4cm, axis lines=middle,
    xlabel={$x$}, title={ReLU},
    xtick={-1, 0, 1}, ytick={0, 1},
    ymin=-0.5, ymax=2, xmin=-2, xmax=2,
    axis line style={->}, ticks=both]
\addplot[thick, blue, domain=-2:0] {0};
\addplot[thick, blue, domain=0:2] {x};
\end{axis}
\begin{axis}[at={(4.5cm,0)}, anchor=origin,
    width=5cm, height=4cm, axis lines=middle,
    xlabel={$x$}, title={Heaviside},
    xtick={-1, 0, 1}, ytick={0, 1},
    ymin=-0.5, ymax=2, xmin=-2, xmax=2,
    axis line style={->}, ticks=both]
\addplot[thick, red, domain=-2:0] {0};
\addplot[thick, red, domain=0:2] {1};
\draw[red, fill=white] (axis cs:0,0) circle[radius=1.2pt];
\draw[red, fill=red] (axis cs:0,1) circle[radius=1pt];
\end{axis}
\begin{axis}[at={(9cm,0)}, anchor=origin,
    width=5cm, height=4cm, axis lines=middle,
    xlabel={$x$}, title={Dirac Delta (symbolic)},
    xtick={-1, 0, 1}, ytick=\empty,
    ymin=-0.5, ymax=2, xmin=-2, xmax=2,
    axis line style={->},
    ticks=both]
\draw[violet, thick, ->] (axis cs:0,0) -- (axis cs:0,2);
\end{axis}
\end{tikzpicture}
\vspace{-1mm}
\caption{Examples of functions where the notion of a weak derivative is useful. The Rectified Linear Unit (ReLU) has a sharp corner at $x=0$, making it not strictly differentiable. But intuitively its derivative should look something like a Heaviside function, which in turn has a weak derivative given by the Dirac Delta. This notion is made rigorous by the machinery of ``generalized" functions (also called ``distributions")~\cite{weak-derivs}, members of which can be written as a \emph{functional} which maps from function to value via an integral against a test function, essentially identical to the form in Galerkin projection.}
\end{figure}
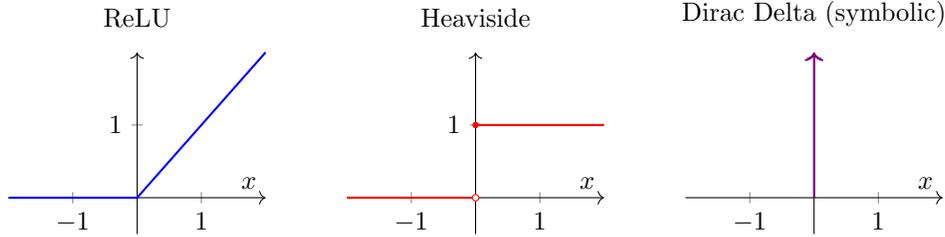

Since $y_h$ is intentionally piecewise, it will have corners and kinks at element boundaries, making it not strictly differentiable at those locations, like the functions in \autoref{fig:weak-derivatives}. So instead of directly taking derivatives to solve a governing equation, the solution is found in the ``weak" (integral) form, which only requires that the inner products of the two sides of the equation with test functions be equal, thereby allowing discontinuities in the first derivative~\cite{comsol}.

This setup matches the Galerkin method naturally, and it comes with the added benefit, thanks to an integration by parts trick\footnote{The same manipulation is famously used in Calculus of Variations to transfer a differential operator from the variation to the Lagrangian while spawning a boundary term, thereby enabling a factorization that yields the Euler-Lagrange equation. In higher dimension, this trick may be called Green's First Identity~\cite{langtangen}, a special case of Divergence Theorem.} (demonstrated below), that the solution need only be (weakly) differentiable one fewer times in space than in the original PDE, thereby enabling us to describe it with a lower-order basis, which means fewer nodes (see \autoref{fig:fe-basis-2}) and a smaller system of equations, especially helpful for saving computation in higher dimensions.\footnote{The weak formulation may not be strong, but it is powerful.}

It is instructive to see an example. In the following we use subscript $x$ and $t$ as a shorthand to mean partial derivatives. The weak solution for the 1D wave equation, often named $u$ as in \autoref{sec:wave-example} rather than $y$, can be formed as:
$$u_{tt} = u_{xx},\quad x \in [0, L]$$
$$u(x,0) = u_0(x),\quad u_t(x, 0) = 0,\quad u(0,t) = u(L,t) = 0$$
$$\rightarrow \int\limits_0^L u_{tt} v\ dx = \int\limits_0^L u_{xx} v\ dx\quad \text{{\hsize=4cm \tight[\raggedright]{\small inner product both sides against a test function}}}$$ 

Using $u_{xx} dx = d(\frac{d}{dx}u)$, we can integrate the right-hand side by parts:
$$\int\limits_0^L v\ du_x = u_xv\Big|_0^L - \int\limits_0^L u_x dv$$

\noindent If we require the test functions $v$ to be 0 at the edges of the domain, then the first term disappears. We can additionally use $dv = \frac{d}{dx}v\ dx$ to obtain:
$$\int\limits_0^L u_{tt} v\ dx = -\int\limits_0^L u_x v_x\ dx$$

Notice the equation now only depends on $u_x$ rather than $u_{xx}$. Imposing a time-space separation, $u(x, t) = \sum_{j=0}^{N-1} c_j(t) \varphi_j(x)$, which further equals $\sum_{j=1}^{N-2} c_j(t) \varphi_j(x)$, because the zero boundary conditions dictate the values of the coefficients of the basis functions at those points, $c_0, c_{N-1} = 0$, leaving us with only $N-2$ unknowns. Choosing $v$ from the basis set ($\varphi$), and now using $\dot{\ }$ to represent a time derivative and $'$ to represent a spatial derivative, we get:
\begin{align*}
&\int\limits_0^L \sum_{j=1}^{N-2} \ddot{c}_j(t)\varphi_j(x)\varphi_k(x)\ dx = -\int\limits_0^L \sum_{j=1}^{N-2} c_j(t) \varphi'_j(x) \varphi'_k(x)\ dx\quad\text{for } k = 1, ..., N-2\\
\rightarrow &\sum_{j=1}^{N-2} \ddot{c}_j(t) \underbrace{\int\limits_0^L \varphi_j(x)\varphi_k(x)\ dx}_{M_{jk}} = -\sum_{j=1}^{N-2} c_j(t) \underbrace{\int\limits_0^L \varphi'_j(x) \varphi'_k(x)\ dx}_{K_{jk}}\quad\text{for } k = 1, ..., N-2\\
&\rightarrow \bM\ddot{\bc}_{1:N-2}(t) + \bK\bc_{1:N-2}(t) = 0
\end{align*}

The $(N-2)\!\times\!(N-2)$ entries of $\bM$ and $\bK$ do not change over time. We have transformed the PDE into a system of ODEs. Initial coefficients \textbf{c} can be found from initial conditions and then forward-solved with a method of lines scheme and steppers like Runge-Kutta. If the PDE instead describes a steady state with unknown $u$, the derivation culminates in an algebraic equation rather than ODEs, similar to the Fourier series example in \autoref{eqn:Fourier-is-Galerkin-algebraic}. For a further example, see the derivation for the 2D Poisson equation in any of~\cite{langtangen, larson, dokken, olver, uiuc-galerkin}, which is also instructive for generalizing to more spatial dimensions.

The overall solution process of Finite Elements is quite involved, from mesh generation, to boundary condition enforcement, to linear inverse problem (or even nonlinear problem~\cite{langtangen, larson}) setup, to numerical solution methods (e.g.~Newton-Raphson), to error estimation and refinement. Thankfully, most of these steps can be automated, which is exactly what is done by commercial packages like COMSOL and Ansys and open-source ones like FEniCSx~\cite{fenicsx}.

\subsubsection{Recommendations} Finite Elements is not well suited to differentiating a typical array of numbers representing a signal, but it is the most versatile method for solving governing PDEs on arbitrary domains with arbitrary boundary conditions, which are common in engineering, especially in two and three dimensions. Due to the use of piecewise functions and the weak form, Finite Elements are also capable of representing true discontinuities like shockwaves,\footnote{It is possible to approximate a discontinuity with Spectral Methods, as demonstrated by Dedalus~\cite{dedalus}, but the solution will require a high number of modes to represent and is still technically smooth, with rounded corners and a steep connecting segment.} which can be present in the solutions of hyperbolic PDEs~\cite{olver}. Due to its complexity, Finite Elements is best done via established software. Commercial packages require a license, while the free, open-source FEniCSx requires formulating a PDE in the weak form~\cite{dokken}.

\section{Noise}
\label{sec:noise}

In practice, most real-world data is corrupted due to imperfect sensing (``measurement noise"), unknown forcing in the dynamics (``process noise"), or statistical or numerical artifacts (such as sampling noise or quantization noise). We model this by considering observations to be the additive result of some unknown true signal and noise, $\eta$:
\begin{equation}\label{eqn:noisy-signal}
y = y_\text{true} + \eta
\end{equation}

Noise can be deadly for sample-based differentiation methods, because small deviations in sample values can lead to large deviations of slope, especially when considering a constant measurement noise level with shrinking $\Delta x$. This effect is only magnified for calculations of curvature and higher-order derivatives, $\frac{d^\nu y}{dx^\nu}$. For example, in Finite Differences $(y_{n+1} - y_n)/\Delta x = (y_{\text{true},n+1} - y_{\text{true},n})/\Delta x + (\eta_{n+1} - \eta_n)/\Delta x$, where the last term is error, which compounds recursively to $O(\Delta x^{-\nu})$ for the $\nu^\text{th}$ derivative.

However, if we have reason to believe the data is generated by an underlying process that evolves \textit{smoothly} along our independent variable of interest (whether time or space), and we sample sufficiently densely to capture these variations, then meaningful changes in the signal can be represented by a sum of relatively low-wavenumber Fourier basis functions. This is a central observation in Signal Processing: For a naturally-occurring signal, ``energy" in the sense of Parseval's Relation~\cite{oppenheim}, \autoref{eqn:parseval}, must be finite, and consequently its spectrum must fall off toward zero as frequency increases~\cite{hamming}.
\begin{equation}\label{eqn:parseval}
\int\limits_{-\infty}^{\infty} |y(x)|^2 dx = \frac{1}{2\pi}\!\!\!\int\limits_{-\infty}^{\infty} |\hat{y}(\omega)|^2 d\omega
\end{equation}

\begin{figure}[!t]\label{fig:svd}
  \centering
  \includegraphics[width=0.99\textwidth]{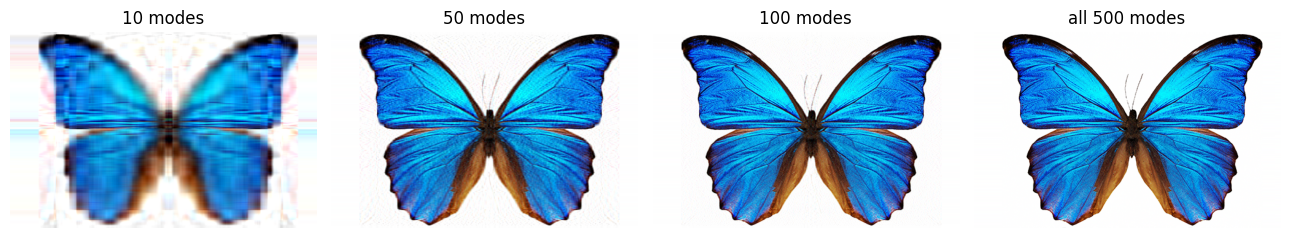}
  \caption{We can stack red, blue, and green channels of an image, perform SVD on the 2D matrix, and reshape $U[:,:r] \cdot \Sigma[:r] \cdot V^\mathrm{T}[:r]$ to produce reconstructions with $r$ modes. In this case, the first 10 modes alone contain more than half the signal energy. Source image from Smithsonian.}
\end{figure}

\noindent It is worth noting that energy also empirically clusters in lower modes when signals are decomposed using other functions, e.g.~polynomials (\autoref{sec:chebyshev}) or neural nets~\cite{driscoll2014chebfun, spectral-bias-nn}. This is even more so the case when a basis is tailored to represent a particular signal, like Singular Value Decomposition (SVD)~\cite{ddse}, as demonstrated in \autoref{fig:svd}.

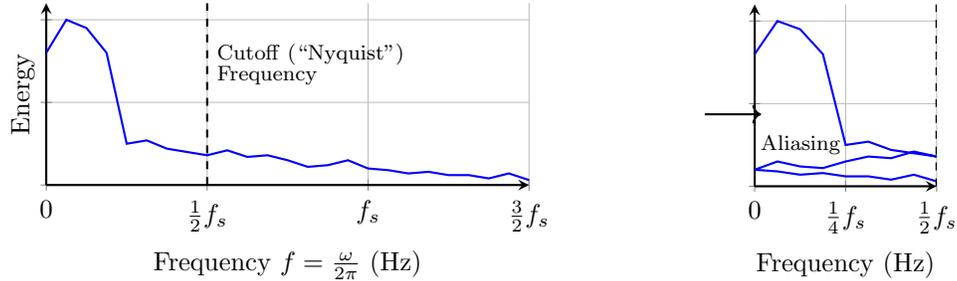
\begin{figure}[!t]\label{fig:noise-spectrum}
\centering
\begin{minipage}{0.6\textwidth}
\centering
\begin{tikzpicture}
  \begin{axis}[width=8cm, height=4cm,
    xlabel={Frequency $f = \frac{\omega}{2\pi}$ (Hz)}, ylabel={Energy}, ylabel style={yshift=-3mm},
    xmin=0, xmax=1.5, xtick={0, 0.5, 1, 1.5}, xticklabels={$0$, $\frac{1}{2}f_s$, $f_s$, $\frac{3}{2}f_s$},
    ymin=0, ymax=1.1, yticklabels={},
    domain=0:1.5,
    samples=200,thick,
    grid=both, axis lines=left]
  \addplot[blue, thick] coordinates {(0.00, 0.8) (0.0625, 1) (0.125, 0.95) (0.1875, 0.8) (0.25, 0.25)
    (0.3125, 0.27) (0.375, 0.22) (0.4375, 0.2) (0.5, 0.18) (0.5625, 0.21) (0.625, 0.17) (0.6875, 0.18)
    (0.75, 0.15) (0.8125, 0.11) (0.875, 0.12) (0.9375, 0.15) (1, 0.1) (1.0625, 0.09) (1.125, 0.07)
    (1.1875, 0.08) (1.25, 0.06) (1.3125, 0.06) (1.375, 0.04) (1.4375, 0.07) (1.5, 0.03)};
  \addplot[dashed, black] coordinates {(0.5,0) (0.5,1.1)};
  \node[black, right, text width=2.5cm] at (axis cs:0.5,0.75) {\footnotesize \tight[\raggedright]{Cutoff (``Nyquist") Frequency}};
  \end{axis}
\end{tikzpicture}
\end{minipage}
\begin{minipage}{0.3\textwidth}
\centering
\begin{tikzpicture}
  \begin{axis}[width=4cm, height=4cm,
    xlabel={Frequency (Hz)},
    xmin=0, xmax=0.5, xtick={0, 0.25, 0.5}, xticklabels={$0$, $\frac{1}{4}f_s$, $\frac{1}{2}f_s$},
    ymin=0, ymax=1.1, yticklabels={},
    domain=0:0.5,
    samples=200,thick,
    grid=both, axis lines=left]
  \addplot[blue, thick] coordinates {(0, 0.8) (0.0625, 1) (0.125, 0.95) (0.1875, 0.8) (0.25, 0.25)
    (0.3125, 0.27) (0.375, 0.22) (0.4375, 0.2) (0.5, 0.18) (0.4375, 0.21) (0.375, 0.17) (0.3125, 0.18)
    (0.25, 0.15) (0.1875, 0.11) (0.125, 0.12) (0.0625, 0.15) (0, 0.1) (0.0625, 0.09) (0.125, 0.07)
    (0.1875, 0.08) (0.25, 0.06) (0.3125, 0.06) (0.375, 0.04) (0.4375, 0.07) (0.5, 0.03)};
  \addplot[dashed, black] coordinates {(0.5,0) (0.5,1.1)};
  \node[black, right] at (axis cs:-0.01,0.25) {\footnotesize Aliasing};
  \end{axis}
\end{tikzpicture}
\end{minipage}
\begin{tikzpicture}[overlay, remember picture]
  \draw[->, thick] (-5.75, 0.5) -- (-4.75, 0.5);
\end{tikzpicture}
\vspace{-2mm}
\caption{Reproduced from~\cite{spectral-derivatives}, typical energy spectrum in the frequency (i.e.~Fourier) domain of a noisy signal, before and after sampling. $f_s$ is a sampling rate of our choosing. Equispaced samples of a sinusoid with frequency just over the Nyquist limit, $\frac{1}{2}f_s + \epsilon$, can be fit equally well by a single sinusoid of frequency just under the Nyquist limit, $\frac{1}{2}f_s - \epsilon$, a phenomenon known as ``aliasing" due to this 1:1 relationship. This causes a folding pattern in the spectrum, here shown explicitly, although the final spectrum of a sampled signal is the sum of these folds up to infinite frequency. The infinite sum must converge, so the signal contains ever-decreasing power at higher frequencies. With increasing $f_s$, the folding noise spectrum's energy is wider-spread, causing less overlap with the signal.}
\end{figure}

Noise must obey this same energy decay, but its Fourier spectrum tends to have less bias toward low frequencies, because noise contains sharp variations that can only be modeled with higher frequencies. Alternative bases tend not to separate high-frequency noise as consistently across the domain~\cite{spectral-derivatives}. \autoref{fig:noise-spectrum} sketches the spectrum of a noisy signal, which can be thought of as the sum of signal and noise spectra, by consequence of \autoref{eqn:noisy-signal} and the linearity of the Fourier transform. It also depicts the phenomenon of aliasing, where a sampled higher frequency looks like and adds to a lower one, because it takes $>\!2$ samples per cycle to truly resolve a sinusoid~\cite{oppenheim}, while a spectral reconstruction is composed of lower-magnitude wavenumbers (see \autoref{sec:fourier}).

The \textit{band-separability} of signal from noise illustrated in \autoref{fig:noise-spectrum} is the spiritual core of noise reduction techniques: At bottom, all accomplish some kind of low-pass filtering. The following sections explore myriad such approaches and present a strategy for choosing algorithm hyperparameters based on the frequency content of measured data.

\subsection{Recommendations} \label{sec:noise-recs}

\begin{figure}[!t]\label{fig:cheb-noise}
  \centering
  \includegraphics[width=0.99\textwidth]{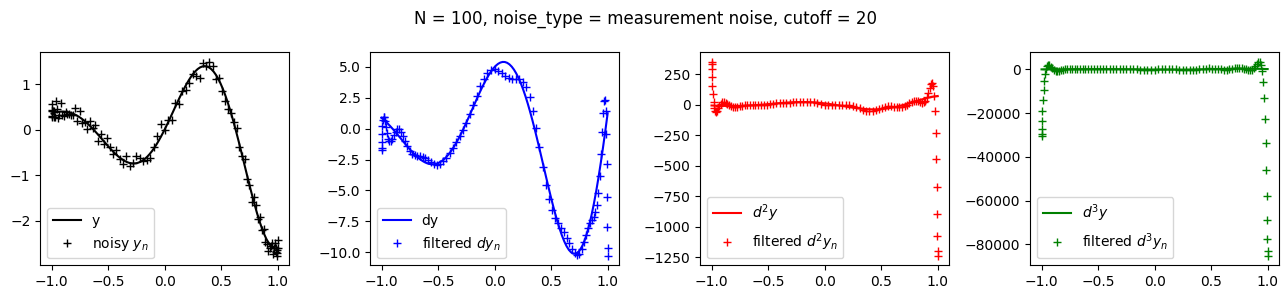}
  \vspace{-2mm}
  \caption{From~\cite{spectral-derivatives}, the function $y(x) = e^x\sin(5x)$ is sampled with noise $\sim \mathcal{N}(\mu=0, \sigma=0.1)$ at 100 points and fit with 100 Chebyshev polynomials, only the first 20 of which are kept after a filtering step. The fit is characteristically prone to overfit near domain edges where the basis functions are steeper and have more wobbles, leading to dramatic errors in derivative estimates. Raising the cutoff makes this phenomenon worse, while lowering the cutoff leads to poor fits by degrading the method's expressive power. This phenomenon is not easily sidestepped by padding, because cosine-spaced samples cannot be padded without breaking cosine-spacing and $O(N \log N)$ runtime.}
\end{figure}

When a signal is known to be periodic, applying an ideal low-pass filter, i.e.~\autoref{algo:deriv-via-fft} with higher modes zeroed out before FFT$^{-1}$, separates signal from high-frequency noise efficiently and optimally---although some noise spectrum will still unavoidably overlap the signal's spectrum. Beware, one should \textit{not} attempt to use the Chebyshev basis in the presence of noise, because polynomials have nonuniform oscillations and therefore do not filter noise equally well throughout the domain. Even worse, noise-induced fit inaccuracies \textit{propagate and compound} from higher modes to lower ones during differentiation due to Chebyshev mode coupling\footnote{By contrast with the aliasing of sinusoids, the best fit to samples from a single higher-degree polynomial requires many different lower-degree polynomials.}~\cite{brown-cheb}, resulting in systematic error toward the domain edges, as depicted in \autoref{fig:cheb-noise}. Likewise, the eigenmodes of a signal, ideal for parsimonious representation~\cite{ddse}, are generally not best for separating noise due to their nonuniform frequency content across the domain. So in practice, with noisy, aperiodic data, one of the many methods discussed in the following two sections is likely to provide the most accurate results.

\section{Differentiating Noisy Data Using Prior Knowledge}
\label{sec:noisy-with-knowledge}

In some cases, additional information may be available, such as a model of the system that produced data, a model of noise, or relevant synchronized data streams. Incorporating this knowledge can help distinguish signal from noise. The most dominant framework for this situation is Kalman filtering/smoothing,\footnote{The filtering problem is to make estimates online, in one direction, whereas smoothing assumes a whole data series is available offline, enabling a backwards pass. See \autoref{sec:rauch-tung-striebel}.} which models the true signal (often along with one or more of its derivatives) as a hidden state and calculates estimates by weighted sum of noisy measurements and dynamics-based predictions.\footnote{This same balance of model against measurements is ubiquitous across methods in the context of uncertainty, although it can take very different forms. For example, Physics-Informed Neural Networks (\autoref{sec:autodiff}) incorporate a system model by fitting noisy data according to a residual loss while attempting to satisfy an analytic governing-equation regularizer. Structurally, PINNs encode a continuous field over space-time and do not natively quantify uncertainty, whereas a Kalman filter's online state-space estimates encode particular quantities as probability distributions at discrete sample times. For the purposes of differentiation, one could read derivatives off a trained PINN, but the dependency runs the other way, because regularizer terms are built \textit{from} derivatives.}

The standard Kalman filter handles only \textit{linear} dynamics, uses a mean squared error (MSE) penalty metric, and is the Maximum \textit{A Posteriori} (MAP) estimator (\autoref{sec:MAP}) only for \textit{zero-mean Gaussian white noise} disturbances. This is fairly limiting and reflects the original formulation's place in history,\footnote{The Kalman filter was meant to generalize the Wiener filter for non-stationary processes~\cite{mit-wiener-kalman}, although Wiener can explicitly handle colored noise (correlated in time) while Kalman requires modification with shaping filters and augmented state to turn such noise white~\cite{crassidis_junkins}.} so we contextualize by surveying a variety of generalizations, including the Luenberger observer~\cite{luenberger} (\autoref{sec:kalman-filter}); robust variants designed to handle more exotic noise distributions, outliers, sparsity, and constraints on the state~\cite{aravkin} (\autoref{sec:robust-estimation}); $H$ filtering~\cite{robust-control} (Hardy space); and a modification for use with nonlinear systems~\cite{julier-uhlmann} (\autoref{sec:UKF}).

Common symbols used in this framework differ from the notation used so far, so we give them in \autoref{ta:kalman-symbols}.

\begin{table}[!t]
\caption{\label{ta:kalman-symbols} Common symbols used with Kalman Filters.}\vspace{-2mm}
\centering
\begin{tabular}{lp{11.4cm}}
  \hline $\bx$ & true signal (or state), not directly observable\\
  $\hbx$ & estimated signal (state)\\
  $\be$ & state estimation error, $\bx - \hbx$\\
  $\by$ & noisy measurement\\
  $\bu$ & known control input or synchronized data stream\\
  $\bA$ & sometimes also called $\mathbf{F}$ or $\boldsymbol\Phi$, linear state evolution matrix or Jacobian of nonlinear state evolution function $\bff$\\
  $\bB$ & linear control matrix or relation of known states to hidden states\\
  $\bC$ & sometimes also called $\mathbf{H}$ or $\boldsymbol\Phi$ (confusingly), linear measurement matrix or Jacobian of nonlinear measurement function $\bh$\\
  $\bw$ & process noise\\
  $\bv$ & measurement noise\\
  $\bQ$ & process noise covariance\\
  $\bR$ & measurement noise covariance\\
  $\bP$ & covariance of state estimate error\\
  $\bK$ & the Kalman gain matrix, filter parameters\\
  $p$ & a probability density function\\
  $\mathcal{N}$ & the normal distribution, univariate or multivariate, parameterized by scalar or vector mean, $\mu$ or $\bmu$, and standard deviation or covariance, $\sigma$ or $\bSigma$\\ \hline
\end{tabular}
\vspace{-4mm}
\end{table}

\subsection{The Classic Kalman Filter}\label{sec:kalman-filter}

The Kalman filter solves the minimum mean squared error optimization problem, for a discrete update from index $n-1$ to $n$:
\begin{equation}\label{eqn:kalman-minimization}
\begin{aligned}
\underset{\bK}{\text{minimize}} \ \ &\mathbb{E}[\|\bx_n - \hbx_n\|_2^2]\\
\text{s.t.}\quad \text{(i)}\ \ &\bx_n = \bA\bx_{n-1} + \bB\bu_n + \bw_n,\ \bw_n \sim \mathcal{N}(0, \bQ)\\
\text{(ii)}\ \ &\by_n = \bC\bx_n + \bv_n,\ \bv_n \sim \mathcal{N}(0, \bR)\\
\text{(iii)}\ \ &\hbx_n = \bA\hbx_{n-1} + \bB\bu_n + \bK(\by_n - \hby_n)\\
\text{(iv)}\ \ &\hby_n = \bC(\bA\hbx_{n-1} + \bB\bu_n)
\end{aligned}
\end{equation}

\noindent That is, assuming linear dynamics and measurement models, perturbed by noise drawn from zero-mean Gaussian probability distributions ((i) and (ii)), form the state estimate by weighted sum of (1) a prediction from the dynamics and (2) the consequent measurement error ((iii) and (iv)). Then find the weight matrix projecting from measurement space to state space that minimizes the Mean Squared Error (MSE), the expected value of the squared $\ell_2$ norm of the state estimate error.

\begin{figure}[!t]\label{fig:covariance-ellipse}
\centering
\begin{tikzpicture}[scale=1.5, >=stealth]
  \draw[->] (-1.5, 0) -- (1.5, 0) node[anchor=west] {};
  \draw[->] (0, -1) -- (0, 1) node[anchor=south] {};
  \draw[thick, dashed] (0,0) ellipse[x radius=1.4, y radius=0.7, rotate=30];
  \draw[->, thick, blue] (0,0) -- ({1.4*cos(30)}, {1.4*sin(30)}) node[blue, anchor=west] {$\sigma_0$};
  \draw[->] (0,0) -- ({cos(30)}, {sin(30)}) node[below] {$\vec{v}_0$};
  \draw[->] (0,0) -- ({cos(120)}, {sin(120)}) node[left] {$\vec{v}_1$};
  \draw[->, thick, red] (0,0) -- ({0.7*cos(120)}, {0.7*sin(120)}) node[red, anchor=north east, xshift=2mm, yshift=-2mm]{$\sigma_1$};
\end{tikzpicture}
\vspace{-1mm}
\caption{Visualization of one standard deviation of a 2D Gaussian distribution with mean $\bmu = 0$ and covariance $\bSigma = \bV\bLambda \bV^\mathrm{T}$.}
\end{figure}
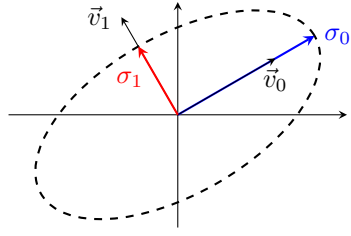

Noise terms in the model are characterized by \textit{vector} means (0 for both) and covariance \textit{matrices}, $\bQ = \text{cov}(\bw, \bw)$ and $\bR = \text{cov}(\bv, \bv)$, because states and measurements may have multiple variables. Covariance of a random vector $\rbchi$ with itself is given as:
$$\text{cov}(\rbchi,\rbchi) = \mathbb{E}[(\rbchi - \mathbb{E}[\rbchi])(\rbchi - \mathbb{E}[\rbchi])^\mathrm{T}] = \mathbb{E}[\rbchi\rbchi^\mathrm{T}] - \mathbb{E}[\rbchi]\mathbb{E}[\rbchi]^\mathrm{T}$$

\noindent where the last step uses the fact $\mathbb{E}[\rbchi]$ is a constant, so the multiplicative cross terms both have the same value as the last term, but with opposing sign: $\mathbb{E}[\mathbb{E}[\rbchi] \rbchi^\mathrm{T}] = \mathbb{E}[\rbchi \mathbb{E}[\rbchi]^\mathrm{T}] = \mathbb{E}[\rbchi]\mathbb{E}[\rbchi]^\mathrm{T}$. Because these covariance matrices are formed as products of vectors with their transposes, they are always symmetric and positive semidefinite and can be eigendecomposed:
$$\bV\bLambda \bV^\mathrm{T} = \begin{bmatrix} \vert & \vert &\\v_0 & v_1 & \cdots\\ \vert & \vert & \end{bmatrix} \begin{bmatrix}\sigma_0^2 & 0 & \cdots\\0 & \sigma^2_1 & \\ \vdots & & \ddots\end{bmatrix} \begin{bmatrix} \rotvert\ v_0\ \rotvert \\ \rotvert\ v_1\ \rotvert\\ \vdots \end{bmatrix} $$

\noindent where $\bV$ is an orthonormal basis aligned with principal directions of variation, and $\bLambda$ is filled with variances, $\sigma \geq 0$. This provides a neat way to visualize the standard deviation: Multiply points on a unit sphere by the matrix $\bV\bLambda^{1/2}$, which performs scaling followed by rotation, to produce points on an ellipsoid, as in \autoref{fig:covariance-ellipse}.

An illuminating reformulation of \autoref{eqn:kalman-minimization}, purely in terms of state estimate error, $\be = \bx - \hbx$, can be found by substituting (ii) and (iv) into (iii) and subtracting from (i):
\begin{equation}\label{eqn:kalman-min-error}
\begin{aligned}
\hspace{4em}\underset{\bK}{\text{minimize}}\ \ &\mathbb{E}[\|\be_n\|_2^2]\\
\text{s.t.}\quad&\be_n = (\mathbb{I} - \bK\bC)\bA\be_{n-1} + \underbrace{(\mathbb{I} - \bK\bC)\bw_n - \bK\bv_n}_\text{still Gaussian}
\end{aligned}
\end{equation}

\noindent The constraint in \autoref{eqn:kalman-min-error} gives the error dynamics. For error to decay, it is necessary that the closed-loop discrete-time map, $(\mathbb{I} - \bK\bC)\bA$, has eigenvalues inside the unit circle (magnitude $\!<\!1$). This is achievable by pole placement when the discrete-time pair $(\bA,\bC\bA)$ is observable~\cite{xu-chen}, giving rise to a Luenberger observer~\cite{luenberger}, of which the Kalman filter is a special case that additionally accounts for noise by solving the minimization problem.

The error distribution is itself also a zero-mean Gaussian, because error is the sum of a deterministic part, that tends toward zero for a stable observer, and scaled, zero-mean Gaussians. Its mean is evaluated as $\mathbb{E}[(\mathbb{I} - \bK\bC)(\bA\be_{n-1} + \bw_n) -\bK\bv_n] = (\mathbb{I} - \bK\bC)(\bA\mathbb{E}[\be_{n-1}] + \mathbb{E}[\bw_n]) -\bK\mathbb{E}[\bv_n]= 0$, and its covariance is:

\begin{equation}\label{eqn:cov-P}
\begin{aligned}
\text{cov}(\be_n, \be_n)&=\mathbb{E}[((\mathbb{I} - \bK\bC)(\bA\be_{n-1} + \bw_n) -\bK\bv_n)((\mathbb{I} - \bK\bC)(\bA\be_{n-1} + \bw_n)-\bK\bv_n)^\mathrm{T}]\\
&=(\mathbb{I} - \bK\bC)(\bA\underbrace{\mathbb{E}[\be_{n-1}\be_{n-1}^\mathrm{T}]}_{\mathclap{\text{cov}(\be_{n-1}, \be_{n-1})}} \bA^\mathrm{T} + \underbrace{\mathbb{E}[\bw_n\bw_n^\mathrm{T}]}_\bQ)(\mathbb{I} - \bK\bC)^\mathrm{T} + \bK\underbrace{\mathbb{E}[\bv_n\bv_n^\mathrm{T}]}_\bR \bK^\mathrm{T}
\end{aligned}
\end{equation}

\noindent where the cross terms are zero because $\be$, $\bv$, and $\bw$ are assumed to be independent, so $\mathbb{E}[\be\bw^\mathrm{T}]\!=\!\mathbb{E}[\bv\be^\mathrm{T}]\!=\!\mathbb{E}[\bv\bw^\mathrm{T}]\!=\!0$. We give this covariance the name $\bP$ and note that it is recursively defined like $\be$.

The minimization problem can be further reformulated in terms of this $\bP$:
\begin{equation}\label{eqn:kalman-min-tr-P}
\begin{aligned}
\hspace{4em}\underset{\bK}{\text{minimize}}\ \ &\text{Tr}[\bP_n]\\
\text{s.t.}\quad&\bP_n = \mathbb{E}[\be_n\be_n^\mathrm{T}]\\
&\be_n = (\mathbb{I} - \bK\bC)\bA\be_{n-1} + (\mathbb{I} - \bK\bC)\bw_n - \bK\bv_n
\end{aligned}
\end{equation}

\noindent because $\mathbb{E}[\|\be\|_2^2] = \mathbb{E}[\text{Tr}[\be\be^\mathrm{T}]]$, and the order of trace and expected value can be flipped, since they are linear operators. This form is solvable with calculus~\cite{lacey-kf} to determine the optimal $\bK$. First manipulate \autoref{eqn:cov-P}:
\begin{align}
\bP_n &=(\mathbb{I} - \bK\bC)\underbrace{(\bA\bP_{n-1}\bA^\mathrm{T} + \bQ)}_{\text{name this }\bP_{n|n-1}}(\mathbb{I} - \bK\bC)^\mathrm{T} + \bK\bR \bK^\mathrm{T} \label{eqn:Pn}\\
&=\bP_{n|n-1} - \bK\bC\bP_{n|n-1} - \bP_{n|n-1}(\bK\bC)^\mathrm{T} + \bK\bC\bP_{n|n-1}(
\bK\bC)^\mathrm{T} + \bK\bR\bK^\mathrm{T} \notag\\ \notag\\
\rightarrow &\text{Tr}[\bP_n] = \text{Tr}[\bP_{n|n-1}] \underbrace{- 2\text{Tr}[\bK\bC\bP_{n|n-1}]}_{\mathclap{\text{because Tr[$\bM$]=Tr[$\bM^\mathrm{T}$] and $\bP$=$\bP^\mathrm{T}$}}} + \text{Tr}[\bK(\bC\bP_{n|n-1}\bC^\mathrm{T} + \bR)\bK^\mathrm{T}] \notag
\end{align}

\noindent where subscript $n|n-1$ indicates an \textit{a priori} estimate at step $n$ based only on information up to step $n-1$, and subscript $n$ (often written $n|n$) indicates the \textit{a posteriori} estimate after the step resolves and new information is factored in. Using the facts that for matrices $\bM$ and $\bW$, $(\bW\bM)^\mathrm{T} = \bM^\mathrm{T}\bW^\mathrm{T}$; $\nabla_\bM \text{Tr}[\bM\bW] = \bW^\mathrm{T}$; and $\nabla_\bM \text{Tr}[\bM\bW\bM^\mathrm{T}] = \bM\bW^\mathrm{T} + \bM\bW$,  which further equals $2\bM\bW$ when $\bW$ is symmetric~\cite{duchi-trace-derivative}:
\begin{align*}
\rightarrow & \nabla_\bK \text{Tr}[\bP_n] = -2(\bC\bP_{n|n-1})^\mathrm{T} + 2\bK(\bC\bP_{n|n-1}\bC^\mathrm{T} + \bR) = 0\\
\rightarrow & \bK_\text{opt} = \bP_{n|n-1}\bC^\mathrm{T} (\bC\bP_{n|n-1}\bC^\mathrm{T} + \bR)^{-1}\quad  \text{\small$\blacksquare$}
\end{align*}

Discrete Kalman filtering, then, is the process of making this closed-form calculation of $\bK$ repeatedly at each step and applying it to the system (constraints of \autoref{eqn:kalman-minimization} or \autoref{eqn:kalman-min-tr-P}) to form updated estimates of the state and error covariance. The framework is flexible enough to accept different system, measurement, and noise covariance matrices at each step, in which case these get subscript $n$, omitted above for brevity. When all matrices are constant, then because the calculation of $\bK$ does not rely on the data itself, the filter gain approaches a steady state that accomplishes the same calculation as a Wiener filter~\cite{mit-wiener-kalman, wunsch} and can be calculated offline and hard-coded to avoid repeated matrix inversion.

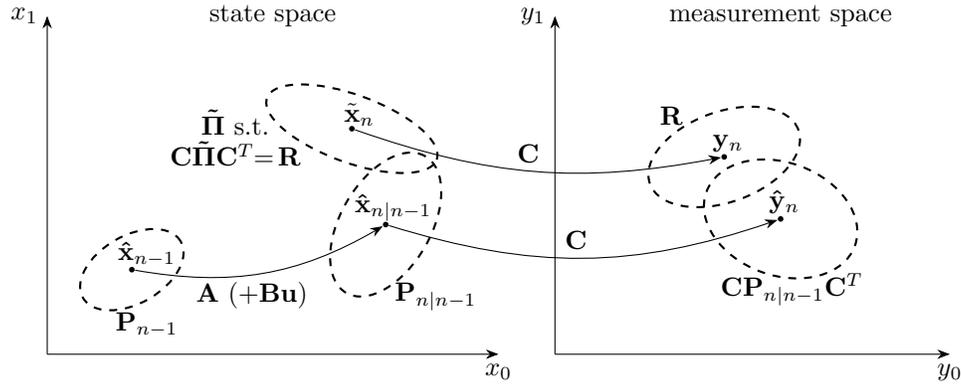
\begin{figure}[!t]\label{fig:kalman-steps}
\centering
\begin{tikzpicture}[scale=1.5, >=Stealth]
  \draw[->] (0, 0) -- (4, 0) node[anchor=north] {$x_0$};
  \draw[->] (0, 0) -- (0, 3) node[anchor=east] {$x_1$};
  \node at (2, 3) {state space};
  \draw[->] (4.5, 0) -- (8, 0) node[anchor=north] {$y_0$};
  \draw[->] (4.5, 0) -- (4.5, 3) node[anchor=east] {$y_1$};
  \node at (6.5, 3) {measurement space};
  
  \draw[thick, dashed] (0.75,0.75) ellipse[x radius=0.5, y radius=0.3, rotate=30];
  \node[fill=black, circle, inner sep=0.75pt, label={[above, xshift=2mm, yshift=-1.25mm]:{$\hbx_{n-1}$}}] (xhat) at (0.75,0.75) {};
  \node[below, xshift=2mm, yshift=-4.5mm] at (0.75,0.75) {$\bP_{n-1}$};
  
  \draw[thick, dashed] (3,1.15) ellipse[x radius=0.7, y radius=0.4, rotate=60];
  \node[fill=black, circle, inner sep=0.75pt, label={[above, xshift=1mm, yshift=-1.25mm]:{$\hbx_{n|n-1}$}}] (Axhat) at (3,1.15) {};
  \node[below right, yshift=-6.5mm] at (3,1.15) {$\bP_{n|n-1}$};

  \draw[thick, dashed] (2.7,2) ellipse[x radius=0.8, y radius=0.3, rotate=-20];
  \node[fill=black, circle, inner sep=0.75pt, label={[above, xshift=1mm, yshift=-1.25mm]:{$\tilde\bx_n$}}] (x) at (2.7,2) {};
  \node[left, yshift=-1.2mm, xshift=-5mm, align=center, text width=18.5mm] at (2.7,2) {$\mathbf{\tilde \Pi}$ s.t. $\bC\mathbf{\tilde \Pi}\bC^\mathrm{T}\!\!=\!\bR$};

  \draw[thick, dashed] (6.5,1.2) ellipse[x radius=0.7, y radius=0.5, rotate=-20];
  \node[fill=black, circle, inner sep=0.75pt, label={[above, xshift=0.5mm, yshift=-1mm]:{$\hby_n$}}] (yhat) at (6.5,1.2) {};
  \node[below, xshift=1.5mm, yshift=-6mm] at (6.5,1.2) {$\bC\bP_{n|n-1}\bC^\mathrm{T}$};

  \draw[thick, dashed] (6,1.75) ellipse[x radius=0.7, y radius=0.4, rotate=20];
  \node[fill=black, circle, inner sep=0.75pt, label={[above, xshift=0.5mm, yshift=-1mm]:{$\by_n$}}] (y) at (6,1.75) {};
  \node[above, xshift=-7mm, yshift=3mm] at (6,1.75) {$\bR$};
  
  \draw[->, bend right=20] (xhat) to node[midway, below, xshift=-1.5mm] {$\bA$ (+$\bB\bu$)} (Axhat);
  \draw[->, bend right=18] (Axhat) to node[midway, above, xshift=-1mm] {$\bC$} (yhat);
  \draw[->, bend right=15] (x) to node[midway, above, xshift=-1mm] {$\bC$} (y);
\end{tikzpicture}
\vspace{-4mm}
\caption{Depiction of the prediction and measurement steps of a Kalman filter. Starting with previous state estimate characterized by $\mathcal{N}(\hbx_{n-1}, \bP_{n-1})$, first predict the state at the current time using only past information as $\mathcal{N}(\hbx_{n|n-1}, \bP_{n|n-1})$. The covariance tends to grow over this step due to the addition of $\bQ$ in $\bP_{n|n-1}$. Then project that prediction into measurement space, $\mathcal{N}(\hby, \bC\bP_{n|n-1}\bC^\mathrm{T})$. Simultaneously, an unknown true state is corrupted by measurement noise with state-space covariance $\mathbf{\tilde \Pi}$ to yield $\tilde\bx_n$ and then projected into measurement space as $\mathcal{N}(\by_n,\bR)$. The new state estimate is then formed by optimally---in both the MSE and MAP (\autoref{sec:MAP}) senses---combining the prediction with the observation to obtain much tighter probability distributions centered where ellipses overlap. Here both state and measurement space are drawn as 2D, but in general they need not have the same dimension.}
\end{figure}

To kick off the filtering process requires: (1) knowledge of the system (what entries of $\bx, \by$ mean; matrices $\bA, \bC$; and additionally $\bu, \bB$ if there are control inputs or synchronized data streams), (2) knowledge of noise covariances ($\bQ, \bR$), and (3) an initial guess of the state and its covariance ($\hbx_0, \bP_0$) to seed the recursion. The filter then proceeds according to the steps of \autoref{algo:kalman-filter-algo}, sketched in \autoref{fig:kalman-steps}.

\begin{algorithm}[!t]
\caption{\textbf{Kalman Filter}}
\label{algo:kalman-filter-algo}
\begin{algorithmic}[1] 
\State assume values for $\hbx_0, \bP_0$ (often treated as \textit{a priori} guesses such that the first two steps of the loop can be skipped on the first iteration)
\For{$n \in \{0, ..., N-1\}$}
    \State calculate predicted next state from past estimate, $\hbx_{n|n-1} = \bA\hbx_{n-1} + \bB\bu_n$
    \State propagate the covariance forward a step, $\bP_{n|n-1} = \bA\bP_{n-1}\bA^\mathrm{T} + \bQ$
    \State calculate predicted measurement, $\hby_n = \bC\hbx_{n|n-1}$
    \State observe true measurement, $\by_n$
    \State calculate the current step's gain, $\bK = \bP_{n|n-1}\bC^\mathrm{T}(\bC\bP_{n|n-1}\bC^\mathrm{T} + \bR)^{-1}$
    \State calculate new state estimate, $\hbx_n = \hbx_{n|n-1} + \bK(\by_n - \hby_n)$
    \State calculate new error covariance, $\bP_n = (\mathbb{I} - \bK\bC)\bP_{n|n-1}$ (by simplification of
\EndFor \hspace{21em}\autoref{eqn:Pn} with $\bK_\text{opt}$)
\end{algorithmic}
\end{algorithm}

\subsubsection{Linear Modeling Example: Cruise Control}
\label{sec:kalman-example}

Before delving further into theory, here is a demonstration of the setup required to use a Kalman filter: states, control inputs, measurement outputs, linear evolution, and noise estimates.

Let $\by$ be noisy position data from simulation of a proportional-integral cruise controller that works to maintain a car's positive constant velocity while driving up and down oscillating hills. To incorporate both the integral feedback control and the action of the hills, state and control inputs are formed as:
\begin{align*}
\bx &= [x_0, x_1, x_2, x_3]^\mathrm{T} = [\text{pos}, \text{vel}, \text{accel}, \text{cumulative pos error}]^\mathrm{T}\\
\bu &= [h, v_d]^\mathrm{T} = [\text{hill slope}, \text{desired velocity}]^\mathrm{T}
\end{align*}

\noindent and the discrete-time system dynamics are agglomerated into the following matrices:
\begin{equation}\label{eqn:cruise_control}
\begin{aligned}
\bA = \begin{bmatrix}
1 & \Delta t & \frac{\Delta t^2}{2} & 0\\
0 & 1 & \Delta t & 0\\
0 & -f_r - \frac{k_p}{\Delta t} & 0 & \frac{k_i}{\Delta t^2}\\
0 & -\Delta t & 0 & 1\\
\end{bmatrix},
&&
\bB = \begin{bmatrix}
0 & 0\\
0 & 0\\
-mg & \frac{k_p}{\Delta t}\\
0 & \Delta t\\
\end{bmatrix},
&&
\bC = \begin{bmatrix}
1 \\ 0 \\ 0 \\ 0\\
\end{bmatrix}^\mathrm{T}
\end{aligned}
\end{equation}

\noindent where $f_r$ represents a friction term, $mg$ is the weight of the car, and $k_p, k_i$ are (unitless) proportional and integral feedback gains, respectively. The first two rows of $\bA$ essentially integrate position and velocity from acceleration. On the last row of $(\bA,\bB)$, estimated and desired velocity at each time step are multiplied by $\Delta t$ to accomplish a simple integration that puts them in units of position, subtracted to create a new position error, and added to the cumulative error. On the third row of $\bA$ and $\bB$, velocity and position error are divided by $\Delta t$ and $\Delta t^2$, respectively, to put them in units of acceleration, then multiplied by appropriate gains. The car's weight is signed negative so that a positive slope (uphill) will cause deceleration and negative slope (downhill) will cause acceleration. Note that friction due to drag is physically quadratic in velocity, but to use the Kalman filter, we have to assume a linear model, perhaps linearizing around a set point.

\begin{figure}[!t]\label{fig:cruise-control-sim}
  \centering
  \includegraphics[width=0.7\textwidth]{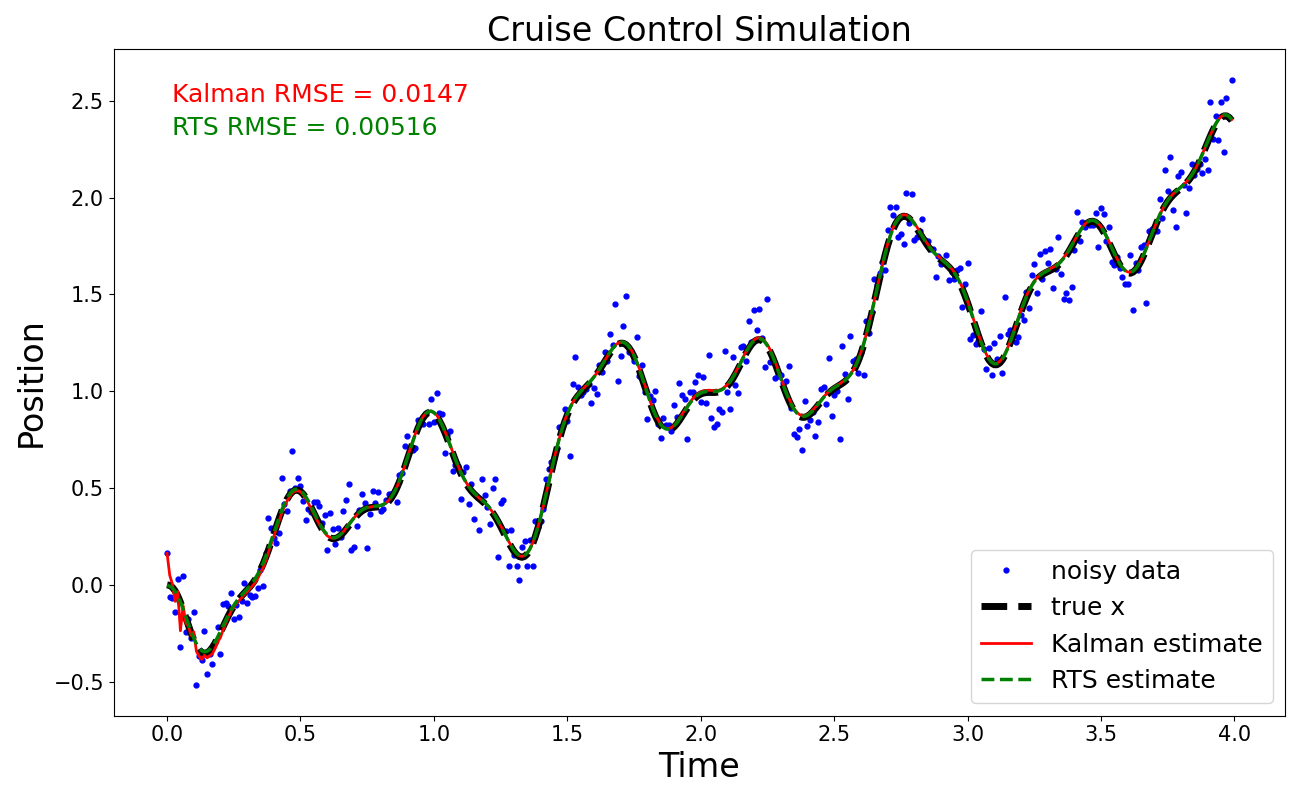}
  \vspace{-3mm}
  \caption{Simulation of \autoref{eqn:cruise_control} with $\Delta t = 0.01$, $mg = 10000$, $f_r = 0.9$, $k_i = 0.05$, $k_p = 0.25$, $v_d = 0.5$, and $h = \big(\sin(2\pi t) + 0.3\sin(8\pi t + 0.5) + 1.2 \sin(3.4\pi t + 0.5)\big)/100$, without noise (black dotted line) and with additive noise $\sim \mathcal{N}(\mu\!=\!0,\ \sigma\!=\!0.1)$ (blue dots). Using the system matrices along with assumed matrices $\bP_0 = 10\cdot\mathbb{I}$, $\bQ = 1000\cdot \Delta t\cdot \text{diag}[(\frac{1}{2}\Delta t^2)^2, \Delta t^2, 1, (\frac{1}{2}\Delta t^2)^2]$, and $\bR = [0.1]$, we run a Kalman filter (red line) and RTS smoother (green dashed line) over the noisy points to produce highly accurate estimates of the true signal (root mean squared errors at upper left). The Kalman filter estimate is thrown around by incoming points, especially at the start.}
\end{figure}

A simulation of this system is given in \autoref{fig:cruise-control-sim}, along with Kalman filter and RTS smoothing (\autoref{sec:rauch-tung-striebel}) estimates using assumed process and measurement noise covariances, $\bQ$ and $\bR$. Gaussian white noise has been added to data points \textit{post hoc}, because process and measurement noise contributions at each step combine to produce a single Gaussian, which is the same for all samples if there is constant step size, $\Delta t$. This leaves the underlying noise covariances unspecified and unknown, so we model them according to a common, rough strategy as diagonal matrices of guessed state dimension variances. This assumes independence between state dimensions, which is inaccurate in this case, but the Kalman filter is robust enough to successfully follow the true signal anyway. See \autoref{sec:kalman-irregular-dt} for discussion of variable step size and derivation of more sophisticated discrete-time $\bQ$.

\subsection{Bayesian Interpretation of the Kalman Filter}
\label{sec:MAP} Denoting probability distributions with $p$, the Maximum Likelihood Estimator (MLE) and Maximum \textit{A Posteriori} (MAP) estimator find parameters or values, $\boldsymbol{\hat\theta}$, given information or measurements, $\by$, according to:
\begin{align*}
\boldsymbol{\hat\theta}_\text{MLE} &= \argmax_{\bth}\ p(\by | \bth)\\
\boldsymbol{\hat\theta}_\text{MAP} &= \argmax_{\bth}\ p(\bth|\by) = \argmax_{\bth}\ p(\by|\bth) p(\bth)
\end{align*}

\noindent where in the latter equation we use Bayes' rule and the fact its denominator does not depend on $\bth$ and therefore makes no difference to the argmax. The MAP is more powerful because it can account for arbitrary priors, and it maps neatly onto the Kalman context.

The Kalman filter can be viewed as fusing evidence from predictions and sensors. We know information from previous steps as well as current measurements and inputs and now wish to find the best estimate for the present state:
\begin{equation}\label{eqn:MAP-kalman-step}
\begin{aligned}
\hbx_n = \argmax_{\bx_n}\ &\underbrace{p(\bx_n | \by_n, \bu_n, \hbx_{n-1})}\\
&=\frac{p(\by_n|\bx_n,\bu_n,\hbx_{n-1})p(\bx_n|\bu_n,\hbx_{n-1})}{p(\by_n|\bu_n,\hbx_{n-1})}\quad\text{by Bayes' rule}
\end{aligned}
\end{equation}

\noindent The probability distribution of $\by_n =  \bC\bx_n + \bv_n,\ \bv_n \sim\mathcal{N}(0,\bR)$ is a Gaussian with covariance $\bR$ centered at $\bC\bx_n$, which is $\propto e^{-\frac{1}{2}\|\bR^{-1/2}(\by_n - \bC\bx_n)\|_2^2}$. The distribution of $\bx_n =\bA\bx_{n-1} + \bB\bu_n + \bw_n,\ \bw_n \sim\mathcal{N}(0,\bQ)$ is recursive and cannot be found directly, but we are given $\bu_n$ and $\hbx_{n-1}$ and thus can construct $\hbx_{n|n-1} = \bA\hbx_{n-1} + \bB\bu_n$. Subtracting, $\bx_n - \hbx_{n|n-1} = \bA(\bx_{n-1} - \hbx_{n-1}) + \bw_n$, which has known covariance $\bP_{n|n-1}$ (derived as first term of \autoref{eqn:cov-P}), so $\bx_n$ is a Gaussian centered at $\hbx_{n|n-1}$, which is ${\propto e^{-\frac{1}{2}\|\bP_{n|n-1}^{-1/2}(\bx_n - \hbx_{n|n-1})\|_2^2}}$. Substituting in \autoref{eqn:MAP-kalman-step} and taking the negative logarithm to get rid of exp and turn the maximization into minimization, we obtain:
\begin{equation}\label{eqn:MAP-optimization}\hbx_n = \argmin_{\bx_n} \|\bR^{-1/2}(\by_n - \bC\bx_n)\|_2^2 + \|\bP_{n|n-1}^{-1/2}(\bx_n - \hbx_{n|n-1})\|_2^2
\end{equation}

\noindent This form is also solvable with calculus. Taking the gradient w.r.t.~$\bx_n$ and setting equal to zero results in $\hbx_n = (\bC^\mathrm{T}\bR^{-1}\bC + \bP_{n|n-1}^{-1})^{-1}(\bC^\mathrm{T}\bR^{-1}\by_n + \bP_{n|n-1}^{-1}\hbx_{n|n-1})$, which can be manipulated into $\hbx_n = \hbx_{n|n-1} + \bK(\by_n-\bC\hbx_{n|n-1})$, although showing equivalence is tedious,\footnote{unless the Woodbury matrix identity and clever factorizations really light your fire} see~\cite{johns-hopkins}.

Indeed, the combination of observation and prior can be seen in the Kalman updates to $\hbx_n, \bP_n$, which calculate the joint probability of $\mathcal{N}(\hbx_{n|n-1},\bP_{n|n-1})$ and $\mathcal{N}(\tilde\bx_n, \mathbf{\tilde \Pi})$ from \autoref{fig:kalman-steps}, but without needing to know (or name) $\mathcal{N}(\tilde\bx_n, \mathbf{\tilde \Pi})$ explicitly. Instead, notice the corresponding joint distribution in measurement space is the product of independent $\mathcal{N}(\hby_n, \bC\bP_{n|n-1}\bC^\mathrm{T})$ and $\mathcal{N}(\by_n, \bR)$, which is a new Gaussian with mean $\bmu_\text{joint} = \bmu_0 + \tilde \bK(\bmu_1 - \bmu_0)$ and covariance $\bSigma_\text{joint} = \bSigma_0 - \tilde\bK\bSigma_0$, where $\tilde \bK = \bSigma_0(\bSigma_0 + \bSigma_1)^{-1}$~\cite{KF-in-pictures, combining-gaussians}. Plugging in $\bSigma_0 = \bC\bP_{n|n-1}\bC^\mathrm{T}$ and $\bSigma_1 = \bR$, we find $\tilde \bK = \bC\bP_{n|n-1}\bC^\mathrm{T}(\bC\bP_{n|n-1}\bC^\mathrm{T} + \bR)^{-1}$, which is $\bC\bK$, so the joint is:
\begin{equation}\label{eqn:meas-space-joint}
\begin{aligned}
\mathcal{N}(\hby_n\!+\!\bC\bK(\by_n-\hby_n),\ \bC\bP_{n|n-1}\bC^\mathrm{T}\!-\!\bC\bK\bC\bP_{n|n-1}\bC^\mathrm{T})\\
=\mathcal{N}\big(\bC(\hbx_{n|n-1}\!+\!\bK(\by_n-\hby_n)),\ \bC(\mathbb{I}\!-\!\bK\bC)\bP_{n|n-1}\bC^\mathrm{T}\big)
\end{aligned}
\end{equation}

\noindent Projection of random vectors is carried out by a left matrix multiply because $\mathbb{E}[\bM\rbchi] = \bM\mathbb{E}[\rbchi]$, while projection of covariance matrices is done by left and right multiply because $\mathrm{cov}(\bM\rbchi, \bM\rbchi) = \bM\mathrm{cov}(\rbchi,\rbchi)\bM^\mathrm{T}$, so the unprojected $\mathcal{N}(\hbx_n,\bP_n)$ back in state space can be found by knocking a few $\bC$s off \autoref{eqn:meas-space-joint}, leaving exactly the last two update equations of \autoref{algo:kalman-filter-algo}.

\subsection{Forward-Backward Kalman Smoothing}
\label{sec:rauch-tung-striebel}

Traditionally, Kalman filters are designed to run online, in the forward direction, but the state, $\hbx$, can contain abrupt changes due to poor \textit{a priori} estimates or discordant observations~\cite{wunsch}. If the goal is to obtain an accurate derivative offline for an entire data series, then filter estimates can be refined by incorporating information from subsequent steps. Such methods can be thought of as solving the MAP problem (\autoref{sec:MAP}) for all steps simultaneously~\cite{aravkin, johns-hopkins}:
\begin{equation}\label{eqn:MAP-kalman-smooth}
\begin{aligned}
\{\hbx_n\} = \argmax_{\{\bx_n\}}\ &\underbrace{p(\{\bx_n\}|\{\by_n\},\{\bu_n\})}\\
&= \frac{p(\{\by_n\}|\{\bx_n\},\{\bu_n\})p(\{\bx_n\}|\{\bu_n\})}{p(\{\by_n\}|\{\bu_n\})}\\
&\propto \Big(\prod_{n=0}^{N-1}\underbrace{p(\by_n|\bx_n,\bu_n)}_{=p_\bv(\bv_n)}\Big)\underbrace{p(\bx_0|\bu_0)}_\text{prior}\prod_{n=1}^{N-1}\underbrace{p(\bx_n|\bx_{n-1},\bu_n)}_{=p_\bw(\bw_n)}
\end{aligned}
\end{equation}

\noindent where braces represent sets of all states, measurements, or inputs across steps, and $p_\bw, p_\bv$ are the noise distributions. In this formulation, the prior of each state transition is conditioned on the raw previous state, not the previous state \textit{estimate} as in the single-step case, because $\bx_{n-1}$ is also being optimized over and is therefore available. Considering the joint probability distribution only up to $\bx_n$ and marginalizing over all $\{\bx_0, ..., \bx_{n-1}\}$ recovers \autoref{eqn:MAP-kalman-step}~\cite{johns-hopkins}.

The evolution of the system over all steps can be set up as a massive sparse linear equation of stacked vectors and block diagonal matrices. These can be used to write joint probabilities and plugged in above to create a formulation similar to \autoref{eqn:MAP-optimization}, but using stacked $\bQ$ instead of the single-step $\bP_{n|n-1}$, because having access to $\bx_{n-1}$ enables direct representation of the probability distribution of $\bx_n$ as a Gaussian with covariance $\bQ$ centered at $\bA\bx_{n-1} + \bB\bu_n$. (See Equation 12 of~\cite{aravkin}.) Taking the gradient w.r.t.~stacked $\{\bx_n\}$ and setting equal to zero to find the minimum yields a block-tridiagonal linear inverse problem~\cite{aravkin, aravkin3, Romberg2016Kalman}.

Rauch-Tung-Striebel (RTS) smoothing~\cite{RTS, wunsch, aravkin} is a classic algorithm to efficiently solve said system one block at a time, working sequentially backward using information from a forward Kalman filter pass:
\begin{equation}\label{eqn:rts}
\begin{split}
    \bL_n &= \bP_n \bA^\mathrm{T} \bP_{n+1|n}^{-1}\\
    \mathbf{\widehat{x}}_{n,\text{RTS}} &= \hbx_n + \bL_n(\widehat{\bx}_{n+1,\text{RTS}} - \hbx_{n+1|n})\\
    \bP_{n,\text{RTS}} &= \bP_n + \bL_n(\bP_{n+1,\text{RTS}} - \bP_{n+1|n})\bL_n^\mathrm{T}
\end{split}
\end{equation}

\noindent Remembering estimates of state and error covariance from every step to support these computations can pose a memory challenge for high-dimensional states or long sequences. When system matrices are constant, $\bL$, like $\bP$, converges to a steady state, which can optionally be precalculated to save work and memory.

If the system dynamics, $\bA$, are invertible, then it is also possible to run the Kalman filter, or whole RTS smoothing, on the data in reverse. Note $\bQ$ is left as is for this procedure, because uncertainty due to entropy should still be modeled as increasing across steps, even if time is reversed. Initial guesses $\hbx_0,\bP_0$ bias Kalman estimates for several iterations, so combining forward and reverse passes can sometimes improve results at domain edges.

 \subsection{Generalized Kalman Smoothing for Robust Estimation}
\label{sec:robust-estimation}

The Bayesian interpretation of the Kalman filter and the consequent MAP optimization problem (Equations \ref{eqn:MAP-kalman-step} and \ref{eqn:MAP-kalman-smooth}) provide a flexible framework, allowing us to consider alternative noise distributions and distance metrics through modifications of the cost function~\cite{aravkin, czech}, e.g.~via loss functions $V$ and $J$:
\begin{equation}\label{eqn:MAP-opt-generalized}
\begin{aligned}
\{\hbx_n\} = \argmin_{\{\bx_n\}}& \sum_{n=0}^{N-1}V(\bR^{-1/2}(\by_n - \bC\bx_n))\\
&+ \sum_{n=1}^{N-1}J(\bQ^{-1/2}(\bx_n - \bA\bx_{n-1}-\bB\bu_n)) - \ln(p(\bx_0))
\end{aligned}
\end{equation}

For example, assuming measurement noise drawn from a Laplace distribution, $\propto e^{-\|\boldsymbol\cdot\|_1}$, the corresponding term of the MAP optimization problem (first term of Equations \ref{eqn:MAP-optimization} and \ref{eqn:MAP-opt-generalized}) gets an $\ell_1$ norm. Or, if the data has outliers, then because the $\ell_2$ norm highly penalizes larger inaccuracies, it can bend the state estimate toward outliers too strongly~\cite{ddse}. In this case we may wish to replace the distance metric around the measurement term with something like a Huber loss~\cite{Huber1964}, which is quadratic only up to some radius, $M$, and then grows linearly:
\begin{equation}\label{eqn:huber}
\text{Huber}\footnote{The Huber loss is typically defined with a scalar input, $x$, and is applied to vectors pointwise. The $\sum\text{Huber}(\,\cdot\,,\, M)$ function can be thought of as an interpolation between $\ell_2$ and $\ell_1$ norm terms, because when properly normalized by multiplication with $c_2 = ([4e^{-M^2/2}\frac{1+M^2}{M^3} + \sqrt{2\pi}(2\Phi(M) - 1)]/[2e^{-M^2/2}/M + \sqrt{2\pi}(2\Phi(M) - 1)])^{1/2}$~\cite{aravkin3}, where $\Phi$ is here the cumulative distribution function of the standard normal distribution, $\sum\text{Huber}$ reduces to $\frac{1}{2}\|\cdot\|_2^2$ as $M \to \infty$ and to $\sqrt{2}\|\cdot\|_1$ as $M \to 0$, where $\frac{1}{2}$ and $\sqrt{2}$ are the ``statistically correct" scales for those terms~\cite{aravkin3}. Normalization constants were irrelevant in earlier sections because the cost function's terms had the same form, but when mixing loss functions of different shapes, their relative scales matter to the optimization.}(x,\, M) = \begin{cases} \frac{1}{2}x^2 & |x| \leq M \\ M|x| - \frac{1}{2}M^2 & |x| > M\end{cases}
\end{equation}

\noindent Or, in the presence of relatively rare, sudden state changes, the $\ell_2$ norm cannot follow fast jumps~\cite{ohlsson}, so it may be desirable to use the $\ell_0$ norm on the process (second term of Equations \ref{eqn:MAP-optimization} and \ref{eqn:MAP-opt-generalized}), or in practice its similarly sparsity-promoting convex relaxation, the $\ell_1$ norm.

\begin{figure}[!t]\label{fig:robust-demo}
\centering
\includegraphics[width=0.7\textwidth]{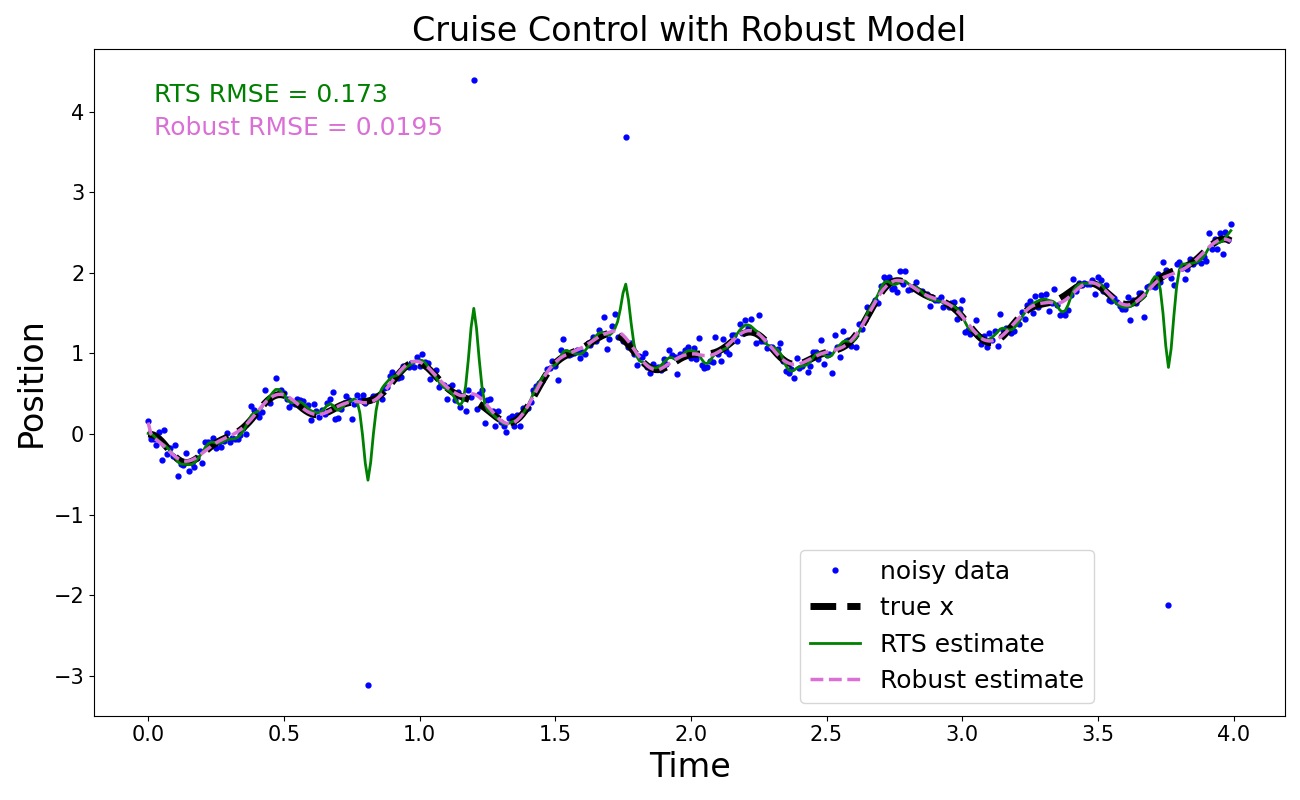}
\vspace{-3mm}
\caption{RTS smoother (green) and robust MAP smoother (dashed pink) estimates on the cruise control data from \autoref{fig:cruise-control-sim}, but this time with 4 randomly chosen points turned into outliers. All matrices remain the same as before except for $\bQ$, which has been scaled up by a factor of $10^5$ to coerce the solution to deform significantly toward outliers, because RTS smoothing does remarkably well even with outliers when the system matrices are accurate. Despite such a poor choice of $\bQ$, a robust smoother with process penalty $\sqrt{2}\|\cdot\|_1$ and measurement penalty $1.04\sum\text{Huber}(\cdot\,,\, 2)$ (constants explained in footnote) is able to ignore outliers and obtain an estimate close to ground truth.}
\end{figure}

These substitutions can make a big difference to performance, as exhibited in \autoref{fig:robust-demo}, but they typically break the ability to find a closed-form recursive solution by calculus. Fortunately, we can use modern tools like convex optimization to solve the problem efficiently anyway~\cite{Boyd_Vandenberghe_2004}---in fact, time complexity is \textit{linear} thanks to sparse, block-tridiagonal systems~\cite{aravkin3, aravkin}, similar to the previous section. We can also exploit techniques from this framework to do things like enforce bounds on the state via inequality constraints~\cite{aravkin}.

There are still limits down this path: Nonlinear models result in nonconvex optimization problems that may or may not be reasonable to solve with alternative approaches, e.g.~Gauss-Newton, and some noise distributions, like Student's t, have merely quasiconvex negative log, necessitating iterative solution procedures~\cite{aravkin2, Boyd_Vandenberghe_2004}. Still, this is a powerful generalization when sophisticated system models are available, and we recommend the review given in~\cite{aravkin}, which includes a tutorial on optimization techniques.

\subsubsection{Alternative Notion of Robustness}

In controls and signal processing, the term ``robust" is also used to mean \textit{guaranteed} performance under input or model (parametric) uncertainty~\cite{robust-control, dullerud2000course, xie}. For unknown inputs, this gives rise to $H_\infty$ filtering/smoothing~\cite{nagpal, batteries}, which is related to $H_2$ filtering~\cite{robust-control}, of which the Kalman filter is the most famous special case.\footnote{To cast the Kalman filter into form of \autoref{eqn:H-G-system}, the salient output and internal state are \textit{both} the error $\be$, with dynamics as derived in \autoref{eqn:kalman-min-error}. Vector-valued disturbance can be split into pieces and multiplied by half-powers of covariance matrices $\bQ$ and $\bR$ to produce appropriately-shaped-but-deterministic (w.r.t.~$\tbw$) process and measurement disturbances:
$$\begin{bmatrix} \be_n \\ \be_{n-1} \end{bmatrix} = \begin{bmatrix} (\mathbb{I} - \bK\bC)\bA & \begin{bmatrix} (\mathbb{I} - \bK\bC)\bQ^{1/2} & -\bK\bR^{1/2}\end{bmatrix} \\ \mathbb{I} & 0 \end{bmatrix} \begin{bmatrix} \be_{n-1} \\ \begin{bmatrix} \tbw_\text{proc} \\ \tbw_\text{meas} \end{bmatrix}_n \end{bmatrix}$$Using the time-domain version of the $H_2$ norm, which is defined with G's impulse response (see~\cite{robust-control}), along with this ($\tbA, \tbB, \tbC, \mathbf{\tilde D}$), we can manipulate the cost function of \autoref{eqn:H-G-optimization} into $\text{Tr}(\tbC\bP\tbC^\mathrm{T})$, where $\bP$ is the controllability Gramian $=\sum_{m=0}^\infty\tbA^m\tbB\tbB^\mathrm{T}(\tbA^m)^\mathrm{T}$, which solves the discrete-time Lyapunov equation, $\tbA\bP_{n-1}\tbA^\mathrm{T} -\bP_n + \tbB\tbB^\mathrm{T} = 0$. Algebra recovers Equations \ref{eqn:kalman-min-tr-P} and \ref{eqn:Pn}, modulo the probabilistic bits. Choosing $\tbw \sim \mathcal{N}(0,\mathbb{I})$ puts the Gramian sum for $\bP$ in the form of an expected value, $\mathbb{E}[\cdot]$, completing the bridge.} These methods, perhaps confusingly, are unrelated to MAP generalizations of the Kalman filter. Instead, they consider a linear system named G, which maps exogenous ``disturbance" $\tbw$ and internal state to output error $\tbz$ and next state (or state derivative in continuous time):
\begin{equation}\label{eqn:H-G-system}
\begin{bmatrix}\mathbf{\tilde x}_{n+1}\\ \tbz_n\end{bmatrix} = 
\setlength{\arrayrulewidth}{0.1pt} {\renewcommand{\arraystretch}{1.3}
\underbrace{\left[\begin{array}{c|c}
\tbA & \tbB \\ \hline
\tbC & \mathbf{\tilde D}
\end{array}\right]}_{\text{\normalsize G}}}
\begin{bmatrix} \mathbf{\tilde x}_n \\ \tbw_n\end{bmatrix}
\end{equation}

\begin{figure}\label{fig:G-H-system}
\centering
\begin{tikzpicture}
  \node[draw, minimum height=1cm] (P) {plant};
  \node[draw, below=5mm of P, minimum height=1cm] (K) {controller};
  \draw[->] ([yshift=-2mm] P.east) -- +(5.5mm,0) |- (K.east);
  \draw[->] (K.west) -- +(-2mm,0) |- ([yshift=-2mm] P.west);
  \node[draw, dashed, fit=(P)(K), minimum width=2.7cm, inner sep=5pt, rounded corners, label=above:G] (G) {};
  \draw[->] +(-2cm,2mm) node[left] {$\tbw$} -- ([yshift=2mm] P.west);
  \draw[->] ([yshift=2mm] P.east) -- +(1.5cm,0) node[right] {$\tbz$};
\end{tikzpicture}
\vspace{-1mm}
\caption{In the framework of H norms, system G takes in a disturbance signal and produces an error signal. Internally, G can have constituent components that affect its behavior, which are factored into the overall \autoref{eqn:H-G-system}.}
\end{figure}
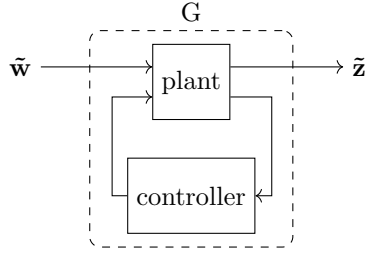

\noindent As illustrated in \autoref{fig:G-H-system}, the system can have its own internal logic, such as a ``plant" with self-contained state and linear dynamics connected to a controller (or observer). Signals in this framework are typically taken to be $\in L_2$ so they have nice inner products, and the complex-valued transfer function that maps between signals is constrained to a Hardy space, $H_p$, to ensure its norm integrals are bounded.

We then form the following optimization problem:
\begin{equation}\label{eqn:H-G-optimization}
\begin{aligned}
\argmin_\bK\ \ &\| \text{G}\|_{H_p}\\
\text{s.t.}\ \ &\text{G's dynamics satisfied}
\end{aligned}
\end{equation}

\noindent where $\bK$ are tunable parameters in the system and $\|\cdot\|_{H_p}$ denotes a Hardy norm. We care about $p=2$, because $\|\text{G}\|_{H_2}$ is the inner product of G's transfer function with itself in Hardy space, which integrates the response to sinusoidal inputs over all frequencies, thereby penalizing the average response~\cite{robust-control}. We also care about $p=\infty$, where $\|\text{G}\|_{H_\infty} \triangleq \sup_\omega \sigma_\text{max}(\text{G})$,\footnote{In practice, this cannot be solved directly, but the Kalman–Yakubovich–Popov lemma provides a way to check whether $\|\text{G}\|_{H_\infty} < \gamma$ using a Linear Matrix Inequality~\cite{boyd-lmis}, which is a convex problem. If infeasible, raise $\gamma$ to loosen, and if feasible, step $\gamma$ down to tighten; it's quasiconvex.} which penalizes the worst-case response, the system's most sensitive mode at its most sensitive frequency~\cite{robust-control}. A result from functional analysis says the $H_\infty$ norm also gives $\sup_{\tbw \neq 0} \frac{\|\tbz\|_2}{\|\tbw\|_2}$, i.e.~the maximal ratio of output magnitude (error) to input magnitude (disturbance), including arbitrary, non-sinusoidal disturbances~\cite{robust-control}, making it a useful design criterion when error amplification must never exceed a hard bound.

However, this mathematically intense framework is often weaker than desired: $H_p$-optimal solutions are optimal \textit{regardless} of any particular form of bounded-energy disturbance. This can be counterintuitive, especially because it treats disturbances as deterministic, whereas from the outset we have considered random \textit{noise} (\autoref{sec:noise}). The classic Kalman filter is always $H_2$-optimal because it minimizes the system norm, but we understand it is also the ``best" estimator in a probabilistic sense under certain conditions and have seen examples in \autoref{sec:robust-estimation} of extending that same notion of best (MAP) to alternative assumptions. This requires filter modifications that make the system into something new, no longer \textit{the} $H_2$-optimal option. Similarly, whether the $H_\infty$ filter is appropriate depends on the scenario. Some systems that can be treated with Kalman filtering have unbounded $H_\infty$ norm, making $H_\infty$ filter design infeasible, and an $H_\infty$ filter places more weight on measurements and less on priors than a Kalman filter, allowing more high-frequency noise through~\cite{crassidis_junkins}. Thus, unless truly nothing is known about the disturbance (a rare situation), this notion of robustness is not the right tool for model-based filtering: A smoother estimate, more amenable to differentiation, will be possible under very bare assumptions, like wide uniform priors or anticipation of outliers.

\subsection{Nonlinear Filtering}
\label{sec:UKF}

The methods covered so far are all limited to linear models, but many real-world systems evolve according to nonlinear equations. This can even include linear systems with nonlinear parametric dependency if parameters are estimated as part of the state. For example, say the gain terms $k_p$ and $k_i$ in the cruise control example from \autoref{sec:kalman-example} are unknown, then the system can be considered to have the following state and dynamics:

\begin{align*}
\bx &= [x_0, x_1, x_2, x_3, k_i, k_p]^\mathrm{T} \\
\bx_n &= \bff(\bx_{n-1},\bu_n),\quad \bff(\bx,\bu)=\begin{bmatrix}
x_0 + \Delta t\ x_1+ \frac{\Delta t^2}{2}x_2\\
x_1 + \Delta t\ x_2\\
-f_r x_1 \textcolor{blue}{\,-\, \frac{k_p}{\Delta t} x_1 + \frac{k_i}{\Delta t^2} x_3} - mg\, u_0 + \frac{k_p}{\Delta t} u_1 \\
x_3 - \Delta t\, x_1 + \Delta t\, u_1 \\
k_i \\
k_p
\end{bmatrix} \\
\by_n &= \bh(\bx_n,\bu_n) = x_{n,0}
\end{align*}
    
\noindent where blue font indicates the nonlinear terms. 

Addressing nonlinear systems in the Kalman framework requires mathematical extension. The oldest and simplest approach is fittingly called the Extended Kalman Filter (EKF), which forward-solves the process dynamics given by $\bff(\bx,\bu)$ to predict the \textit{a priori} state estimate and uses the measurement function, $\bh(\bx,\bu)$, to get the current measurement, but then uses \textit{Jacobians}\footnote{If $\bff$ and $\bh$ are linear functions, expressible as products of coefficient matrices and variable vectors, then their Jacobians are simply the matrices, and the EKF reduces to the KF.} evaluated at the current state estimate in place of system matrices in the Kalman filter equations to find error covariance, $\bP$, and gain, $\bK$:

$$\bA_n = \frac{\partial \bff(\bx,\bu)}{\partial \bx} \Bigg|_{\hat{\bx}_{n-1}, \bu_n},\quad
\bB_n = \frac{\partial \bff(\bx,\bu)}{\partial \bu} \Bigg|_{\hat{\bx}_{n-1}, \bu_n},\quad
\bC_n = \frac{\partial \bh(\bx,\bu)}{\partial \bx} \Bigg|_{\hat{\bx}_{n|n-1}, \bu_n}$$

\noindent Unlike the linear Kalman filter, the EKF is not optimal and not guaranteed to converge; its first-order linearizations may not be able to keep up with highly nonlinear dynamics.

\begin{figure}[!t]\label{fig:ukf}
  \centering
  \includegraphics[width=0.6\textwidth]{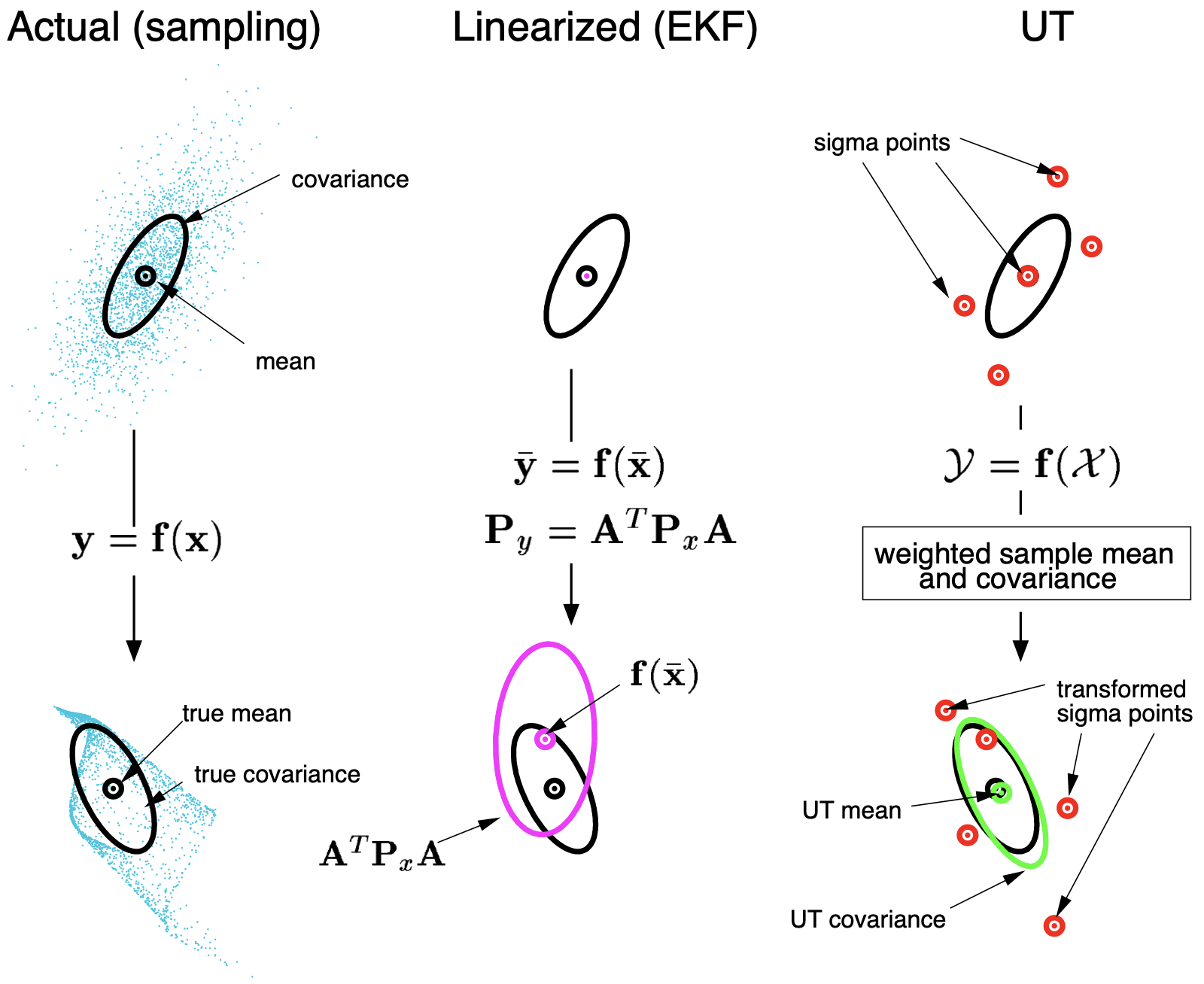}
  \vspace{-3mm}
  \caption{From~\cite{oregon}. Left: Particles sampled from a Gaussian are passed through a nonlinearity to produce a non-Gaussian distribution with calculable mean and covariance. Middle: The EKF uses a local linear approximation of the nonlinearity to propagate mean and covariance inaccurately. Right: Specially-chosen points are passed through the dynamics, and the results are weighted according to the UKF scheme to produce a much better approximation of the true mean and covariance, although non-Gaussianity of the final distribution is not captured.}
\end{figure}

Further innovations have been made to better handle acute nonlinearities, the most famous being the Unscented Kalman Filter (UKF)~\cite{julier-uhlmann,julier-uhlmann2}, which predicts subsequent error covariance\footnote{In the original UKF, the error covariance matrix, $\bP$, can lose positive definiteness due to numerical issues. This can be largely resolved by using the square-root implementation, which guarantees the error covariance remains positive definite by estimating and updating the Cholesky factor, $\mathbf{S}$ such that $\bP=\mathbf{SS}^\mathrm{T}$, instead of $\bP$ itself~\cite{merwe2001sqrt}.} and state by propagating specially-chosen (``sigma") points through the actual nonlinear dynamics. This can be thought of as a kind of particle filter~\cite{russell-norvig}, but the UKF enforces a Gaussian state estimate and exploits this structure to make next estimates using a linear (in the state dimension) number of points, whereas a generic particle filter requires an exponential number of particles to represent arbitrary probability distributions, per the curse of dimensionality. \autoref{fig:ukf} offers a visual comparison. This method ``captures the posterior mean and covariance accurately to the $3^\text{rd}$ order for \textit{any} nonlinearity"~\cite{oregon}. For improved estimates in the \textit{post hoc} regime, an RTS-like smoothing algorithm based on the UKF was derived by~\cite{sarkka}.

\subsection{Recommendations}

Due to Gauss-Markov theorem, the standard Kalman filter is the Best Linear Unbiased Estimator (BLUE) for noise distributions with zero mean and known covariance that are uncorrelated through time (white) and independent of the system state (homoscedastic)~\cite{aravkin, Romberg2016Kalman}. For linear systems, the Kalman filter is also the most \textit{efficient} unbiased estimator of any kind, as the error covariance converges to the Cramér-Rao lower bound~\cite{crassidis_junkins,Cramer1946}. Coupling these facts with the central limit theorem, which states the combined influence of many random sources tends toward a Gaussian, the standard Kalman filter is surprisingly effective across a wide array of practical contexts.

Even so, the classical approach has limits: If data is contaminated with outliers or non-Gaussian noise, or when the system does not evolve linearly, the most optimal estimator becomes biased or nonlinear. Sophisticated generalizations and extensions of Kalman filtering can be complicated to set up and solve, especially as concerns like nonlinearity and robustness stack together, making the mathematics ever more involved. This remains an active area of research with a dizzying number of filter variants~\cite{chilean, friston} and novel approaches~\cite{martin, variational-kalman}.

At an even more basic level, system modeling can be a challenge, and a model-based approach may not confer an advantage versus the model-free paradigm (\autoref{sec:noisy-without-knowledge}), depending on accuracy of the model. Automatic model discovery techniques exist but are difficult to apply to the differentiation problem, because they themselves rely on good derivative estimates~\cite{sindy, pde-find, ddse}. In some cases, blending a partially-known model with the performance metrics of \autoref{sec:performance-metrics} may help in estimating unknown parameters like noise levels.

\section{Differentiating Noisy Data Without Prior Knowledge}
\label{sec:noisy-without-knowledge}

In the absence of a model, estimating the derivative of a function, $x(t)$,\footnote{It is common in the context of noisy data to keep with the Kalman filter notation of \autoref{sec:noisy-with-knowledge}, where the symbol $y$ is reserved for noisy measurements and underlying signal is called $x$, so the independent variable must take on another symbol, $t$. This contrasts with earlier sections of this review, where functions are usually named $y$, taking input $x$.} given noisy samples, $\by = \bx + \boldsymbol{\eta} = x(\{t_n\}) + \eta(\{t_n\})$, where $\{t_n\}$ is a set of sample points and $\eta$ is a noise process, is an ill-posed problem. In other words, the data are ambiguous, so it is not possible to solve for the derivative without imposing prior assumptions. Nevertheless, practitioners must solve this ill-posed problem routinely, often without much guidance.

In data-driven science, for instance, model discovery via SINDy and its PDE extension~\cite{sindy, pde-find} is accomplished by regressing, with sparse selector $\mathbf{\Xi}$, a library of candidate terms $\mathbf{\Theta}$, dependent on data $\mathbf{X}$, against derivatives $\mathbf{\dot{X}}$, according to $\mathbf{\dot{X}} = \mathbf{\Theta}(\mathbf{X}, ...) \mathbf{\Xi}$, where the ellipsis denotes optional dependence on temporal or spatial derivatives of $\mathbf{X}$. Accuracy is contingent on estimates of true quantities of interest and their derivatives, which must be found from noisy samples alone.\footnote{It is sometimes possible to sidestep direct estimates: Weak SINDy integrates both sides against smooth, compactly-supported test functions and shifts derivatives onto those test functions by parts, exactly the move of \autoref{sec:weak-form}, so that integration averages noise rather than amplifying it~\cite{weak-sindy}. But not all models can be rearranged to accommodate the trick.}

\subsection{Performance Metrics in the Presence and Absence of Ground Truth}
\label{sec:performance-metrics}

A whole cornucopia of methods have been developed to estimate derivatives of noisy data, each with a set of tunable hyperparameters that affect its behavior, such as polynomial degree, window width, and number of iterations, which we collectively call $\Phi$. Differentiating in the absence of prior knowledge comes down to choosing a method and hyperparameter settings that achieve a good answer, where ``good" entails two major concerns: closeness of fit and complexity~\cite{apriori-boulder}. Following~\cite{floris2020}, we use the following metrics:

\vspace{2mm}\begin{enumerate}
  \item Accuracy: The derivative estimate, $\hdbx$, should be faithful to the true derivative, $\dbx$, resulting in a small root mean squared error:
  \begin{equation}\label{eqn:rmse}
  \text{RMSE} = \sqrt{\frac{1}{N}\sum{(\hdbx-\dbx)^2}}
  \end{equation}
  \item Bias: The derivative estimate should be unbiased in the sense its accuracy should not depend on the underlying true value of the derivative. Bias occurs in low-pass filtering, for instance, because it not only removes noise but also dulls sharp slopes, causing systematic under-prediction of larger-magnitude derivative values. We can measure bias with the squared correlation coefficient between estimate error, $\hdbx-\dbx$, and the actual derivative:
  \begin{equation}\label{eqn:ec-r2}
  R^2 = \frac{\text{cov}(\hdbx-\dbx,\; \dbx)^2}{\text{var}(\hdbx-\dbx)\,\text{var}(\dbx)}
  \end{equation}
\end{enumerate}

Both the above metrics require knowledge of the true derivative, $\dbx$, which in real-world applications is usually unknown, so the authors of~\cite{floris2020} suggest a proxy loss function that uses only measurements $\by$ and the assumption that noise is largely band-separable from signal (\autoref{sec:noise}) to balance fidelity and smoothness:
\begin{equation}\label{eqn:floris-cost}
L(\Phi) = \text{RMSE}\big(\text{trapz}(\hdbx(\Phi)) + \mu,\ \by\big) + \gamma\, \text{TV}\big(\hdbx(\Phi)\big)
\end{equation}

\noindent where vector $\hdbx(\Phi)$ is the estimate returned by a differentiation method parameterized by $\Phi$; trapz is numerical integration with trapezoidal rule, which is less sensitive to noise than higher-order quadrature rules~\cite{pynumdiff}; $\by$ are observed noisy measurements; $\mu$ corrects for a lost constant of integration by solving $\argmin_\mu \| \text{trapz}(\hdbx) + \mu - \by \|_2$, which is just the mean difference between the two vectors; TV is normalized Total Variation (\autoref{sec:tvr}), equal to $\frac{1}{N}\|\hdbx_{0:N-1} - \hdbx_{1:N}\|_1$; and $\gamma$ is a hyperparameter, with larger values of $\gamma$ favoring smoother derivative estimates.

In the presence of outliers, RMSE can weight spurious points too highly, so we further suggest a generalization of \autoref{eqn:floris-cost} using the Huber loss (\autoref{eqn:huber}):
\begin{equation}\label{eqn:pavel-cost}
L(\Phi) = \sqrt{\frac{2}{N}\sum\text{Huber}\Big(\text{trapz}(\hdbx(\Phi)) + c -\ \by,\, M\sigma_\text{MAD}\Big)} + \gamma\,\text{TV}\big(\hdbx(\Phi)\big)
\end{equation}

\noindent where Huber is applied pointwise to vector inputs; $\sigma_\text{MAD}$ is a measure of scatter in the residuals (first Huber argument) using median absolute deviation (MAD), normalized by the median absolute value of the standard normal distribution, $\mathcal{N}(0, 1)$, so its scale matches that of standard deviation for Gaussian errors; $M$ parameterizes where the Huber loss changes from quadratic to linear, in units of $\sigma_\text{MAD}$; and $c$ is an integration constant, this time estimated robustly by solving $\argmin_c \sum\text{Huber}(\text{trapz}(\hdbx) + c - \by,\, M\sigma_\text{MAD})$. As $M \to \infty$, \autoref{eqn:pavel-cost} reduces to \autoref{eqn:floris-cost}, but for even moderate $M$ the vast majority of data points are treated with RMSE; e.g., assuming Gaussian inliers, the portion beyond $M\sigma$ can be found with \texttt{2*(1-scipy.stats.norm.cdf(M))}, which is about 0.0455 for $M = 2$ and 0.0027 for $M = 3$.

Due to the complex relationship between $\Phi$ and $\hdbx$ through a differentiation method, Equations \ref{eqn:floris-cost} and \ref{eqn:pavel-cost} are unfortunately not convex in $\Phi$, nor are gradients with respect to $\Phi$ easy to find. But they can still be optimized via strategies such as Nelder-Mead~\cite{nelder-mead}, a downhill simplex direct search, starting from multiple initial conditions~\cite{pynumdiff, floris2020}.

\begin{figure}[!t]\label{fig:pareto}
  \centering
  \includegraphics[width=0.3\textwidth]{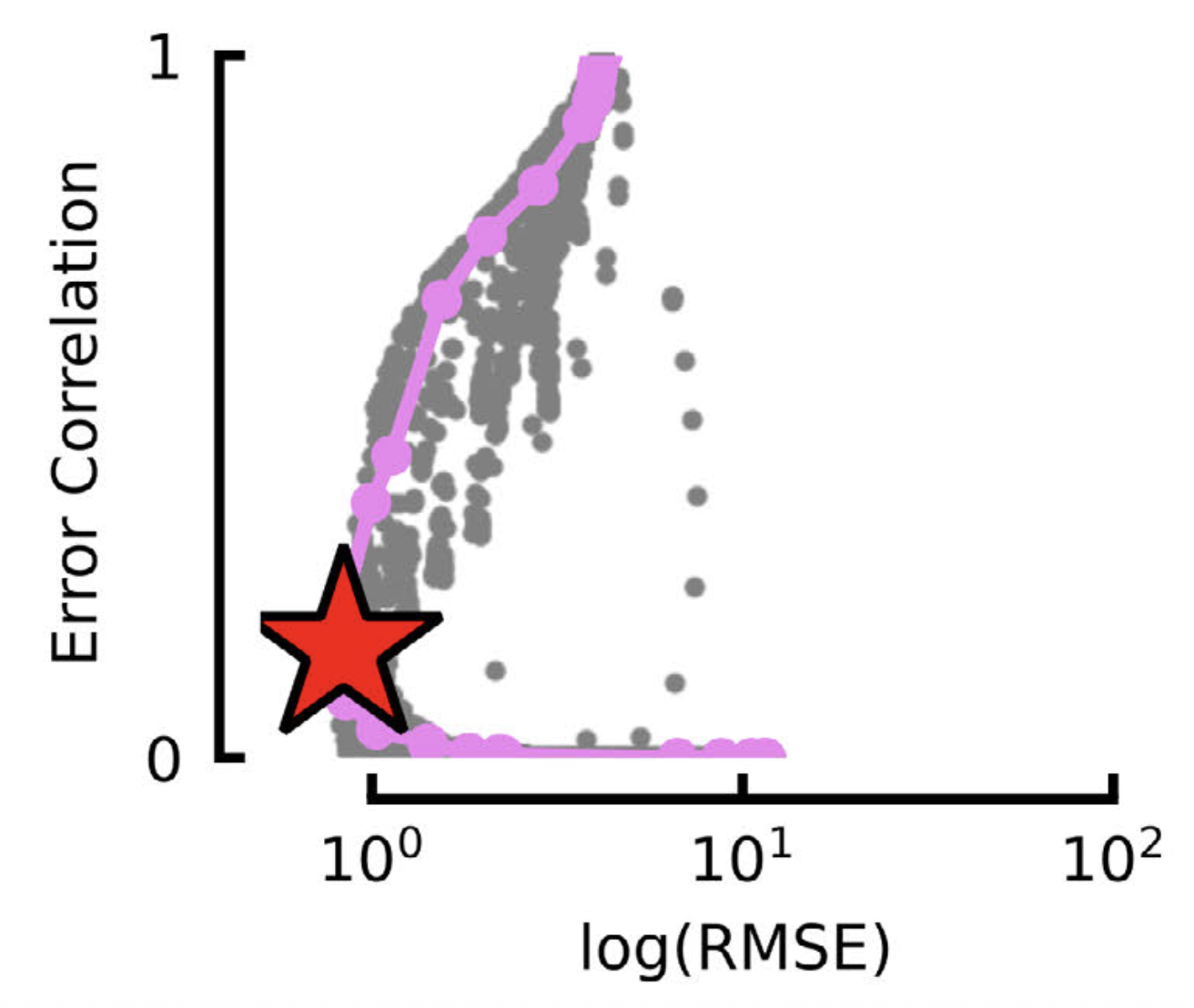}
  \caption{From~\cite{floris2020}, the (RMSE,$R^2$)-plane. Grey dots are evaluations of Savitzky-Golay (\autoref{sec:savitzky-golay}) estimates of the derivative for various choices of hyperparameters, $\Phi$. Choosing $\gamma$ and running Nelder-Mead optimization over \autoref{eqn:floris-cost} results in solutions along the pink line. Choosing $\gamma$ according to the recommended heuristic results in the solution at the red star.}
\end{figure}

To justify \autoref{eqn:floris-cost}, the authors of~\cite{floris2020} apply a variety of differentiation methods to a suite of synthetic problems, where $\dbx$ can be known and used as ground truth. Different hyperparameter settings, $\Phi$, yield estimates that fall at different points of the (RMSE,$R^2$)-plane, and these data points exhibit a definite Pareto front, as shown in \autoref{fig:pareto}. Assuming values for $\gamma$ and optimizing \autoref{eqn:floris-cost} yields hyperparameters that produce derivatives very near the Pareto front, and this result holds across methods and problems.

Choosing $\gamma$ to achieve both low RMSE and $R^2$ is found to depend on the temporal resolution and frequency content of the data and not on function shape, noise level, or differentiation method, and a heuristic formula is provided to determine its value~\cite{floris2020}:
\begin{equation}\label{eqn:gamma-heuristic}
\gamma = f^{-1.6}\Delta t^{-0.71} e^{-5.1}
\end{equation}
where $f$ is the signal bandlimit in Hz, i.e.~the highest frequency expected to come from signal rather than from noise, and the constants come from a best-fit regression. The idea is that when a signal contains higher frequencies, its derivative cannot be smoothed as aggressively.

Armed with the knowledge that smoothing is indeed the right approach, the methods in this section can broadly be thought of as implementing smoothing in a number of different ways.

\subsection{Prefiltering Methods}
\label{sec:prefilt}

The simplest way to smooth data is as an independent first step, before calculating a derivative. This is a heavy-handed operation, so although global differentiation methods offer an error-bound advantage versus local ones~\cite{ahnert}, prefiltering erases a lot of information and nullifies the difference. As such, simple Finite Difference is ideal for this case, cheap to compute and imposing no extra assumptions.

\begin{figure}[!t]\label{fig:freq-responses}
  \centering
  \includegraphics[width=\linewidth]{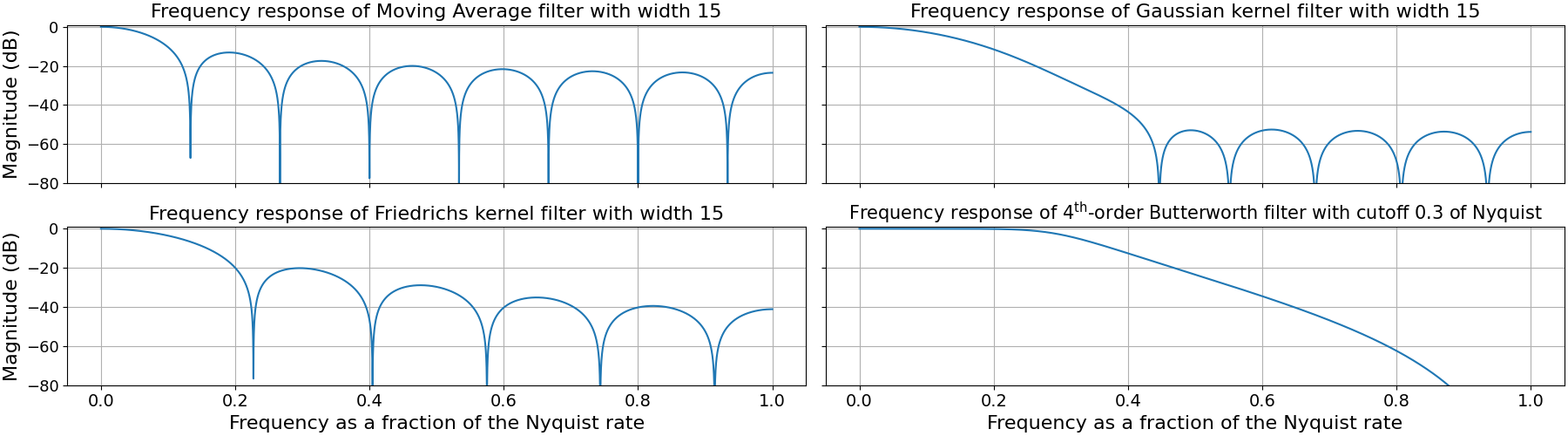}
  \vspace{-5mm}
  \caption{Frequency responses of possible prefilters: moving average (top left), Gaussian kernel (top right), Friedrichs kernel (lower left), and Butterworth filter (lower right). Plots are generated by \texttt{scipy.signal.freqz}, with frequency given as a fraction of the Nyquist rate by scaling angular frequency $\omega \in [0, \pi]$ by $\frac{1}{\pi}$. Frequency caps at the Nyquist rate, $\frac{f_s}{2}$, because higher frequencies alias down into this band (as in \autoref{fig:noise-spectrum}). All methods demonstrate low-pass character, but only Butterworth is purpose-designed, while the others exhibit undesirable humps in higher frequencies. The bandwidth is said to stop where power output is halved, about -3 decibels per $10 \log_{10}(\frac{1}{2})$.}
\end{figure}

A common method of smoothing is to simply take a moving average, which is equivalent to convolving with a discrete kernel of uniform values that sum to 1. We can also take the median over windows, which can remove spurious values while preserving sharp transitions if noise levels are not too high, or we can convolve with alternative normalized kernels sampled from Gaussian or Friedrichs\footnote{a bump function, e.g.~$e^{-1/(1-x^2)}$, $x \in (-1, 1)$} functions. These kernel methods constitute crude digital low-pass filters, with frequency responses (magnitude Bode plots~\cite{oppenheim}) shown in the first three panels of \autoref{fig:freq-responses}. The median filter is not visualized, because it is nonlinear and lacks a well-defined response.

\begin{figure}[!t]\label{fig:kernel-cc}
  \centering
  \includegraphics[width=0.99\textwidth]{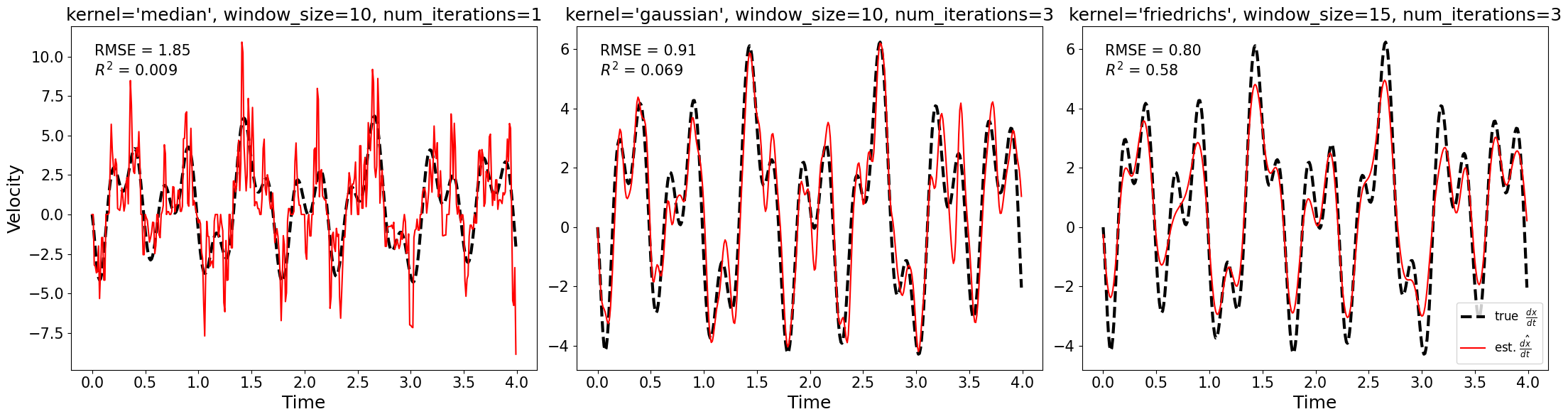}
  \vspace{-3mm}
  \caption{Derivative estimates for noisy cruise control data (\autoref{fig:cruise-control-sim}) found by applying a kernel in a sliding window followed by second-order Finite Difference, using three different choices of hyperparameters.}
\end{figure}

\begin{figure}[!t]
  \centering
  \includegraphics[width=0.99\textwidth]{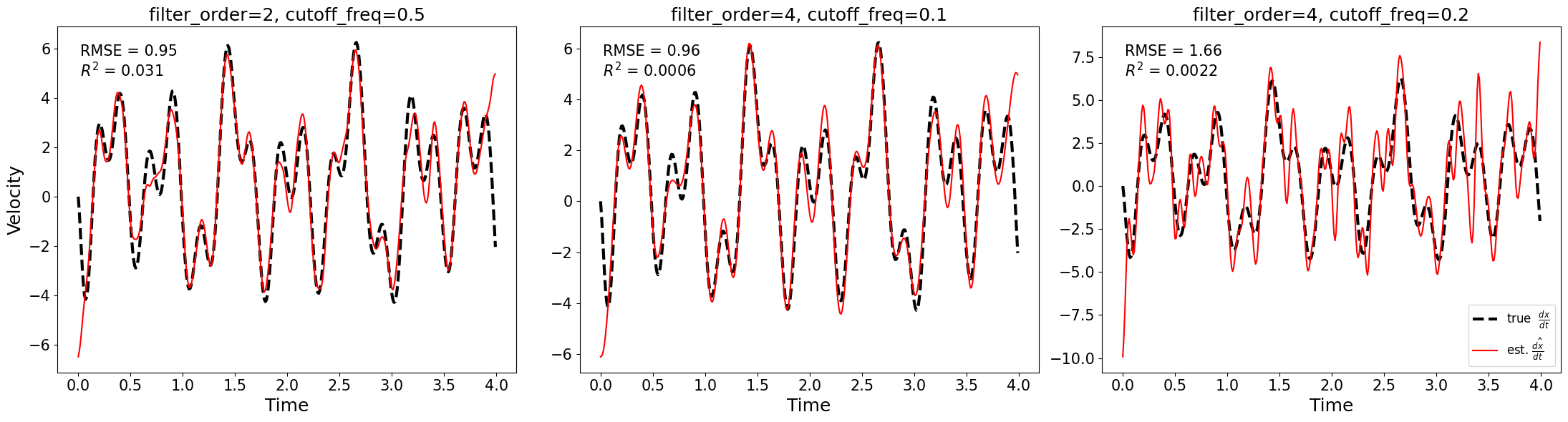}
  \vspace{-3mm}
  \caption{Derivative estimates for noisy cruise control data (\autoref{fig:cruise-control-sim}) found by Butterworth smoothing with three different choices of hyperparameters followed by second-order Finite Difference.}
\end{figure}

Instead of relying on such approximations, we can design and apply a true low-pass filter directly. A simple option is the Butterworth filter, with gain (complex transfer function magnitude) given by:
\begin{equation}\label{eqn:butter-transfer}
G(\omega) = |H(i\omega)| =  \frac{1}{\sqrt{1+\big(\frac{\omega}{\omega_c}\big)^{2n}}}
\end{equation}

\noindent where $\omega$ is frequency, $\omega_c$ is the ``cutoff" where power is halved, and $n$ is the filter order. This response, visualized in \autoref{fig:freq-responses}, is intended to stay as flat as possible in the passband, thereby introducing minimal distortion. In contrast to the global Fourier-Spectral Method (\autoref{sec:fourier}), which can be used as an ideal low-pass filter by zeroing out higher modes followed by inverse transforming, the Butterworth filter is local and can only achieve gradual cutoff. But as compensation, Butterworth filters are applicable to any signal, not only those that are periodic.

\subsection{Iterated Finite Difference}\label{sec:iteratedfinitedifference}

\begin{figure}[!t]\label{fig:second-order}
  \centering
  \includegraphics[width=0.8\textwidth]{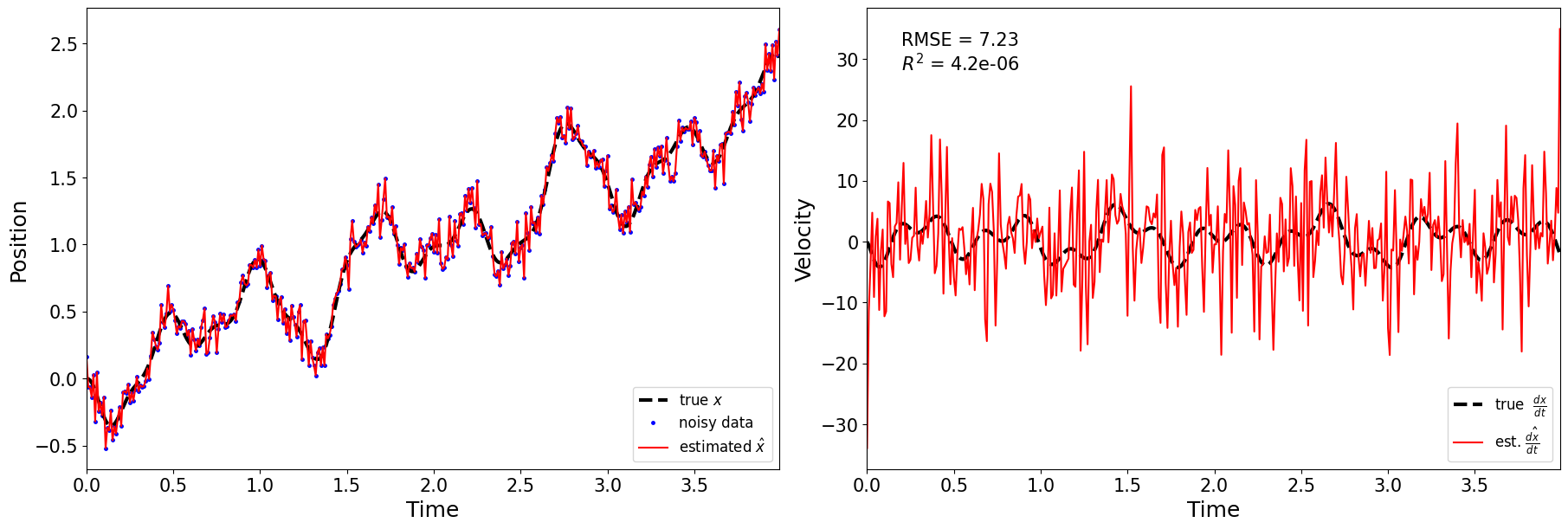}
  \vspace{-3mm}
  \caption{Second-Order Finite Difference on the cruise control data of \autoref{sec:kalman-example}. The unmodified method has no mechanism to estimate smoothed $\hbx$ from the data, resulting in exceptionally jagged derivative estimates. Domain edges are found with formulas that further intensify noise.}
\end{figure}

In the presence of noise, Finite Difference's division by powers of the (small) step size, as in the formulas of \autoref{ta:fd-schemes}, greatly amplifies deviations in noisy $\by$, resulting in poor derivative estimates of the underlying signal, $\bx$. In addition, forward- and backward-differencing equations tend to have large constants, which magnify noise and cause characteristic blowup at the domain edges, as seen in \autoref{fig:second-order}.

The first problem, division by $\Delta t$, can be counterbalanced by performing integration, which multiplies by $\Delta t$. These operations may appear to pointlessly cancel out; indeed the most naive scheme, first-order Finite Difference followed by cumulative summation, essentially just shifts the data:
\begin{align*}
\by = [y_0, y_1, y_2, y_3] \rightarrow \text{first-order-FD}\big(\by \big) &= [y_1 - y_0, y_2 - y_1, y_3 - y_2]/\Delta t\\
\rightarrow \text{cumsum}\big( \text{first-order-FD}\big(\by \big) \big) \cdot \Delta t &= [y_1, y_2, y_3] - y_0
\end{align*}

\noindent However, second-order center-difference followed by cumulative trapezoidal integration (trapz) introduces a factor of $\frac{1}{2}$. Ignoring the edges by \texttt{numpy}-slicing $[1\!:\!-1]$, we find:
\begin{align*}
\text{second-order-FD}\big([y_0, y_1, y_2, y_3] \big)[1\!:\!-1] &= [y_2 - y_0,\ y_3 - y_1]/2\Delta t\\
\rightarrow \text{trapz}\big( [y_2 - y_0,\ y_3 - y_1]/2\Delta t\big) &=\!\!
\begin{tikzpicture}[baseline=(filt.base)]
  \node (filt) {$[0,\ \frac{y_2 - y_0 + y_3 - y_1}{4}]$};
  \node[below left=0mm and -9mm of filt] (n1) {\footnotesize $n\!=\!1$};
  \node[below right=0mm and -12mm of filt] (n2) {\footnotesize $n\!=\!2$};
  \draw[->] (n1.center)+(0,1mm) -- +(-1mm, 3.5mm);
  \draw[->] (n2.center)+(0,1mm) -- +(-1mm, 3.5mm);
\end{tikzpicture}
\end{align*}

\begin{figure}[!t]\label{fig:iterated-fd-freq-response}
  \centering
  \includegraphics[width=0.7\textwidth]{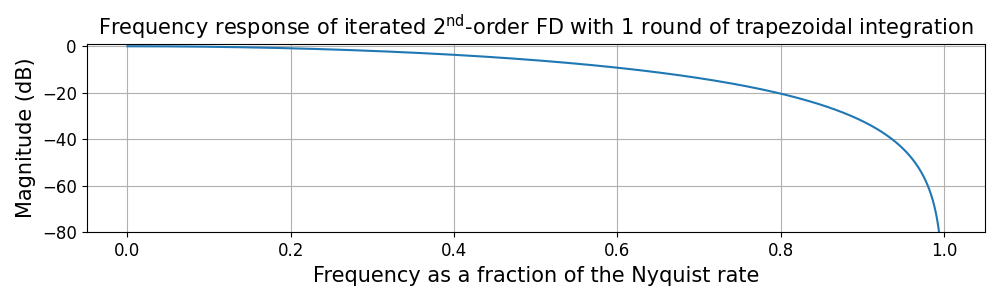}
  \vspace{-3mm}
  \caption{The frequency response of the IIR filter equivalent to one round of second-order center-differencing followed by trapezoidal integration. Iterating more times causes the curve to multiply itself and bend down more rapidly, reaching -80dB by about 1/3 Nyquist for 30 repetitions. The passband is also curved, not flat.}
\end{figure}

In the language of digital filters, this is effectively implementing the infinite impulse response (IIR) filter, $\hat{x}[n] = \hat{x}[n-1] + \frac{1}{4}(y[n+1] + y[n] - y[n-1] - y[n-2])$, where, in a stunning reversal of ordinary digital filtering notation, $y$ are inputs and $\hat{x}$ are outputs. This is yet another low-pass filter, with response shown in \autoref{fig:iterated-fd-freq-response}. Iterating the procedure is like passing through this filter multiple times.

Higher-order differencing schemes smooth across more samples, although the center-differencing coefficients of immediately-neighboring samples tend to become greater, concentrating energy more locally and necessitating more iterations to adequately smooth. It is also possible to use higher-order quadrature rules, but even Simpson's 1/3 rule, which integrates over $2\Delta t$ at a time with $\frac{1}{3}(y_{n-2} + 4y_{n-1} + y_n)$, contains a coefficient $>\!\!1$, which risks amplifying noise.

Edge blowup can be avoided by simply using lower-order differencing equations, but one can also solve the Finite Difference linear inverse problem (\autoref{eqn:fd-generalized}) with a spaced-out stencil, e.g.~$[0, 2, 4]$, to yield coefficients that do not exceed unity.

\begin{figure}[!t]
  \centering
  \includegraphics[width=0.99\textwidth]{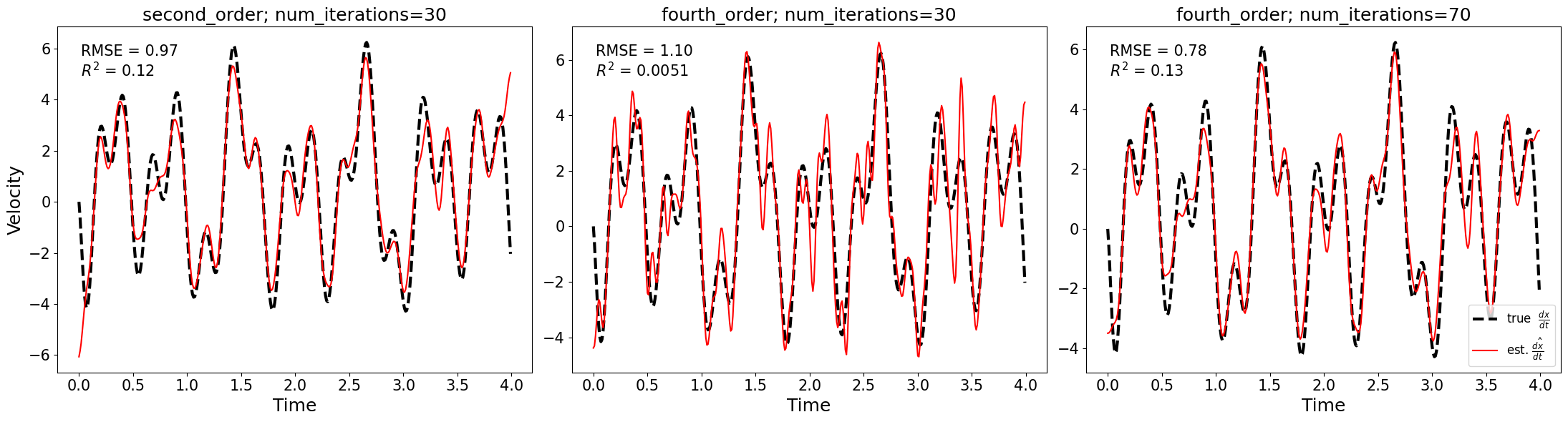}
  \vspace{-3mm}
  \caption{Iterated Finite Difference on noisy cruise control data (\autoref{fig:cruise-control-sim}) with three different choices of hyperparameters.}
\end{figure}

\subsection{Polynomial Fits}
\label{sec:polyfit}

This section covers methods that set up and solve regression problems to approximate data with smooth polynomial functions, which can then be differentiated analytically by power rule and evaluated at sample points.

\subsubsection{Sliding Window Polynomial Fits}\label{sec:sliding-poly}

A function $x(t_i) = c_0 + c_1t_i + ... + c_dt_i^d$, where $d$ is some degree, can be fit to data $\by$ by solving the linear inverse problem:
\begin{equation}\label{eqn:polynomial-least-squares}
\begin{bmatrix}
1 & t_1 & t_1^2 & \cdots & t_1^d \\
1 & t_2 & t_2^2 & \cdots & t_2^d \\
\vdots & \vdots & \vdots & \ddots & \vdots \\
1 & t_n & t_n^2 & \cdots & t_n^d \end{bmatrix}
\begin{bmatrix} c_0 \\ c_1 \\ \vdots \\ c_d \end{bmatrix} \approx \by
\end{equation}

If the Vandermonde matrix is square, the two sides of \autoref{eqn:polynomial-least-squares} can be made equal, and $x(t)$ goes \textit{through} all data points. To obtain a smoother function, we limit degree $d$, making the system overdetermined, and minimize the $\ell_2$ norm of the difference between the two sides.

\begin{figure}[!t]\label{fig:sliding-polyfits}
  \centering
  \includegraphics[width=0.6\textwidth]{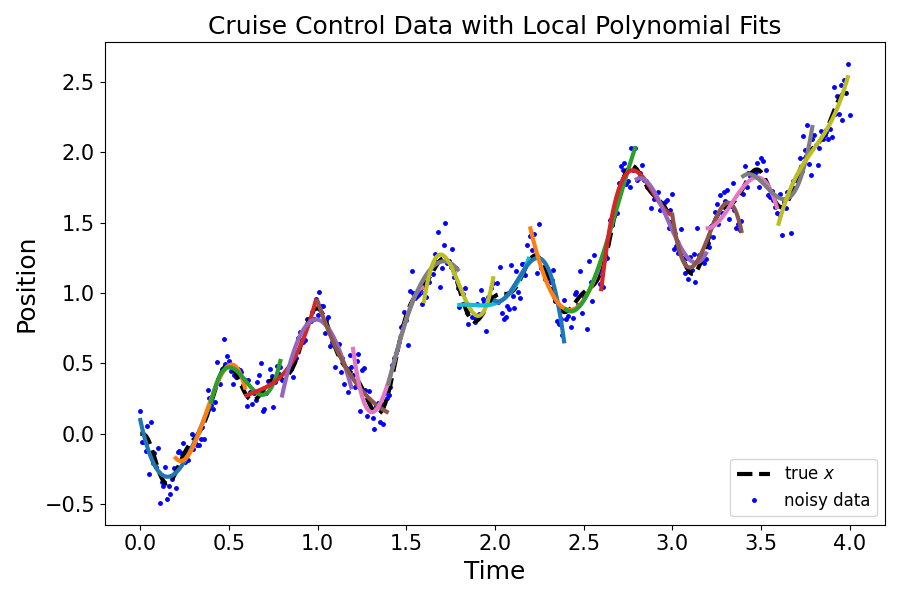}
  \vspace{-4mm}
  \caption{Polynomial fits, each with its own color, to batches of 40 samples with a stride of 20, on cruise control data from \autoref{sec:kalman-example}.}
\end{figure}

However, neither a low nor high choice of $d$ may be acceptable, because low-degree polynomials have very little expressive power while higher-degree fits tend to exhibit Runge's phenomenon (\autoref{fig:runge-phenomenon}). To circumvent this, we can restrict a fit's scope to some local subset of data and keep degree fairly low, because local subsets do not have as much opportunity for complex behavior. Covering the whole domain then requires many such fits, which can be thought of as happening inside a sliding window that moves by some stride, as illustrated in \autoref{fig:sliding-polyfits}.

\begin{figure}[!t]\label{fig:poly-cc}
  \centering
  \includegraphics[width=0.99\textwidth]{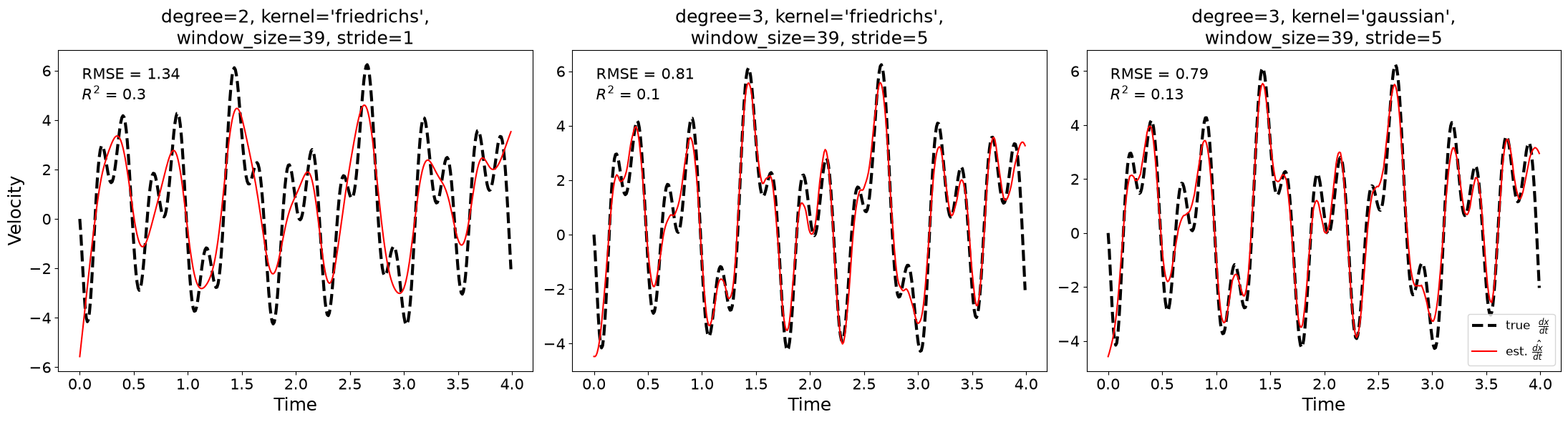}
  \vspace{-3mm}
  \caption{Sliding-window polynomial differentiation of noisy cruise control data (\autoref{fig:cruise-control-sim}) with three different choices of hyperparameters.}
\end{figure}

Samplings of these fits and their derivatives can then be combined by averaging their overlaps or, if a kernel is applied to weight fit samples, weighted averaging. Example outputs are shown in \autoref{fig:poly-cc}.

\subsubsection{Savitzky-Golay Filtering}
\label{sec:savitzky-golay}

In the special case samples are equispaced, no weighting kernel is desired, and each output is found as a single sample from its own fit (moving-polynomial fit, stride of 1), output values can be calculated without explicitly solving linear inverse problems, because the solution can be written as a fixed linear combination of window inputs~\cite{savitzky1964smoothing, schafer-savgol}.

This is seen by an argument very similar to the derivation of Finite Difference equations in \autoref{sec:finite-difference}. Specifically, the moving window implies a stencil of $S$ samples, $[s_0, s_1, ..., s_{S-1}]$, e.g.~$[-1, 0, 1]$ to sample 3 equispaced points centered on a point of interest. We then wish to find a $d^\text{th}$ order polynomial of the form $c_0 + c_1s + c_2s^2 + ... + c_ds^d$ that when sampled at the stencil-points approximates the data, $\by_\text{window}$, at the corresponding points. Paralleling Equations \ref{eqn:fd-generalized} and \ref{eqn:polynomial-least-squares}, this can be set up as:

$$\begin{bmatrix}
s_0^0 & \cdots & s_0^d\\
\vdots & \ddots & \vdots \\
s_{S-1}^0 & \cdots & s_{S-1}^d
\end{bmatrix} \begin{bmatrix} c_0 \\ \vdots \\ c_d \end{bmatrix} = \begin{bmatrix} \vert \\ \by_\text{window} \\ \vert\end{bmatrix}$$

\begin{figure}[!t]\label{fig:savgol-freq-response}
  \centering
  \includegraphics[width=0.7\textwidth]{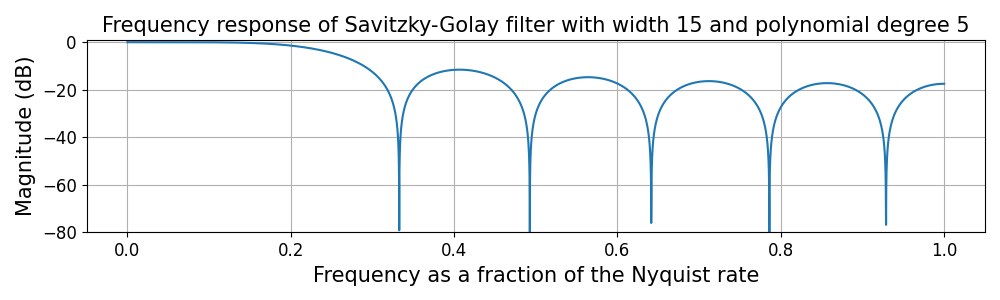}
  \vspace{-3mm}
  \caption{The frequency response of a Savitzky-Golay filter. It is very flat in the passband, which minimizes distortion~\cite{schafer-savgol}. When polynomial degree is 1, it is identical to the moving average filter (\autoref{fig:freq-responses}).}
\end{figure}

At the point of interest, $s=0$, the fit's value is given by $c_0$ alone, which can be found as the dot product of $\by_\text{window}$ with the first row of the pseudoinverse of the stencil Vandermonde matrix. Because said matrix depends only on the chosen stencil and polynomial degree $d$, the coefficients multiplying $\by_\text{window}$ are independent of the data and can be found offline. Moreover, the vector of coefficients constitutes a digital filter, with frequency response shown in \autoref{fig:savgol-freq-response}, which can be applied to data efficiently via convolution.\footnote{Data edges, where the stencil may not fit, must be treated specially, often by extending the data with repeats or mirroring prior to convolution or by explicitly calculating edge polynomial fits so they can be evaluated beyond single points.}

Savitzky-Golay can also deliver values of the derivative, by

$$\frac{d}{dt} (c_0 + c_1s + ... + c_ds^d) = \frac{d}{ds} (c_0 + c_1s + ... + c_ds^d) \frac{ds}{dt} = \frac{c_1}{\Delta t}$$

\noindent because $s = \frac{t}{\Delta t}$ is just a scaled variable, so $\frac{ds}{dt} = \frac{1}{\Delta t}$; and $c_1$, like $c_0$, is also merely a dot product of $\by_\text{window}$ with a vector of coefficients that can be found offline.

In practice, this method's results can be somewhat thrown around by the independence of implicit fits at subsequent data points, so it can be helpful to smooth the results in a postprocessing step by convolving with a Gaussian kernel.

\begin{figure}[!t]
  \centering
  \includegraphics[width=0.99\textwidth]{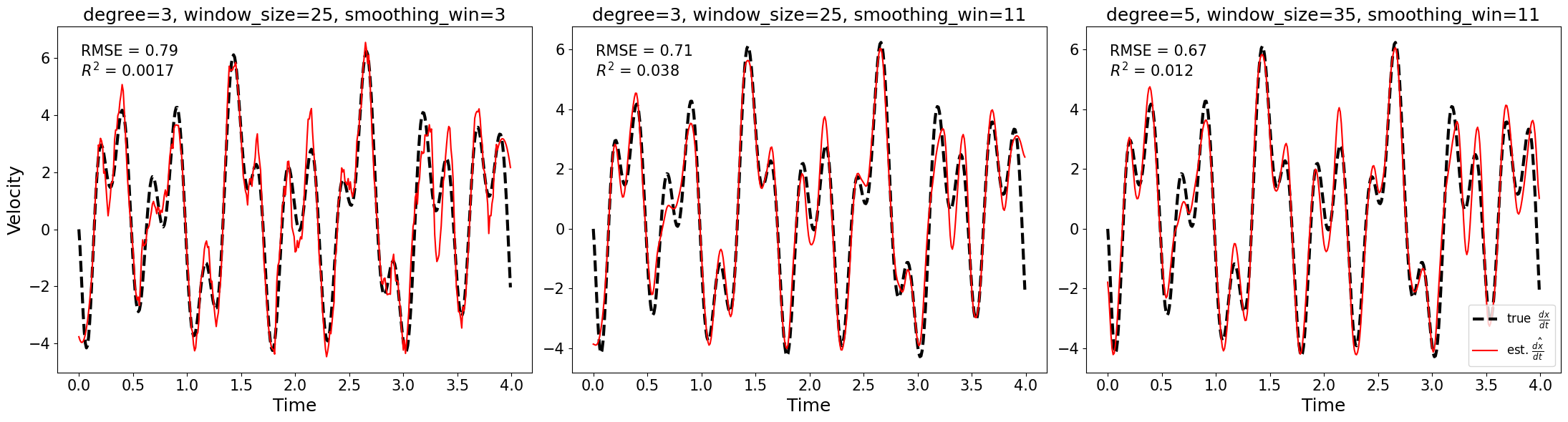}
  \vspace{-3mm}
  \caption{Savitzky-Golay differentiation followed by Gaussian smoothing of noisy cruise control data (\autoref{fig:cruise-control-sim}) with three different choices of hyperparameters.}
\end{figure}

\subsubsection{Spline Smoothing}
\label{sec:splinesmoothing}

Splines are piecewise polynomials, constructed to be continuous and to have continuous derivatives up to some order. Joins between the pieces are known as knots. With a sufficient number of knots, a spline can achieve a high degree of expressive power, but because each piecewise component only needs to consider local information, it can be fit with a low-degree polynomial, which averts Runge's phenomenon (\autoref{fig:runge-phenomenon}).

Each component of a spline can be given a definite formula, for example a cubic:

$$S_k(t) = c_{k,0} + c_{k,1} (t-t_k) + c_{k,2} (t-t_k)^2 + c_{k,3} (t-t_k)^3,\quad t\in[t_k,t_{k+1}]$$

\noindent By enforcing continuity of the function and its first derivative at the left and right boundary points of each spline component, we get the following four equations to determine the four unknown coefficients:

$$\left\{
\begin{aligned}
S_k(t_k)&=S_{k-1}(t_k)\\
S_k(t_{k+1})&=S_{k+1}(t_{k+1})\\
\dot{S}_k(t_k)&=\dot{S}_{k-1}(t_k)\\
\dot{S}_k(t_{k+1})&=\dot{S}_{k+1}(t_{k+1})
\end{aligned}
\right.$$

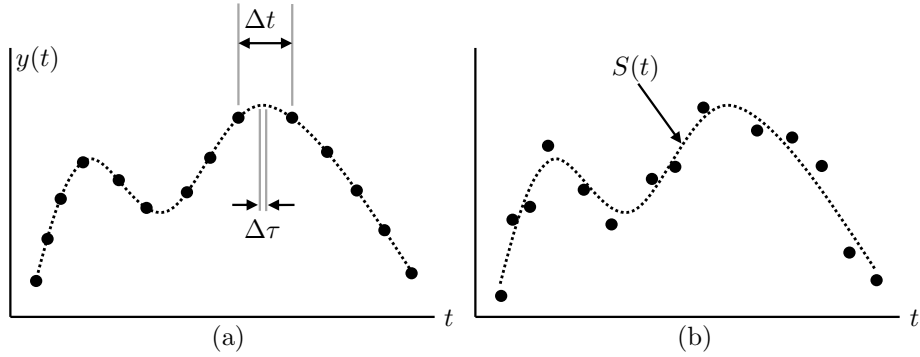
\begin{figure}[!t]\label{fig:spline_smoothing}
\centering
\begin{overpic}[width=0.75\textwidth]{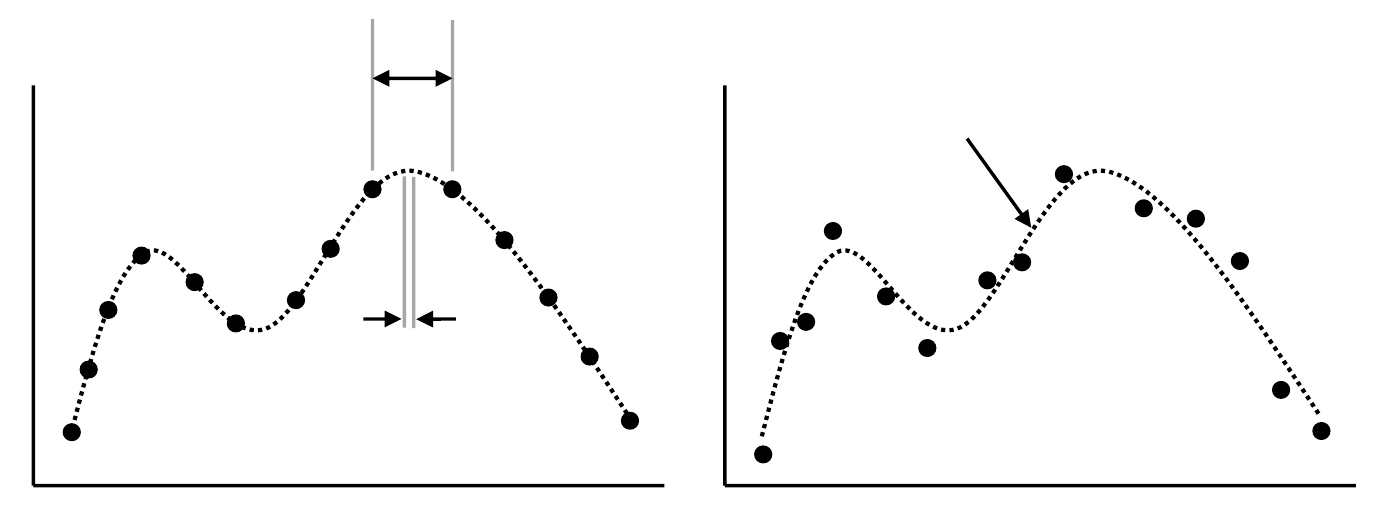}
\put(27.5,34){$\Delta t$}
\put(27.5,11.5){$\Delta \tau$}
\put(67,29){$S(t)$}
\put(49,2){$t$}
\put(99,2){$t$}
\put(3,30){$y(t)$}
\put(24,0){(a)}
\put(74,0){(b)}
\end{overpic}
\vspace{-1mm}
\caption{(a) Cubic spline fitting to noiseless data for refinement of a grid from $\Delta t$ to $\Delta \tau$. (b) Cubic spline smoothing of noisy data. In this case, the spline is not forced to go through the samples, rather a balance between data fitting and smoothness is sought.}
\end{figure}

\noindent At the left and right boundaries of the domain, a neighboring component is no longer available in one direction, so alternative constraints must be imposed, such as $\ddot{S}_0(t_0)=0$ and $\ddot{S}_{n-1}(t_n)=0$. All these equations can ultimately be written as a sparse linear inverse problem with target $\by$, similar to \autoref{eqn:spectral-linear-inverse}. If there are as many knots as sample points, the spline intersects every data point, which can be useful for refining coarsely-spaced, noiseless data, as shown in \autoref{fig:spline_smoothing}(a).

For noisy data, the spline-fitting problem can be modified so the solution passes between data points, as depicted in \autoref{fig:spline_smoothing}(b). There are a couple different ways to do this. The first~\cite{wahba1975smoothing,craven1978smoothing} is to add a term that penalizes the integral of the second derivative:
\begin{equation}\label{eqn:spline-optimization}
  \argmin_{\bth} \sum_{n=0}^{N-1} \big(y_n - S(t_n;\bth)\big)^2  + \lambda\!\! \int\limits_{t_0}^{t_{N-1}} \!\!\!\big(S''(u;\bth)\big)^2du
\end{equation}

\noindent where $S$ is the overall spline, $\bth$ represents all coefficients collected together, and $\lambda$ is a hyperparameter which balances fitting the data with keeping the second derivative small. The second possibility~\cite{dierckx1975} is to simply require the sum of squares of the best fit be within a bound, $s$, which corresponds to the minimization of \autoref{eqn:spline-optimization} with $\lambda = 0$ and an additional condition:
\begin{equation}\label{eqn:spline-scipy}
\text{s.t.}\quad \sum_{n=0}^{N-1} (y_n - S(t_n))^2 \leq s
\end{equation}

\noindent where larger $s$ allows for more wiggle room and therefore more smoothness. This will be either feasible or infeasible, depending on the number of knots.

\begin{figure}[!t]\label{fig:bsplines}
  \centering
  \includegraphics[width=0.99\textwidth]{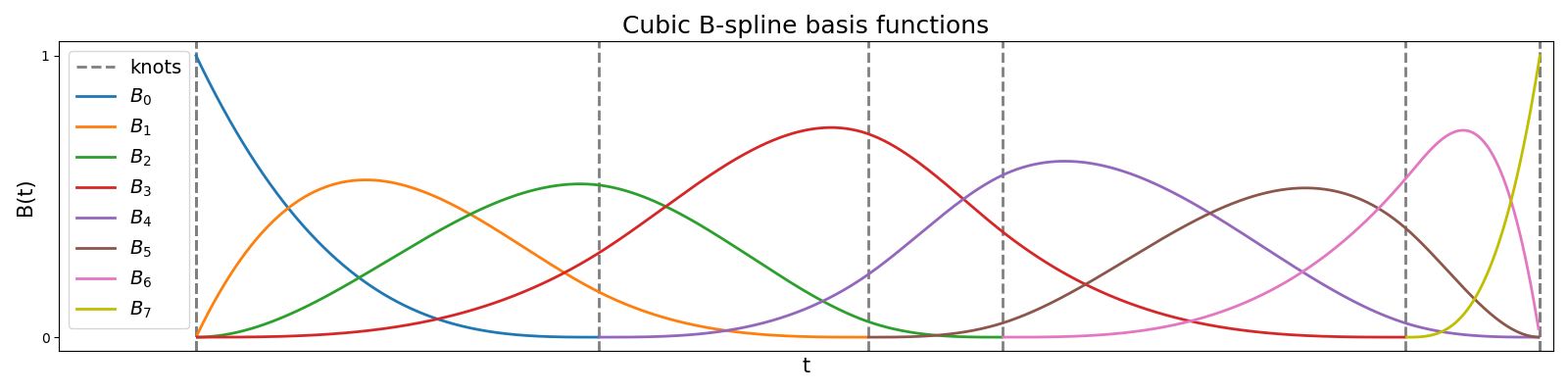}
  \vspace{-3mm}
  \caption{Visualization of cubic basis splines. Knot locations are shown as vertical dashed lines. Each basis function has compact support, i.e.~takes nonzero value only on at most degree\! +\! 1 between-knot intervals. Within each between-knot interval, the basis splines form a partition of unity, i.e.~sum to one.}
\end{figure}

One route to solve these involves representing $S$ as a sum of cubic B-splines, B for ``basis".
\begin{equation}\label{eqn:b-spline-reconstruction}
S(t) = \sum_{j=1}^m \alpha_j B_j(t)
\end{equation}

\noindent where $B_j(t)$ is the $j^\text{th}$ B-spline, $\alpha_j$ is a coefficient, and $m$ is the total number of basis functions, which depends on the number of knots and smoothness we wish to achieve. B-splines, visualized in \autoref{fig:bsplines}, are designed to be nonzero over short intervals, redolent of the basis functions from the Finite Element Method (\autoref{sec:finite-elements}).

\autoref{eqn:b-spline-reconstruction} can be substituted in \autoref{eqn:spline-optimization} to give:

$$\argmin_{\boldsymbol\alpha} \sum_{n=0}^{N-1} \Big(y_n - \sum_j^m \alpha_j B_j(t_n)\Big)^2 + \lambda\!\! \int\limits_{t_0}^{t_{N-1}} \!\!\! \Big( \sum_j^m \alpha_j B_j''(u) \Big)^2 du$$

\noindent where $t_n$ is a sample point, not a knot point as sometimes denoted when working with splines. This can be further written in matrix form:
\begin{equation}\label{eqn:spline-opt-vec-form}
\argmin_{\boldsymbol\alpha} \big\| \mathbf{y} - \bB \boldsymbol\alpha \big\|^2 + \lambda \, \boldsymbol\alpha^\mathrm{T} \bR \boldsymbol\alpha
\end{equation}

\noindent where $\bB \in \mathbb{R}^{N \times m}$ has entries $B_{nj} = B_j(t_n)$ and $\bR \in \mathbb{R}^{m \times m}$ is the roughness penalty matrix, with entries $R_{jk} = \int B_j''(\tau) B_k''(\tau)\ d\tau$.

Taking the gradient with respect to $\boldsymbol\alpha$ and setting equal to zero yields a linear inverse problem:
$$(\bB^\mathrm{T} \bB + \lambda \bR) \boldsymbol\alpha = \bB^\mathrm{T} \mathbf{y}$$

\noindent Due to the compactness of B-splines, only $O(d)$ can overlap any point in the domain, where $d$ is the degree of $S$, and therefore each row of $\bB$ contains only $O(d)$ nonzero elements, so $\bB^\mathrm{T} \bB$ can be found in $O(N d^2)$ instead of $O(Nm^2)$. Both $\bB^\mathrm{T} \bB$ and $\bR$ are \textit{banded}, i.e.~nonzero on only $O(d)$ diagonals, allowing the system to be solved efficiently in $O(md^2)$.

If using smoothness hyperparameter $s$ instead of $\lambda$ (\autoref{eqn:spline-scipy}), the second term of \autoref{eqn:spline-opt-vec-form} disappears, breaking regularization. However, smoothness can still be imposed by keeping the number of knots low, because a spline with fewer knots has less expressive power. So we can find a good choice of parameters by refining a global, low-degree polynomial fit by iteratively adding knots at locations of greatest disagreement between the spline and data.\footnote{This is the approach favored by \texttt{scipy}~\cite{scipy-smoothing-splines}.}

\begin{figure}[!t]
  \centering
  \includegraphics[width=0.99\textwidth]{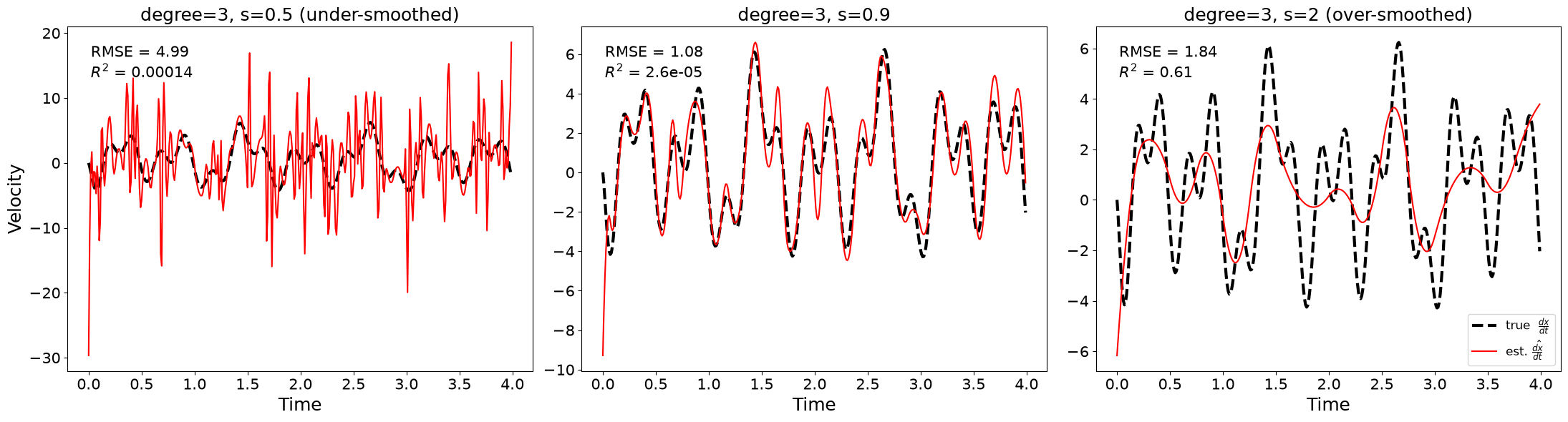}
  \vspace{-3mm}
  \caption{Spline differentiation of noisy cruise control data (\autoref{fig:cruise-control-sim}) with three different choices of hyperparameters.}
\end{figure}

\subsection{Basis Fits}
\label{sec:basis-fit}

This section covers methods that build function approximations from smooth basis functions, both global and local. B-splines, presented in the previous section, are an example of a local basis, so splines can also straddle this category.

\subsubsection{Fourier Spectral with Tricks}

\begin{figure}[!t]\label{fig:spectraldiff-hacks}
\centering
\includegraphics[width=0.99\textwidth]{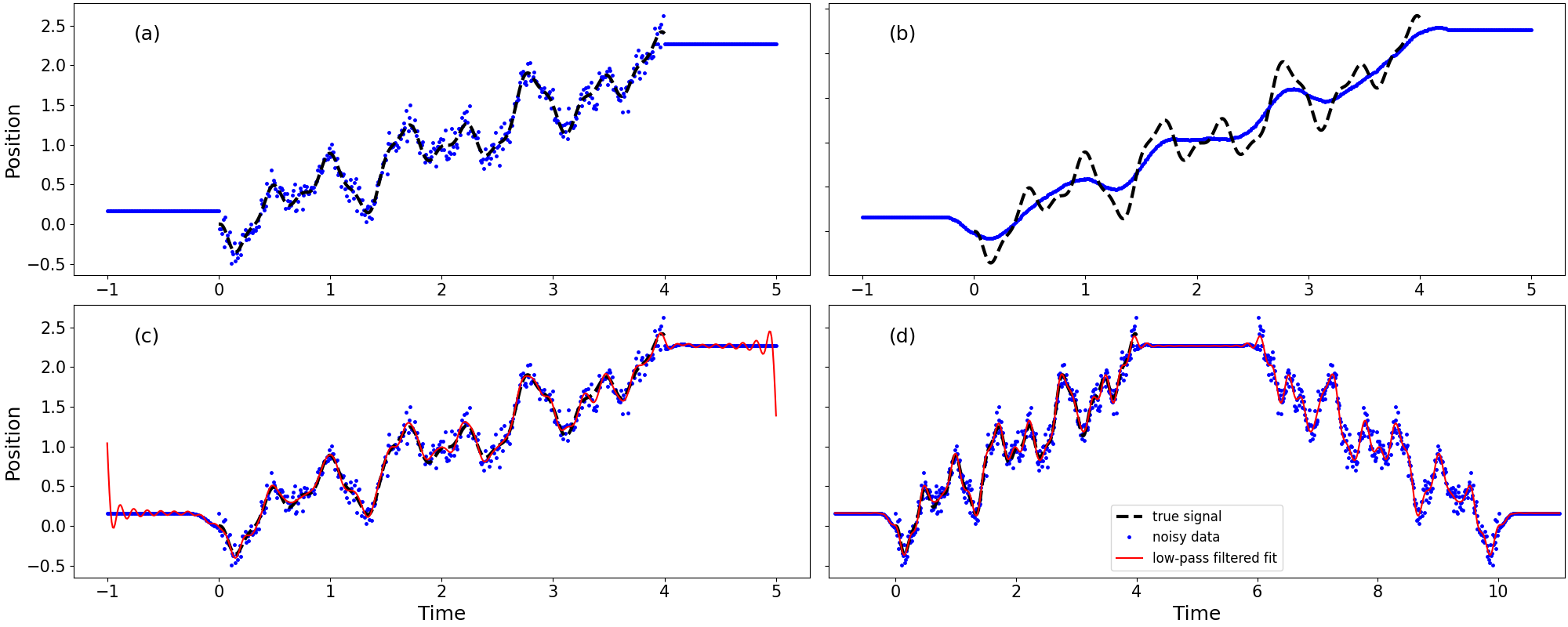}
\vspace{-2mm}
\caption{Preprocessing can make Fourier spectral filtering workable with aperiodic data. (a) Noisy data is padded with repeats of its end values. (b) The padded data is convolved with a moving average filter. (c) The original data is placed back in the center of the smoothed data from (b), and the combination is Fourier transformed. All modes except the lowest 50 are zeroed out, and the inverse transform is shown in red. Gibbs phenomenon is apparent but most extreme in the padded regions, which can be discarded. (d) The even extension of the data from (c) is Fourier transformed, filtered to the first 70 modes, and inverse transformed to produce the red curve. Gibbs phenomenon is much reduced.}
\end{figure}

Even though global Spectral Methods (\autoref{sec:spectral}) are best applied in the absence of noise, due to Gibbs phenomenon in the case of the Fourier basis (\autoref{sec:fourier}) and overfitting in the case of polynomial bases (\autoref{sec:chebyshev}), noise is fundamentally such a frequency-dependent phenomenon (\autoref{sec:noise}), and the FFT is so efficient, that fitting data with sinusoids is still attractive. However, to get around some of the method's weakness and work with arbitrary, likely-aperiodic data requires clever preprocessing.

The key is to recognize Gibbs phenomenon is caused by discontinuities and corners within a signal and at its edges (\autoref{fig:gibbs-phenomenon}), and is most profound \textit{near} these triggers. One strategy, then, is to take the signal's even extension (\autoref{fig:periodic-extensions}), which mirrors across an end, getting rid of discontinuity but resulting in a corner if the signal edges are not flat. It is also possible to take the odd extension, reflecting a signal through its endpoint, which preserves slopes but causes a greater discontinuity. Thankfully, we can account for this by linear trend removal, to put the signal's edges on the same level, easily adding the trend's effect back at the end of the algorithm. Another approach, which can be used in tandem with extension or independently, is to insulate the signal by padding its edges, being careful to smoothly transition between signal and pads so as not to subtly create more corners.

The preprocessed data can then be Fourier transformed, ideal-low-pass filtered by zeroing out higher modes, optionally differentiated by multiplying by appropriate $ik$ (\autoref{eqn:fourier-times-ik}), and inverse transformed. The smoothed signal or derivative is then found by slicing out part of the result. This whole process is illustrated with even extension and padding in \autoref{fig:spectraldiff-hacks}.

\begin{figure}[!t]
  \centering
  \includegraphics[width=0.99\textwidth]{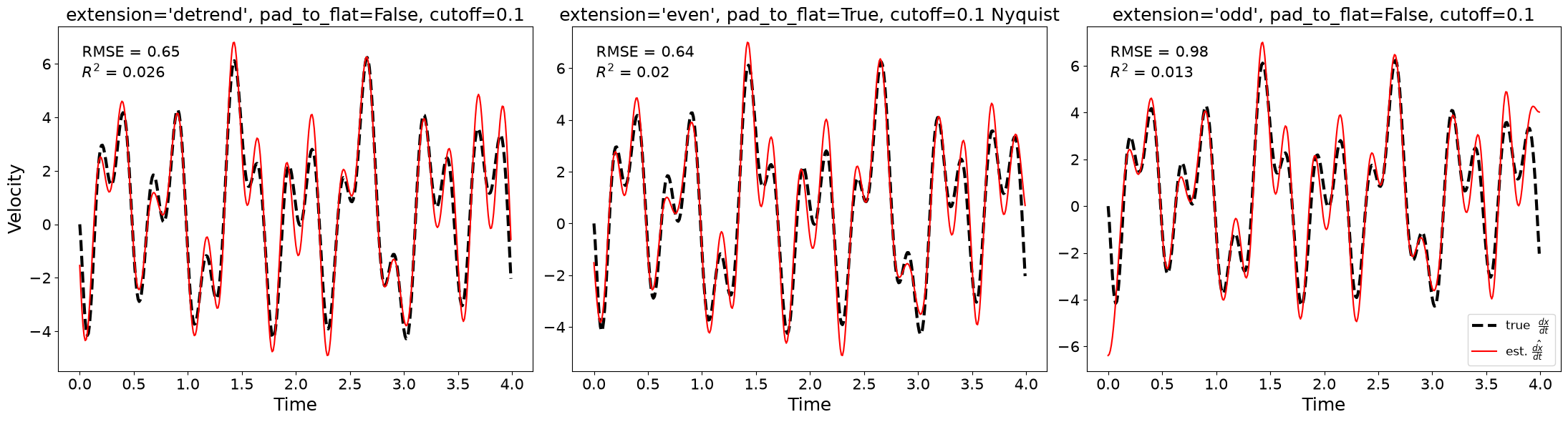}
  \vspace{-3mm}
  \caption{Spectral derivatives of noisy cruise control data (\autoref{fig:cruise-control-sim}) with three different preprocessing strategies.}
\end{figure}

\subsubsection{Radial Basis Functions}
\label{sec:rbf}

\begin{figure}[!t]\label{fig:rbf-A}
  \centering
  \includegraphics[width=0.8\textwidth]{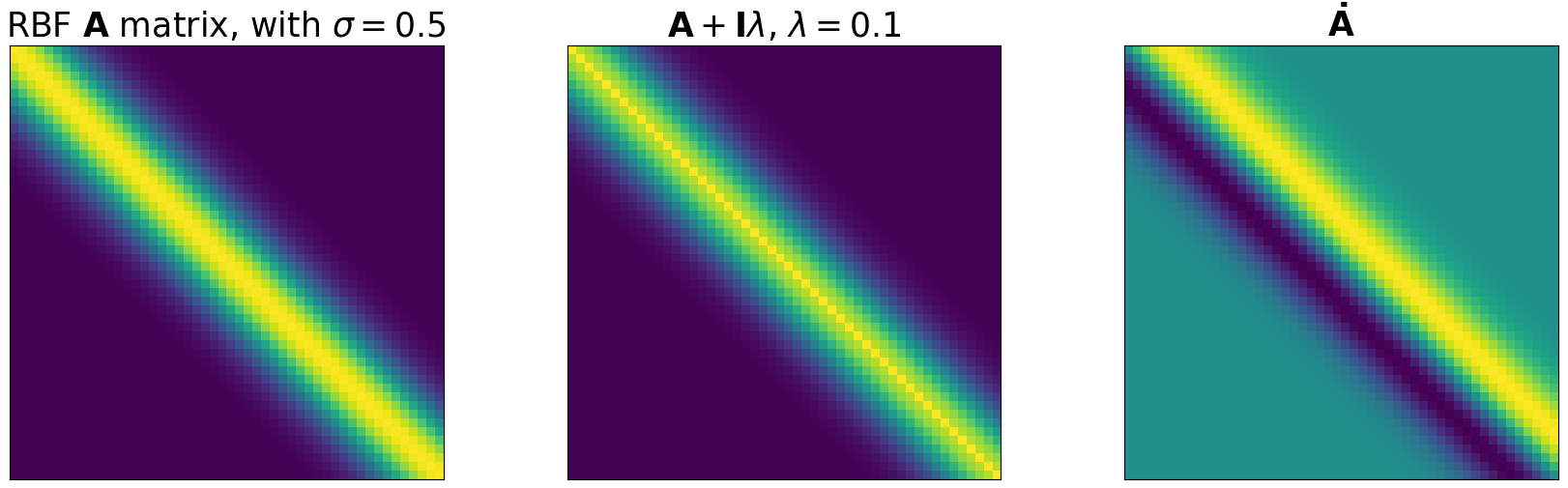}
  \vspace{-2mm}
  \caption{Visualizations of sparse, banded matrices used to fit radial basis functions to noisy data. Brighter colors indicate relatively higher values, and darker colors lower values. At left is a matrix filled with basis function samples, truncated when values get below 1e-4. This example has condition number about 1.2e9. Middle is a damped version of the same with smaller condition number, about 120. Right is a matrix filled with samples of the derivatives of basis functions.}
\end{figure}

Local functions with independent variable given by a radial distance, $r$, which evaluate to zero beyond some distance, $\rho$---such as truncated Gaussians, $e^{-r^2/(2\sigma^2)},\ r < \rho$---can also constitute a linearly independent basis, making the matrix in the linear inverse problem of \autoref{eqn:spectral-linear-inverse} banded, symmetric, and invertible in $O(N\rho^2)$. We call the matrix formed from columns of sampled basis functions $\bA$ and visualize an example in the first panel of \autoref{fig:rbf-A}. Unfortunately, $\bA$ is typically ill-conditioned, having eigenvalues very close to 0, which makes coefficient solutions highly sensitive to noise in the target, $\by$. Unlike the Fourier case where sensitive higher modes, associated with tiny eigenvalues of $\bA$, can be zeroed out, discarding any element of a local basis significantly deforms function reconstruction in some region, so we need to keep them all.

One possible way to handle this dilemma is Tikhonov Regularization~\cite{Tikhonov}, which attempts to keep coefficient values from blowing up by penalizing their $\ell_2$ norm:
$$\underset{\bc}{\text{minimize}}\; \big\|\bA\bc - \by\big\|_2^2 + \lambda\|\bc\|_2^2$$

\noindent where $\lambda$ is a chosen parameter. By taking the gradient w.r.t.~$\bc$ and setting equal to 0, the above is minimized by solving the linear inverse problem $(\bA^\mathrm{T}\bA + \lambda\mathbb{I})\bc = \bA^\mathrm{T}\by$.

Penalizing $\|\bc\|_2^2$ treats the basis functions as independent, but neighbors overlap: Coefficients of $+c$ and $-c$ on adjacent centers result in bumps that nearly cancel when the function is assembled. Since $\bA$ records basis function positions, not only values, $\|\bc\|_\bA^2 = \bc^\mathrm{T}\bA\bc$ penalizes the fit instead of the coefficient list, and minimizing $\big\|\bA\bc - \by\big\|_2^2 + \lambda\|\bc\|_\bA^2$ yields:
$$(\bA + \lambda\mathbb{I})\bc = \by$$
This effectively shifts all the eigenvalues of $\bA$ up away from 0, ``damping" the inverse problem. It also has the benefit of avoiding computation of $\bA^\mathrm{T}\bA$, which is not as sparse.

This latter scheme can also be thought of as manipulating the basis functions themselves, topping every ``hill" with a ``tower" to reach noisy data points. If the radial basis function quantifies how strongly values at nearby points move together, $\bA$ can also be interpreted as a covariance matrix, and $\lambda$ becomes the assumed noise variance, so $(\bA + \lambda\mathbb{I})\bc = \by$ is Gaussian process regression.

To obtain a smoothed estimate, we raze the towers and keep only the contributions of the original basis functions, using the classic reconstruction sum (\autoref{eqn:reconstruction-sum}):

$$\hat{x}(t) = \sum_{k=0}^{N-1}c_k \xi_k(t),\quad \xi_k(t) = \begin{cases} e^{-\frac{(t_k-t)^2}{2\sigma^2}} & |t_k-t| < \rho\\ 0 & \text{otherwise} \end{cases}$$

\noindent or in vector form: $\hbx = \bA\bc$. From this point, differentiation is extremely easy, because the differential operator is linear and can therefore be moved inside the sum: Just perform the reconstruction with $\dot\xi_k$, or in vector form, fill $\mathbf{\dot{A}}$ with samples from basis functions' derivatives, as visualized in the third panel of \autoref{fig:rbf-A}, and find $\hdbx = \mathbf{\dot{A}}\bc$.

\begin{figure}[!t]
  \centering
  \includegraphics[width=0.99\textwidth]{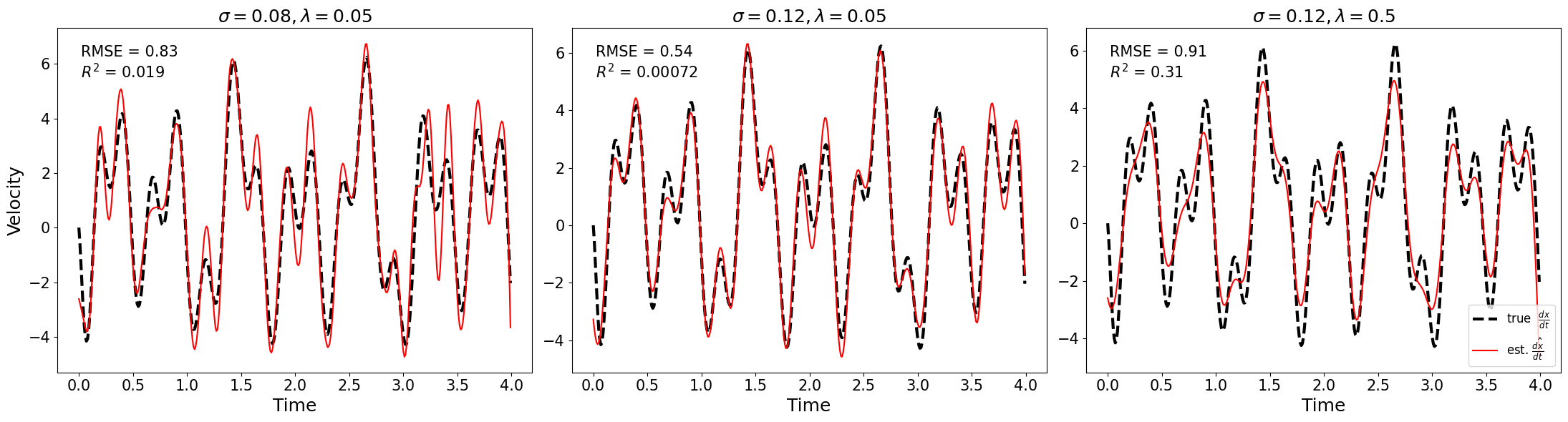}
  \vspace{-3mm}
  \caption{Radial basis function fit derivatives of noisy cruise control data (\autoref{fig:cruise-control-sim}) with three different choices of hyperparameters.}
\end{figure}

\subsubsection{Wavelets} It is also possible to decompose a function with local \textit{scaling functions} that have a self-similarity relation:
$$x(t) \approx \sum_{n=0}^{N-1} c_n\, \varphi(\tfrac{t}{\Delta t} - n),\quad\varphi(t) = \sqrt2\, \sum_{k=0}^{L-1} h_k\, \varphi(2t-k)$$ 
where the $\{h_k\}$ constitute a set of $L$ coefficients. That is, $\varphi$ can be composed by integer-offset, half-width, scaled versions of itself. Decomposing a signal this way relates coefficients at the finest resolution to those at longer distances, in theory capturing multi-scale structure. The factor of $\sqrt{2}$ keeps the power contribution from all levels of this recursion consistent, since halving a function's width halves its squared $L_2$ norm: $\int \varphi^2(2t-k)dt = \tfrac12\int\varphi^2(u)du$.

\begin{figure}[!t]\label{fig:scaling_function_refinement}
  \centering
  \includegraphics[width=0.99\textwidth]{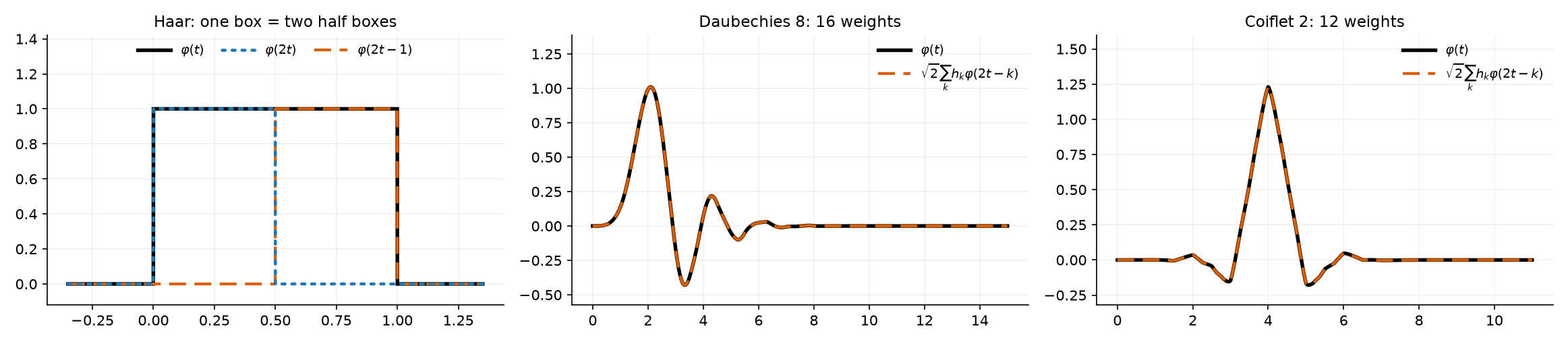}
  \vspace{-3mm}
  \caption{Several scaling functions and their refinements. The Haar family's scaling function can be composed of only two half-width copies but is discontinuous and unsuitable as a basis for differentiation. More complicated shapes like the Daubechies and Coiflet families~\cite{daubechies1992} are composed of more copies.}
\end{figure}

\autoref{fig:scaling_function_refinement} shows the scaling functions for three different wavelet types. Note these are not yet the wavelet itself. For most types, the scaling functions do not actually have a closed analytic form, rather they are numerically computed at sample locations down to arbitrary refinement. Consider sampling only at integers:
$$\begin{aligned}
\varphi[n] &= \sqrt2\, \sum_{k=0}^{L-1} h_{k}\,\varphi[2n - k] = \sqrt2\, \sum_{k=0}^{L-1} h_{2n-k}\,\varphi[k],\quad n, k \in \mathbb{Z} \\
&= (\mathbf M\boldsymbol\varphi)_n,\quad \text{where } \mathbf M \in \mathbb{R}^{L \times L},\ \mathbf M[n,k] = \begin{cases}\sqrt2\,h_{2n-k} & 2n\!-\!k \in [0, L\!-\!1]\\ 0 & \text{otherwise}\end{cases}
\end{aligned}$$
where square brackets indicate discrete sampling. That is, the basis function's values at integers, collected in vector $\boldsymbol\varphi$, is an eigenvector of matrix $\mathbf M$~\cite{strang1996}, which depends on weights $\{h_k\}$, with eigenvalue $1$. Likewise, we can differentiate both sides to find $\boldsymbol{\dot\varphi} = 2\mathbf{M} \boldsymbol{\dot\varphi}$, which makes the sampled derivative also an eigenvector, this time with eigenvalue $\frac{1}{2}$. Since eigenvectors can be of arbitrary magnitude, by convention we normalize so $\sum \boldsymbol\varphi = 1$, equivalent to $\int_a^b\varphi(x)dx = 1$, where $[a, b]$ is $\varphi$'s compact support, because, using the refinement relation (in practice, repeatedly convolving with $\sqrt{2}\{h_k\}$ spaced out with ever more 0s), finer samplings of $\varphi$ have the same spacing-weighted (i.e.~Riemann) sum as coarser ones. By a similar argument, we require $\sum \boldsymbol{\dot \varphi} = \int_a^b \dot \varphi(x) dx = \varphi(x)\big|_a^b = 0$ because $\varphi$ starts and ends at 0, and:
$$\mathbf{k} \cdot \boldsymbol{\dot \varphi} = \sum_{k=0}^{L-1} k\,\boldsymbol{\dot \varphi}_k = \int\limits_a^b x\, \dot \varphi(x) dx = \underbrace{x\,\varphi(x)\big|_a^b}_{\textstyle 0} - \int\limits_a^b\varphi(x)dx = -1$$
We only ever need to compute $\boldsymbol\varphi = \varphi[n]$ and $\boldsymbol{\dot\varphi} = \dot\varphi[n]$ without refinement, because $\Delta t$ between equispaced samples can be cleanly mapped to the space between integers.

To obtain nondegenerate solutions to the above ($\boldsymbol\varphi, \boldsymbol{\dot\varphi} \not\equiv 0$), two conditions must be imposed on $\{h_k\}$~\cite{strang1996}: (1) \textit{Accuracy} of order $p$, so-called because it means integer shifts of $\varphi$ can perfectly reconstruct polynomials of degree $< p$, satisfied when eigenvalues $1, \frac{1}{2}, \dots, 2^{-(p-1)}$ are in the spectrum of $\mathbf M$,\footnote{$\mathbf M$ also has eigenvalues aside from this ladder, with the largest giving an upper bound on $\varphi$'s smoothness as measured by the H\"older exponent, $\alpha \leq -\log_2|\lambda|$~\cite{daubechies1992}. For some types, notably \texttt{db2}, \texttt{sym2}, and \texttt{coif1}, an eigenvalue $\in (\frac{1}{2}, 1)$ warns that $\varphi$ is \textit{rough} enough to be technically non-differentiable almost everywhere, very much like the fractal structure of continuous Brownian motion. Counter-intuitively, these bases can still perfectly recapitulate degree $<p$ polynomials, because the roughness cancels when assembled. The $\lambda=\frac{1}{2}$ eigenvector remains an exact fixed point of $\mathbf M$ on the integers, but refining $\boldsymbol{\dot\varphi}$ in between diverges. In these cases we are therefore differentiating not the basis expansion, but the smooth signal its shifts reproduce, which can seem like hair-splitting since those two coincide in the regular, $\varphi \in C^1$ cases.} which is equivalent to requiring $\sum_k h_k = \sqrt2$ and $\sum_k (-1)^k k^m h_k = 0$ for $m \in \{0,\dots,p-1\}$. (2) \textit{Convergence} of the refinement recursion, satisfied when $\lambda_0 = 1$ is simple (lone) and every subsequent eigenvalue is strictly inside the unit disk.\footnote{$\mathbf M$ is sufficient for convergence on the integers, but convergence of $\varphi$ to a function in $L_2$ space requires these eigenvalue checks on the ``transition matrix", formed similarly but from autocorrelations $r_k = \sum_j h_j h_{j+k}$~\cite{strang1996}.}

\begin{figure}[!t]\label{fig:wavelet_filters}
  \centering
  \includegraphics[width=0.65\textwidth]{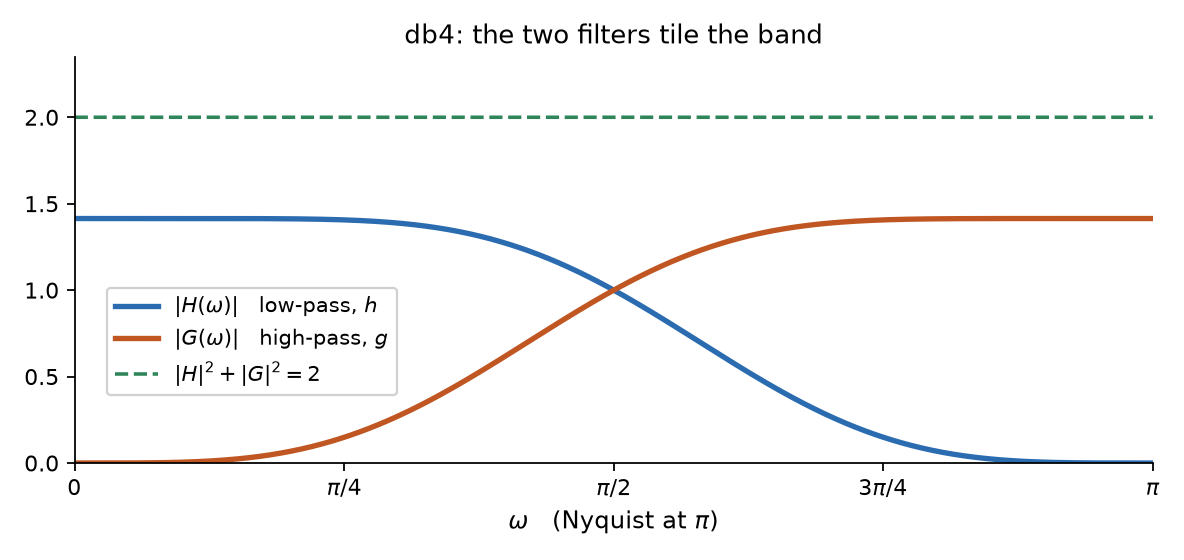}
  \vspace{-3mm}
  \caption{Magnitude plot of complementary low-pass and high-pass filters for the Daubechies type with $p=4$. They are mirrors of one another across $\frac{\pi}{2}$ and have total power summing to constant 2, which becomes 1 when downsampling reduces power by half.}
\end{figure}

We can view $\mathbf h = \{h_k\}$ as a low-pass filter: frequency transform $H(\omega) = \sum_k h_k e^{-i\omega k} \rightarrow H(0) = \sqrt{2}$, and $H(\pi) = 0$ with multiplicity $p$ from the accuracy requirements above. To capture both high and low frequency information from samples $\by$, we can apply the theory of \textit{filter banks}, where a signal is deconstructed and losslessly reconstructed using half the samples to represent each band~\cite{strang1996, mallat1989, smith1984}. By frequency-shifting, conjugating, and delaying $H$ (\autoref{fig:wavelet_filters}) we obtain $G(\omega) = e^{-i(\omega + \pi)(L-1)}\sum_k h_k e^{+i(\omega+\pi)k} \rightarrow g_k = (-1)^k h_{L-1-k}$, a high-pass filter with complementary power to $H$ between 0 and Nyquist as well as conjugate, delayed phase, specifically designed so aliases from downsampling exactly cancel when channels are recombined.

And now at last we meet the wavelet itself, which is formed by adding scaling functions according to $\mathbf g$ rather than $\mathbf h$:
$$\psi(t) = \sqrt{2}\, \sum_{k=0}^{L-1}g_k\varphi(2t-k)$$
\noindent This construction has $p$ \textit{vanishing moments}, meaning $\int_{-\infty}^\infty t^m \psi(t) dt= \langle t^m, \psi(t)\rangle = 0\ \forall\ m \in \{0,\dots,p-1\}$. That is, the same degree $< p$ polynomials $\varphi$ can reproduce perfectly are exactly annihilated in a basis of $\psi$,\footnote{even though $\psi$ is in turn built from a basis of $\varphi$} another way to see they carry complementary information.

The idea behind the \textit{wavelet transform}, then, is to repeatedly separate out higher-order details at many scales by encoding more and more in a wavelet basis, since wavelets are blind to low-order, smooth information. The Discrete Wavelet Transform (DWT) proceeds with:
$$c_n^{(j+1)} = \sum_k h_{k-2n} c_k^{(j)},\quad d_n^{(j+1)} = \sum_k g_{k-2n}c_k^{(j)}$$
where $j$ represents subsequent levels, and $\{d_n\}$ are ``detail coefficients". In linear algebra, this looks something like:
$$\begin{bmatrix}
\bc^{(j+1)}\\\mathbf d^{(j+1)}\end{bmatrix}
=\underbrace{\begin{bmatrix}
\textcolor{gray}{h_0} & \textcolor{gray}{h_1} & h_2 & h_3 & & & &\\
& & h_0 & h_1 & h_2 & h_3 & & \\
& & & & h_0 & h_1 & h_2 & h_3 \\
& & & & & \ddots & \ddots & \ddots \\
\textcolor{gray}{g_0} & \textcolor{gray}{g_1} & g_2 & g_3 & & & & \\
& & g_0 & g_1 & g_2 & g_3 & & \\
& & & & g_0 & g_1 & g_2 & g_3 \\
& & & & & \ddots & \ddots & \ddots
\end{bmatrix}}_{\textstyle \mathbf S} \underbrace{\begin{bmatrix} \textcolor{gray}{\text{ext}} \\\bc^{(j)} \\ \textcolor{gray}{\text{ext}}\end{bmatrix}}_{\textstyle \mathbf{E \, c}^{(j)}}$$
where here filter length $L = 4$ and grey denotes entries that extend $\mathbf c^{(j)}$, so the full-overlap (``valid" mode in \texttt{numpy}) convolution gives $\bc^{(j+1)}$ and $\mathbf d^{(j+1)}$ with $\lfloor (N^{(j)}+L-1)/2 \rfloor$ samples each, with $N^{(j)}$ the length of $\mathbf c^{(j)}$. As we continue down levels, the extension becomes a larger percentage of the input, and boundary effects worsen.

The inverse transform requires the combined extension, $\mathbf E$, and analysis, $\mathbf S$, operations to have a left inverse, i.e.~$\mathbf{\tilde S S E} = \mathbb{I}_{N^{(j)}}$:
$$\mathbf c^{(j)} = \mathbf{\tilde S} \begin{bmatrix} \mathbf c^{(j+1)} \\ \mathbf d^{(j+1)}\end{bmatrix}$$
The condition for $\mathbf{\tilde S}$ to exist is $\sum_k h_k \tilde h_{k+2n} = \delta_n$ for some \textit{synthesis} filter $\mathbf{\tilde h}$, a relationship called biorthogonality.\footnote{For biorthogonal wavelet types, high-pass filter $\mathbf g$ is built from $\mathbf{\tilde h}$, and $\mathbf{\tilde g}$ from $\mathbf h$.} When $\mathbf{\tilde h} = \mathbf h$ the transform is orthogonal, which preserves lengths and therefore keeps the energy in coefficients equal to energy in the signal (Parseval's Theorem applied to wavelets).

\begin{figure}[!t]\label{fig:wavelet_coefficients}
  \centering
  \begin{tikzpicture}[x=12cm,y=0.5cm]
    \draw (0,0) -- (1,0);
    \foreach \pos/\lbl in
      {0/{$0$}, 0.125/{$\frac{\pi}{8}$}, 0.25/{$\frac{\pi}{4}$}, 0.5/{$\frac{\pi}{2}$}, 1/{$\pi$}}
      {\draw (\pos,0.15) -- (\pos,-0.15); \node[below] at (\pos,-0.15) {\scriptsize \lbl};}
    \foreach \lo/\hi/\lbl in 
      {0/0.125/{$\mathbf c^{(3)}$}, 0.125/0.25/{$\mathbf d^{(3)}$}, 0.25/0.5/{$\mathbf d^{(2)}$}, 0.5/1/{$\mathbf d^{(1)}$}}{
      \draw[fill=blue!7] (\lo,0.4) rectangle (\hi,1.5);
      \node at ({(\lo+\hi)/2},1) {\lbl};}
    \node[below] at (0.53,-0.78) {\scriptsize angular frequency, with Nyquist ($f = \frac{f_s}{2} = \frac{1}{2\Delta t}$) at $\omega = \frac{2\pi f}{f_s} = \pi$};
  \end{tikzpicture}
  \vspace{-3mm}
  \caption{Depiction of DWT coefficients' dominant frequency bands, all drawn in terms of the top level's sampling frequency. At the first split, $\mathbf d^{(1)}$ holds high-pass information, primarily over $[\frac{\pi}{2}, \pi]$. Decimation effectively stretches $[0,\pi]$ onto $[0,2\pi]$ by making lower frequencies appear twice as fast per sample index, then aliases by reflecting the $[\pi, 2\pi]$ portion around Nyquist at $\pi$. Blind to this, the next level spans its own $[0, \pi]$, each covering half the bandwidth of its predecessor. Recursing on the low-pass side results in a telescoping dyadic structure, with each further decimation approximately halving the number of coefficients in play, so that all bands together are described in $O(N)$.}
\end{figure}

Each level of the DWT puts half the $\mathbf c$ coefficients in a wavelet basis and then recurses only on the low-pass branch:
$$\begin{array}{c@{\;}c@{\;}c@{\;}c@{\;}c@{\,}c@{\;}c@{\;}c}
\bc^{(0)} & \longrightarrow & \bc^{(1)} & \longrightarrow & \bc^{(2)} & \cdots& \longrightarrow & \bc^{(J)}\\
& \kern-4mm\searrow & & \kern-4mm\searrow & &  & \kern-4mm\searrow  & \\[-2pt]
& & \kern-5mm\mathbf d^{(1)} & & \kern-5mm\mathbf d^{(2)} & &  & \kern-5mm\mathbf d^{(J)}
\end{array}$$
where $J$ is the total number of repetitions. The leaves of this tree can be collected as in \autoref{fig:wavelet_coefficients} and together cover all information from the original signal.

We reach a purely structural stop condition (max $J$) based on the lengths of the signal and filter, $N$ and $L$, because eventually the convolution of $\mathbf h$ with a very small set of coefficients is overwhelmed by boundary effects. But we can choose a more intelligent stopping point based on coefficient magnitude across scales: Because we expect most signal energy to be in the low band (\autoref{sec:noise}) and expect energy from white noise to be scattered through all the $\mathbf d$, we can walk finest to coarsest looking for where wavelet coefficients start to increase sharply, suggesting they are beginning to run up against signal, and back off.

At the bottom of the recursion, we can clean up the signal by thresholding $\mathbf d$, knocking some higher-order, noise-like components out of the reconstruction. We use the Donoho-Johnstone universal\footnote{Universal in the sense that the expected maximum of $N$ standard normal draws approaches $\sqrt{2\ln N}$, so the chance that noise alone clears this threshold tends to zero, and what survives is very likely signal.} hard threshold~\cite{donoho1994}:
$$\mathbf{\hat d} = \begin{cases} \mathbf d & |\mathbf d| > \lambda_{DJ}, \qquad \lambda_{DJ} = |\eta|\sqrt{2\ln N}\\[2pt] 0 & \text{otherwise}\end{cases}$$
where $|\eta|$ is a noise scale, which can be estimated or searched over, and $N$ is the full signal length.

In practice, the DWT with thresholding is sensitive to starting index, because scaling functions and wavelets have specific powers-of-2 offsets, so a coefficient that clears $\lambda_{DJ}$ at one alignment may fall below it at the next. The remedy is to average reconstructed $\mathbf{\hat c}^{(0)}$ over $2^J$ sequential start points, since alignments repeat over this period, an approach known as ``cycle-spinning"\footnote{not to be confused with my favorite method of locomotion}~\cite{coifman1995}.

With cleaned-up coefficients in hand, we can finally do just as with Radial Basis Functions (\autoref{sec:rbf}) and sample the smooth function and derivative estimate with:
$$\hbx = \bA \mathbf{\hat c}, \qquad \hdbx = \mathbf{\dot{A}}\mathbf{\hat c}$$
where $\bA$ is filled with shifts of $\boldsymbol\varphi$ in its columns and $\mathbf{\dot{A}}$ is filled with shifts of $\boldsymbol{\dot\varphi}/\Delta t$, division because $\frac{dx}{dt} = \frac{1}{\Delta t} \sum_n c_n \dot \varphi(\frac{t}{\Delta t} -n)$. Picture it like \autoref{fig:rbf-A} but with more interesting shapes. The matrices are sparse and Toeplitz, so the products are equivalent to convolutions: $\hbx = \texttt{conv}(\mathbf{\hat c},\ \boldsymbol\varphi)$ and $\hdbx = \texttt{conv}(\mathbf{\hat c},\ \boldsymbol{\dot\varphi}/\Delta t)$.

One last wrinkle remains: We have been working with basis coefficients from the start but never specified how to get them. It is common practice to simply set $\bc^{(0)} = \by$, known as the ``wavelet crime"~\cite{strang1996}. To commit it is to describe an interpolant that \textit{does not pass through the data}, throwing off derivative estimates substantially.\footnote{For denoising purposes alone, it is often possible to get away with this crime, because the basis still concentrates low-order information in few coefficients and captures high-order wiggles with many wavelets, so separation, thresholding, and reconstruction can still yield cleaner point-values.} Instead we really should solve $\bA\bc^{(0)} = \by$, which can be done in frequency domain:
$$\by = \boldsymbol\varphi\, *\, \bc \overset{\mathcal F}{\longrightarrow} \hat{y}(\omega) = \hat\varphi(\omega) \hat c(\omega) \rightarrow \hat c(\omega) = \frac{\hat y(\omega)}{\hat\varphi(\omega)} \overset{\mathcal{F}^{-1}}{\longrightarrow} \bc = \by\, *\, \mathcal F^{-1}\Big\{\frac{1}{\mathcal F\{ \boldsymbol \varphi\}}\Big\}$$
where $*$ denotes convolution and $\mathcal F$ is the Fourier transform. The inverse response is not compactly supported, but it decays rapidly, so we can achieve sufficiently accurate answers by truncating where its values become tiny. This means we can properly solve for $\bc$ with a reusable FIR filter depending only on the wavelet type.

\begin{figure}[!t]
  \centering
  \includegraphics[width=0.99\textwidth]{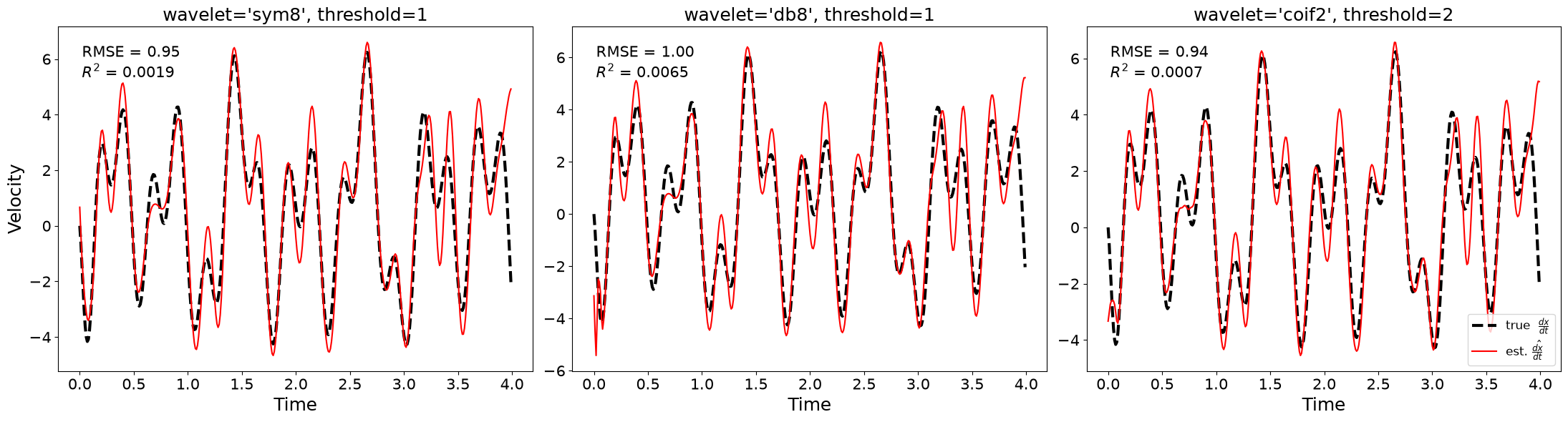}
  \vspace{-3mm}
  \caption{Derivatives in a scaling function basis, cleaned by thresholding in a wavelet basis, on noisy cruise control data (\autoref{fig:cruise-control-sim}).}
\end{figure}

\subsection{Total Variation Regularization}
\label{sec:tvr}

(TVR) Another way to think about noise in discrete series is that it causes extra disagreement between subsequent values, beyond what is necessary to trace the underlying true signal. We can measure the average distance between consecutive samples of a vector $\bv$ with normalized Total Variation:

$$\text{TV}(\bv) = \frac{1}{N}\big\|\bv_{0:N-1} - \bv_{1:N}\big\|_1$$

One way to quell noise, then, is to drive down this quantity, an approach first introduced in~\cite{OsherFatemi92} in the context of noisy images and later applied by~\cite{Chartrand11} to the derivative estimation problem. Penalizing a noisy derivative's total variation, in balance with its fidelity to the data, turns the ill-posed estimation problem into a well-posed one, i.e.~\textit{regularizes} it. We can describe this with the following optimization:
\begin{equation}\label{eqn:tvr-loss}
\underset{\hbx}{\text{minimize}}\ \big\|\by - \hbx\big\|_2^2 + \gamma\, \text{TV}\Big(\frac{d^\nu \hbx}{dt^\nu}\Big)
\end{equation}

\noindent where $\by$ are noisy data, $\hbx$ is an estimate of the smoothed signal, $\gamma$ balances the terms, $\nu$ is the order of the derivative we wish to flatten, and the derivative, which cannot literally be taken on a vector of samples, is found as a vector of Finite Differences. With the insertion of Finite Difference to make the cost precise, the minimization problem becomes convex and is directly solvable with tools like CVXPY~\cite{cvxpy}. We can replace the sum of squares in the first term with a Huber loss (\autoref{eqn:huber}) to obtain some outlier-robustness.

\begin{figure}[!t]\label{fig:tvr-cc}
  \centering
  \includegraphics[width=0.99\textwidth]{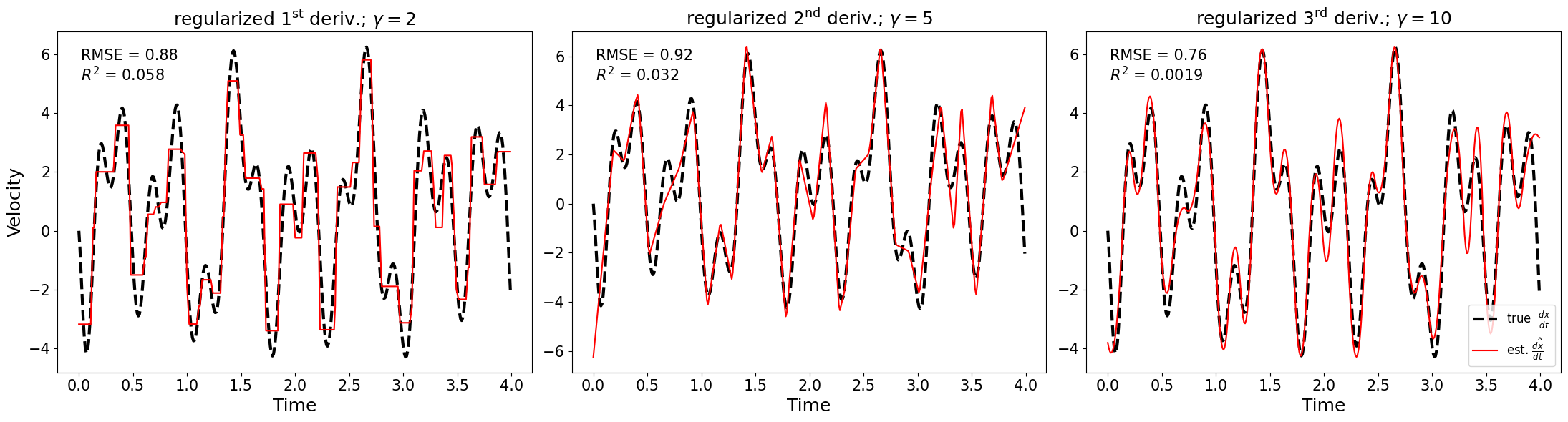}
  \vspace{-3mm}
  \caption{Total Variation Regularized derivatives on noisy cruise control data (\autoref{fig:cruise-control-sim}). Notice the piecewise step, linear, and quadratic character as regularized derivative order increases.}
\end{figure}

Because the second term in \autoref{eqn:tvr-loss} is essentially an $\ell_1$ penalty, the solution is encouraged to have mostly zeros in its $\nu^\text{th}$ finite difference derivative, especially as $\gamma$ grows. This tends to result in a sparse selection of step changes in the $\nu^\text{th}$ derivative, giving TVR derivatives a particular character: piecewise constant for $\nu = 1$, piecewise linear for $\nu = 2$ (the integral of piecewise constant $2^\text{nd}$ derivative), and piecewise quadratic for $\nu = 3$. Examples are shown in \autoref{fig:tvr-cc}. These stark shapes can optionally be convolved with a Gaussian kernel to soften corners.

\subsection{Kalman Smoothing with a Naive Model}
\label{sec:kalman-constant-deriv}

If a model is unknown, then it is possible to assume a model that, like Total Variation Regularization (TVR, \autoref{sec:tvr}), attempts to stabilize the $\nu^\text{th}$ derivative, where $\nu$ is some order, and then apply model-based techniques (\autoref{sec:noisy-with-knowledge}) to estimate the derivative.

For example, if we assume a system model of the form:
\begin{align*}
\bx_n &= \bA\bx_{n-1} + \bB\bu_n + \bw_n,\ \ \mathbb{E}[\bw\bw^\mathrm{T}] = \bQ\\
\by_n &= \bC\bx_n + \bv_n,\ \ \mathbb{E}[\bv\bv^\mathrm{T}] = \bR
\end{align*}
then we can stabilize the second derivative by choosing:
\begin{align*}
\bx_n &= [x_{n,0}, x_{n,1}, x_{n,2}]^\mathrm{T} = [\text{smooth sample},\ 1^\text{st}\text{ deriv. sample},\ 2^\text{nd}\text{ deriv. sample}]^\mathrm{T} \\
\by_n &= y_n,\ \text{the $n^\text{th}$ noisy sample}
\end{align*}
$$\bA = \begin{bmatrix}1 & \Delta t & \frac{\Delta t^2}{2}\\ 0 & 1 & \Delta t\\ 0 & 0 & 1\end{bmatrix},\quad \bB = \begin{bmatrix} 0 \\ 0 \\ 0 \end{bmatrix},\quad \bC = \begin{bmatrix} 1 & 0 & 0\end{bmatrix}$$

\noindent This performs Taylor integrations for the smoothed data and first derivative estimates, while holding the second derivative---sometimes called ``acceleration" regardless of what the data represents---constant. Noisy observations flow into the compatible first dimension of the state and do not directly affect the latent derivatives.

Unlike in TVR, where stability of the $\nu^\text{th}$ derivative is maximized, the amount of stability in Kalman smoothing is squishy, depending on the noise level in the corresponding dimension of the process noise covariance matrix, $\bQ$. If we assume process noise forces solely the last dimension, acting indirectly on its integrals over time, and note the measurement noise covariance matrix, $\bR$, is $1 \times 1$, then the covariance matrices can be treated as two scalar hyperparameters, $q$ and $r$:
$$\bQ = q \begin{bmatrix} \frac{\Delta t^5}{20} & \frac{\Delta t^4}{8} & \frac{\Delta t^3}{6} \\ \frac{\Delta t^4}{8} & \frac{\Delta t^3}{3} & \frac{\Delta t^2}{2} \\ \frac{\Delta t^3}{6} & \frac{\Delta t^2}{2} & \Delta t \end{bmatrix},\quad \bR = \begin{bmatrix} r \\ \end{bmatrix}$$

\noindent where the entries of the discrete-time $\bQ$ matrix (see \autoref{sec:kalman-irregular-dt} for derivation) perform integration similar to $\bA$.\footnote{$\bA$ and $\bQ$ together encode smoothness rather than physics, fixing the covariance between sample values at any two times. In other words, they specify a Gaussian process~\cite{sarkka-sdes}. Most covariance functions do not admit a recursive form like this and must instead be handled by forming and factoring a covariance matrix directly.}

\begin{figure}[!t]\label{fig:rts-freq-response}
  \centering
  \includegraphics[width=0.75\textwidth]{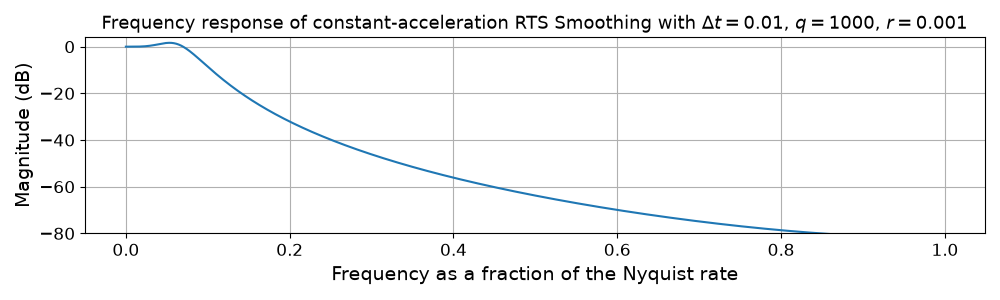}
  \vspace{-3mm}
  \caption{The frequency response of a constant-derivative RTS smoother (derivation in \autoref{ft:rts-filt}) is essentially yet another low-pass filter. Raising $q$ relative to $r$ shifts the peak rightward and broadens it out.}
\end{figure}

We can now solve the model-based Maximum \textit{A Posteriori} (MAP) optimization problem from before (Equations \ref{eqn:MAP-kalman-smooth} and \ref{eqn:MAP-opt-generalized}) to produce a smoothed signal and derivatives. In fact, the $\ell_2$ minimizer is yet another explicit low-pass filter, with response\footnote{\label{ft:rts-filt}To derive, choose $\Delta t$, $q$, and $r$; iterate \autoref{algo:kalman-filter-algo} to convergence to find steady-state gain, $\bK_{ss}$; form the closed-loop discrete LTI system, $G_\text{cl} = ((\mathbb{I} - \bK_{ss}\bC)\bA, \bK_{ss}, \bC, 0)$ using standard $(\bA, \bB, \bC, \mathbf{D})$ notation~\cite{robust-control}; turn this into a transfer function with the formula $G(z) = \bC(z\mathbb{I}-\bA)^{-1}\bB+\mathbf{D}$~\cite{robust-control}, and find its frequency response. A similar process for the RTS backward pass (\autoref{sec:rauch-tung-striebel}) yields $\bL_{ss}$, which, by treating the whole Kalman state as input and RTS state as output, can be used with \autoref{eqn:rts} to form the Multiple-Input Multiple-Output (MIMO) system $G_\text{cl,RTS} = (\bL_{ss}, \mathbb{I} - \bL_{ss}\bA,\mathbb{I},0)$, whence we can extract a frequency response for the first dimension. The response of the two applied in sequence is the product.} plotted in \autoref{fig:rts-freq-response}.

\begin{figure}[!t]
  \centering
  \includegraphics[width=0.99\textwidth]{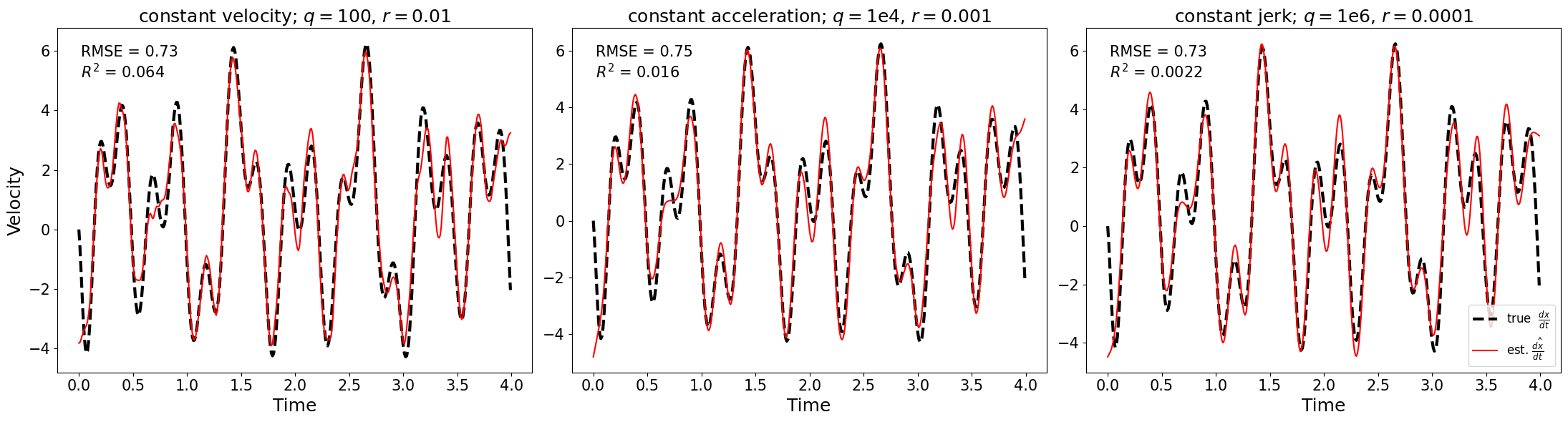}
  \vspace{-3mm}
  \caption{Constant-derivative-model RTS smoothing derivatives of noisy cruise control data (\autoref{fig:cruise-control-sim}) with three different choices of hyperparameters. Derivatives are plucked directly from the second dimension of the state.}
\end{figure}

\begin{figure}[!t]\label{fig:outlier-robust}
  \centering
  \includegraphics[width=0.99\textwidth]{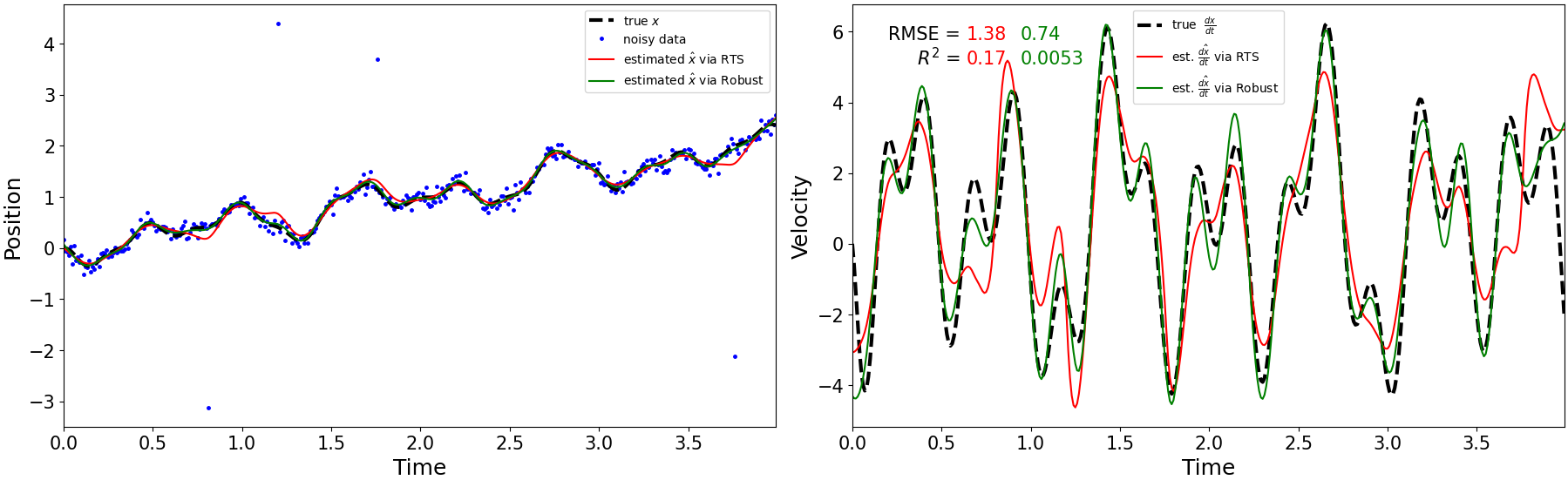}
  \vspace{-3mm}
  \caption{Constant-derivative-model RTS ($\ell_2$ norm loss) and Robust (Huber loss) solutions to the Kalman MAP problem (\autoref{eqn:MAP-opt-generalized}) on noisy cruise control data with outliers (blue points in left panel). Derivatives are optimized for RMSE directly, demonstrating how an ``accurate" fit to outliers can corrupt the smoothed estimate and derivative. Performance deteriorates less with the Huber loss, because it is better able to ignore outliers by less intensely penalizing their neglect. The solution is not as accurate as that in \autoref{fig:robust-demo}, because the model does not match the true underlying system.}
\end{figure}

Just as in the model-based case (\autoref{sec:robust-estimation}), we can also optimize the MAP problem with alternative distance metrics like the $\ell_1$ norm or Huber loss to attain desired qualities like sparse solutions or robustness to outliers, as demonstrated in \autoref{fig:outlier-robust}, although this adds Huber-$M$ hyperparameters.

\subsection{Local Linear Dynamics Fits}
\label{sec:local-linear}

Assume the data is generated by an autonomous linear system, so that the true, underlying $x$ obeys an ODE of order $\nu$ with constant coefficients, $\{a_j\}$, depending only on its own state:
$$x^{(\nu)}(t) = \sum_{j=1}^{\nu} a_j\, x^{(\nu-j)}(t)$$

If derivatives are gathered into a state vector and a few tautological rows stacked in to shift orders' positions between $\bx$ and $\dbx$, this is our old friend $\dbx = \bA\bx$, a class of model encompassing many physical systems, e.g.~the damped oscillator. As written, it relies on the very derivatives we are trying to estimate, however the \textit{relationship} still holds if we take the repeated integral of both sides until no derivative terms remain:
$$x(t) = \Big(\sum_{j=1}^{\nu} a_j \,{\textstyle\int^{(j)}} x(t)\Big) + c(t)$$
\noindent where $\textstyle\int$s are taken w.r.t.~$t$ and $c(t)$ is a polynomial of degree $\nu-1$, a constant of repeated integrations. This can be stacked to form a more linear-algebraic expression:
$$x(t) = \mathbf{a}^\mathrm{T}\begin{bmatrix} \int^{(\nu)} x(t) \\ \vdots \\ \int x(t) \end{bmatrix} + \mathbf{c}^\mathrm{T} \begin{bmatrix} \frac{t^{\nu-1}}{(\nu-1)!} \\ \vdots \\ 1 \end{bmatrix}$$
\noindent where $\mathbf{a} \in \mathbb{R}^\nu$ contains all the $\{a_j\}$, the non-tautological row of $\bA$, and the leftover polynomial has been broken into a Taylor Series basis, with constants collected in a $\nu$-vector.

This we can use, because integrals average noise rather than amplify it. Substituting samples $\by$ as a row vector in place of $x$, we can simply use the trapezoidal rule to accomplish sampled integrations, and form:
$$\hbx = \mathbf{a}^\mathrm{T}\begin{bmatrix} \texttt{trapz}^{(\nu)}(\by) \\ \vdots \\ \texttt{trapz}(\by) \end{bmatrix} + \mathbf{c}^\mathrm{T} \begin{bmatrix} \frac{\mathbf{t}^{\nu-1}}{(\nu-1)!} \\ \vdots \\ \mathbf{1} \end{bmatrix}$$
where $\mathbf{t}$ is a row vector of sample times and $\mathbf{1}$ is a whole row of 1s. Let's call the matrix of data integrals $\int \mathbf{Y}$ and the matrix encoding the Taylor basis $\bB$, both $\in \mathbb{R}^{\nu \times N}$. If we take one fewer integrals in the setup, then we can likewise obtain a derivative estimate:
$$\hdbx = \mathbf{a}^\mathrm{T}\mathbf{Y} + \mathbf{c}^\mathrm{T} \mathbf{\dot B}$$
\noindent where $\mathbf{Y}$ has one fewer \texttt{trapz} applied to all rows, and $\mathbf{\dot B}$ is just a row-shift (upward) of $\bB$ due to the Power Rule. In practice  we replace $\by$ in the last row of $\mathbf Y$ with $\hbx$ to keep noise from sneaking back in.

Thus if we can solve for $\mathbf{a}$ and $\mathbf{c}$ that fit the data, then we have a method to find a model, as opposed to assuming one as in naive Kalman smoothing (\autoref{sec:kalman-constant-deriv}). Set up an optimization:
$$\underset{\mathbf{a},\, \mathbf{c}}{\text{minimize}}\ \big\| \by - \big(\mathbf{a}^\mathrm{T}\textstyle\int\mathbf{Y} + \mathbf{c}^\mathrm{T}\bB\big) \big\|_2^2\, +\, \gamma \|\mathbf{c}\|_1 \,+\,\epsilon \|\mathbf{a}\|_1$$
\noindent Here the first term ensures fidelity to the data. The second term penalizes the constants of integration, forcing the model's dynamics to account for more of the variation, allowing less to leak out and be explained by the polynomial basis. The third term is for solution uniqueness, because for smooth signals integrals of the data can look a lot like the polynomials in $\bB$, making $\mathbf{a}$ somewhat interchangeable with $\mathbf{c}$. We keep its constant small so as not to bias $\mathbf{a}$ and throw off the actual dynamics. We use 1-norm penalties because we want the shrinkage force to be consistent all the way to 0 and across model terms.\footnote{Partly we expect sparse dynamics and constants as a prior, but we also eschew the 2-norm because it can over-penalize legitimate high frequencies. E.g., an oscillator with frequency $\omega$ solves $\ddot x = -\omega^2 x$, so with state $[x, \dot x]^\mathrm{T}$, $\mathbf{a} = [-\omega^2, 0]^\mathrm{T}$ and $\nabla_\mathbf{a} \|\mathbf{a}\|_2^2 = [2 a_1, 0]^\mathrm{T} = [-2\omega^2, 0]^\mathrm{T}$, which shrinks high frequencies more, while $\nabla_\mathbf{a} \|\mathbf{a}\|_1 = [\mathrm{sign}(a_1), 0]^\mathrm{T} = [-1, 0]^\mathrm{T}$ irrespective of frequency.}

This problem is convex and can be solved with the same off-the-shelf tools behind TVR and Robust Kalman. To aid the convex solver, we can make the matrices better-conditioned by scaling and centering the data and working in normalized time units $\tau \in [0,1]$.

For all but the simplest systems, a linear dynamical fit, even this data-driven one, can only approximate the right answer \textit{locally}, in a short window. Dynamics can change, and nonlinearities become more apparent over greater periods. So to apply this method to longer data, we take a sliding-window approach, exactly as with local Polynomial Fits (\autoref{sec:sliding-poly}), combining strided solutions according to a weighting kernel.

\begin{figure}[!t]
  \centering
  \includegraphics[width=0.99\textwidth]{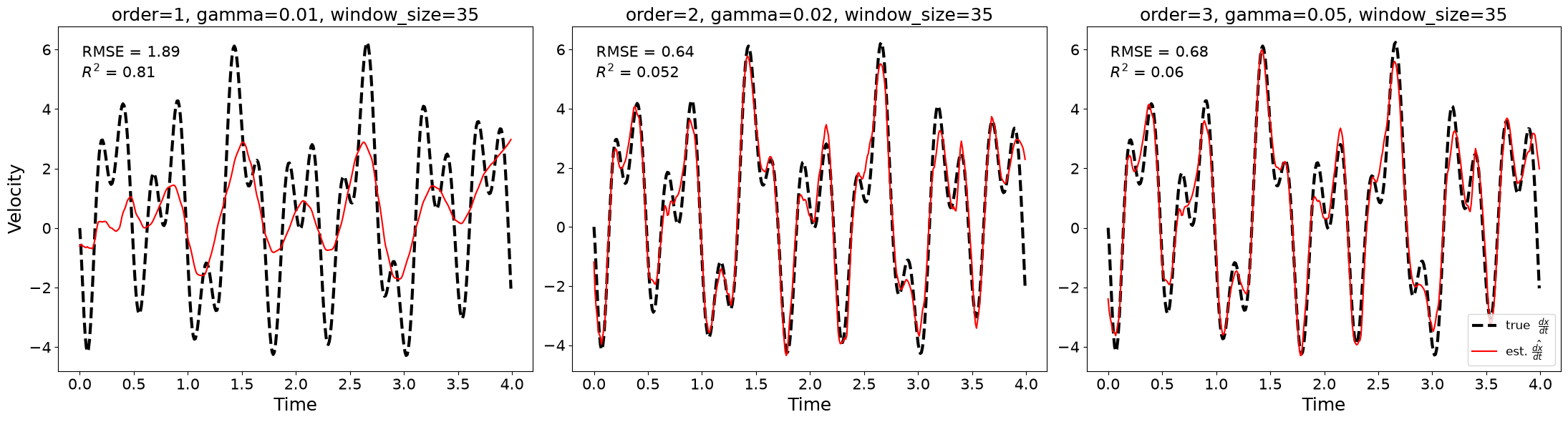}
  \vspace{-3mm}
  \caption{Sliding-window local linear dynamics-based differentiation of noisy cruise control data (\autoref{fig:cruise-control-sim}) with three different choices of hyperparameters.}
\end{figure}

\subsection{Performance Comparison}
\label{sec:performance-comparison}

This section compares smoothing-based numerical differentiation methods, examining relative costs and advantages.\footnote{Experiments in this section are implemented in \path{notebooks/4_performance_analysis.ipynb} of the \texttt{PyNumDiff} repository, at version v0.3. Datasets are simulated, so the whole study reruns from the notebook.}

\subsubsection{Computational Complexity}

The first consideration is computational cost, addressed in \autoref{ta:complexity}. Each method has some complexity in itself,\footnote{Complexity bounds in \autoref{ta:complexity} assume equispaced data. Run times can be worse for irregularly-spaced data. See \autoref{sec:practical}.} which usually takes the form $O(N \cdot f(\Phi))$, where $\Phi$ is the set of hyperparameters, interpreted in each method's explanation above and in \texttt{PyNumDiff} code documentation~\cite{pynumdiff}. The only exceptions are Fourier Spectral, which is log-linear in $N$, and methods using a convex solver. Also of note is the \textit{size} of each method's hyperparameter set, $|\Phi|$, because although the choice of hyperparameters is in principle reduced to the single $\gamma$ from Equation \ref{eqn:floris-cost} or \ref{eqn:pavel-cost}, governing smoothness, we still need to run several Nelder-Mead optimizations to find the optimal $\Phi$, and this is more challenging in higher dimension.

\begin{table}[!t]
\centering
\begin{threeparttable}
\caption{\label{ta:complexity} Computational complexity and hyperparameter set sizes of smoothing methods.} 
\begin{tabular}{p{6cm}p{5.5cm}c}
  \hline \textbf{Method} & \textbf{Complexity} & $\mathbf{|\Phi|}$ \\ \hline
  Filtering $\rightarrow$ Finite Difference & $O(N \cdot \text{\footnotesize window size})$ & 3\\
  Iterated Finite Difference & $O(N \cdot \text{\footnotesize order} \cdot \text{\footnotesize iterations})$ & 2\\
  Sliding Polyfit & $O(N \frac{\text{(window size)}^3}{\text{stride}})$ & 4\\
  Savitzky-Golay & $O(N \cdot \text{\footnotesize window size} \cdot \text{\footnotesize smoothing size})$ & 3 \\
  Splines & $O(N \cdot \text{\footnotesize degree}^2)$ & 2 \\
  Fourier on Extension & $O(N\log N)$ & 3\\
  Radial Basis Fit & $O(N \rho^2)$ & 2\\
  Wavelets & $O(NL\, 2^J)$\tnote{$\ast$} & 3\tnote{$\dagger$} \\ 
  Total Variation Regularization & $O(N^{1.5}\nu^2 \log(\frac{N}{\epsilon}))$\tnote{$\ddagger$} & 3\\
  RTS Smoothing & $O(N \nu^3)$ & 2\tnote{\S}\\
  Robust MAP Smoothing & $O(N^{1.5}\nu^3 \log(\frac{N}{\epsilon}))$\tnote{\P} & 5\\ 
  Sliding Local Linear Dynamics & $O(N \frac{\text{window size}\,\cdot\,\nu^2}{\text{stride}} \sqrt\nu\log (\frac{\nu}{\epsilon}))$\tnote{$\|$} & 5 \\[2pt] \hline
\end{tabular}
\begin{tablenotes}
  \item[$\ast$] \footnotesize The DWT recurses to operate on a single half-length problem at each new level, making it linear in $N$ by convergence of the geometric series $1 + \frac{1}{2} + \frac{1}{4} + ... = 2$. $L$ is filter length, which depends on the wavelet type, and all is repeated across cycle-spinning start locations.
  \item[$\dagger$] Level $J$ is excluded, since it can be automatically selected.
  \item[$\ddagger$] Interior-point methods handle absolute values from $\ell_1$, Huber, and TV terms by processing into inequalities then softening into --log ``barrier" penalties, which blow up to $\infty$ with argument 0~\cite{tibshirani2019convex}. With $m = O(N)$ total barrier terms, these algorithms reach duality gap (which upper-bounds suboptimality) $\epsilon$ in $O(\sqrt m \log(\frac{m}{\epsilon}))$ Newton steps~\cite{Boyd_Vandenberghe_2004, socp-complexity-mosek}. Here each step costs a linear solve, $O(N \nu^2)$ because the system is banded. In practice, solvers like CLARABEL~\cite{goulart-clarabel} take a nearly-constant 10--50 steps, so observed run time is closer to linear in $N$.
  \item[\S] Scaling the two noise hyperparameters, $q$ and $r$ (\autoref{sec:kalman-constant-deriv}), by the same factor does not change the RTS frequency response (\autoref{fig:rts-freq-response}), because the center of the joint distribution found at the Kalman filter's final step (\autoref{fig:kalman-steps}) is invariant to absolute scale~\cite{pynumdiff}. Thus only the ratio between $q$ and $r$ makes a difference. Note this simplification is not possible for generalized MAP smoothing, because the absolute scales of $q$ and $r$ interact with distinct process and measurement Huber-$M$s.
  \item[\P] Piecewise linear quadratic (PLQ) losses~\cite{aravkin3} like $\ell_1$ and Huber can be optimized by CLARABEL, where the number of barrier penalties is $O(N)$ like for TVR. When applied to measurements in series, the full system is sparse as in RTS (\autoref{sec:rauch-tung-striebel}), allowing linear per-iteration cost in $N$, rather than cubic as for general, dense Second Order Cone Programs (SOCPs)~\cite{socp-matmul, socp-complexity-mosek, tibshirani2019convex}.
  \item[$\|$] Each fit is a Quadratic Program in $O(\nu)$ unknowns over a $|\text{window size}|$ number of residuals, so forming the reduced KKT system at each step costs $O(\text{window size} \cdot \nu^2)$~\cite{Boyd_Vandenberghe_2004}. Only the coefficients are subject to $\ell_1$, so the number of barrier penalties is $O(\nu)$, and an interior point method reaches tolerance $\epsilon$ in $O(\sqrt\nu \log(\frac{\nu}{\epsilon}))$ iterations.
\end{tablenotes}
\end{threeparttable}
\end{table}

\subsubsection{Accuracy and Bias}
\label{sec:accuracy-and-bias}

The other major consideration is method accuracy in terms of RMSE (\autoref{eqn:rmse}) and bias in terms of Error Correlation (\autoref{eqn:ec-r2}) against the true, withheld derivative. We are interested in whether the choice of method should change with different data characteristics: outliers, noise type, noise level, and the sampling rate. We are also interested in whether the idiosyncrasies of particular data series favor particular differentiation algorithms and would like to know how sensitive methods are to the choice of bandlimit, which affects the smoothness of the optimized answers through the loss function's $\gamma$, set according to the heuristic in \autoref{eqn:gamma-heuristic}. In the presence of outliers, we optimize the robust loss function in \autoref{eqn:pavel-cost} with Huber parameter $M = 2$, and in their absence choose $M = 6$.\footnote{We have observed empirically that, under conditions of Gaussian noise, optimizing the robust loss with $M\!=\!6$ essentially always produces identical answers to optimizing the original RMSE-based loss, \autoref{eqn:floris-cost}. This makes sense theoretically, because inlier errors beyond $6\sigma$ from the mean error are exceedingly rare: \texttt{2*(1-scipy.stats.norm.cdf(6)) = 1.97e-9}. Lowering $M$ to 2 provides a slight advantage in the presence of outliers, weighting them lighter, but in the absence of outliers $M\!=\!2$ begins to degrade performance, yet only very slightly, because the portion of inliers within $2\sigma$ is more than 95\%.}

\begin{figure}[!t]\label{fig:simulations}
  \centering
  \includegraphics[width=0.99\textwidth]{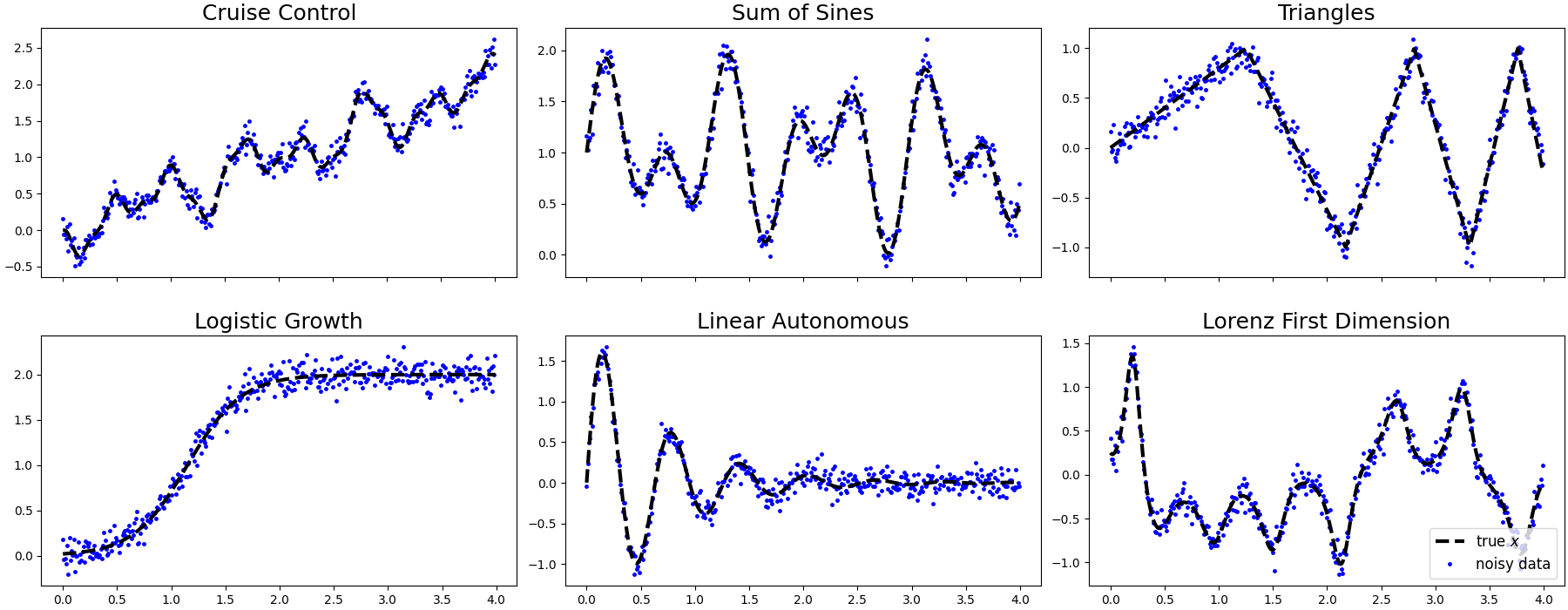}
  \vspace{-2mm}
  \caption{Six simulations are used to produce noisy data and true underlying signal, as well as true derivatives (not shown). In the pictured examples, $t \in [0, 4)$, $\Delta t = 0.01$, and noise is Gaussian, $\eta \sim \mathcal{N}(\mu=0, \sigma=0.1)$. Both $\Delta t$ and $\eta$ are varied to test differentiation methods.}
\end{figure}

To examine these questions, we use six simulated examples, shown in \autoref{fig:simulations}, each exhibiting its own unique attributes. Two are chosen by shape, to be composed of periodic or piecewise linear components, and the others arise from systems: two linear, as might be encountered in control theory, one nonlinear and chaotic, and one from mathematical biology. All simulations are on a similar scale so signal-to-noise ratio is comparable (though not identical), and all can be readily recomputed with different step sizes, $\Delta t$.

\begin{figure}[!t]\label{fig:noise-distributions}
  \centering
  \includegraphics[width=0.6\textwidth]{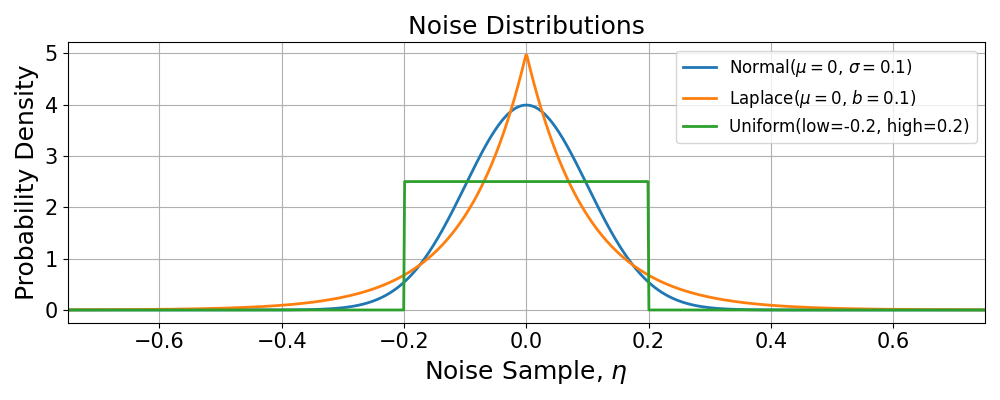}
  \vspace{-3mm}
  \caption{Noise can be sampled from a variety of probability density functions. For the purposes of our experiment, we use these three. Increasing noise scale is equivalent to spreading density horizontally while keeping integral 1. Energy content (variance $= \int \eta^2 \rho(\eta) d\eta$) is different for each: $\sigma^2 = 0.01$ for Gaussian, $2b^2 = 0.02$ for Laplace, and $\frac{1}{3}\text{low}^2 = 0.01\bar{3}$ for symmetric Uniform.}
\end{figure}

We add noise from one of three different distributions, visualized in \autoref{fig:noise-distributions}: Normal Gaussian with $\mu,\sigma\!=\!(0, 0.1)$, Laplace with $\mu,b\!=\!(0, 0.1)$, or Uniform $\in [-0.2, 0.2]$. The widths of these distributions are varied by multiplying their parameters by a noise scale. We further try corrupting a randomly-chosen 1\% of samples with outlier values sampled uniformly from $\pm[50\%, 150\%]$ of the max -- min of their data series (e.g.~left panel of \autoref{fig:outlier-robust}).

\begin{figure}[!t]\label{fig:power-spectra}
  \centering
  \includegraphics[width=0.99\textwidth]{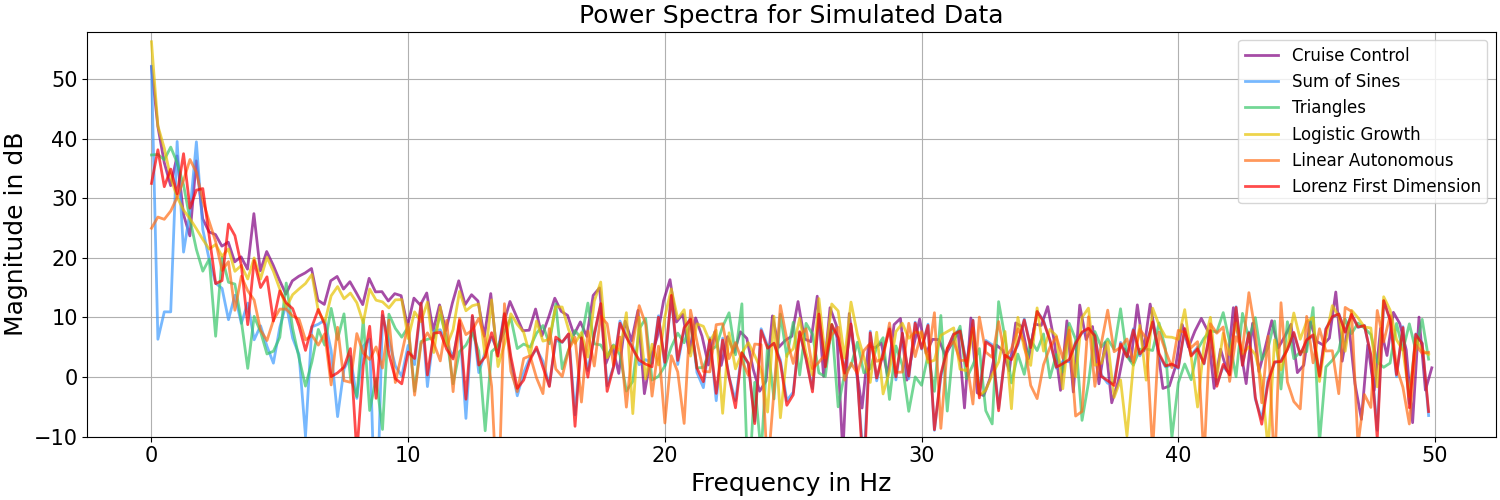}
  \vspace{-3mm}
  \caption{The power spectra of data $\by = \bx + \boldsymbol\eta,\ \boldsymbol\eta \sim \mathcal{N}(\mu\!=\!0,\, \sigma\!=\!0.1)$, calculated as the positive frequencies of $10\log_{10}|\text{FFT}(\by)|^2$. Spectrum extends up to the Nyquist frequency $ = \frac{f_s}{2} = \frac{1}{2\Delta t}$, which for $\Delta t = 0.01$ is 50 Hz.}
\end{figure}

To determine a reasonable bandlimit, we examine signals' power spectra in \autoref{fig:power-spectra} and observe a drop-off as in \autoref{fig:noise-spectrum}. This occurs more rapidly for some than for others, but the same general shape holds, even with other choices of noise type, noise scale, and $\Delta t$. We choose to center our examination around 3 Hz.

\begin{figure}[!t]\label{fig:search-space}
\centering
\begin{tikzpicture}[scale=0.6]
  \draw[thick](-3,0) -- (6,0) node[right] {$|\Delta t\text{s}| = 4$};
  \foreach \x in {-3, 0, 3, 6}{
  \fill (\x,0) circle(2pt);}
  
  \draw[thick](0,0) -- (35:3) node[above] {$|\text{whether outliers}| = 2$};
  \fill (35:3) circle (2pt);
  
  \draw[thick](165:3) -- (-15:6) node[right] {$|\text{noise scales}| = 4$};
  \foreach \r in {-3, 3, 6}{
  \fill (-15:\r) circle(2pt);}

  \draw[thick](-165:6) -- (15:6) node[right] {$|\text{bandlimits}| = 5$};
  \foreach \r in {-6, -3, 3, 6}{
  \fill (15:\r) circle(2pt);}

  \draw[thick](145:3) -- (-35:3) node[below] {$|\text{noise types}| = 3$};
  \foreach \r in {-3, 3}{
  \fill (-35:\r) circle(2pt);}
\end{tikzpicture}
\vspace{-2mm}
\caption{Notional look at slices in the performance analysis configuration space. Figures \ref{fig:vary-outliers}--\ref{fig:vary-bandlimit} each correspond to one of the lines sketched here. The common central point is \{bandlimit=3 Hz, $\Delta t$=0.01, noise type=normal, noise scale=1, outliers=False\}.}
\end{figure}
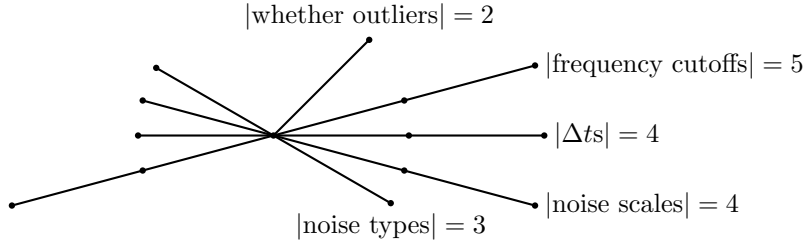

The analysis space now has 7 dimensions: differentiation method, simulation, step size, noise type, noise scale, whether there are outliers, and bandlimit. There is also a bonus dimension: To get a sense of performance center and scatter, we need several rounds of simulation followed by optimization with different random seeds. Some dimensions must be explored exhaustively: 12 methods, 6 simulations, and we turn the number of random seeds all the way up to 52. To save compute and target our analysis at specific questions, we sweep along slices of the remaining orthotope, only varying one parameter at a time. We choose slices to all pass through a common central point, \{bandlimit=3 Hz, $\Delta t$=0.01, noise type=normal, noise scale=1, outliers=False\}, as symbolized in \autoref{fig:search-space}.

\begin{figure}[!t]\label{fig:vary-outliers}
  \centering
  \includegraphics[width=0.99\textwidth]{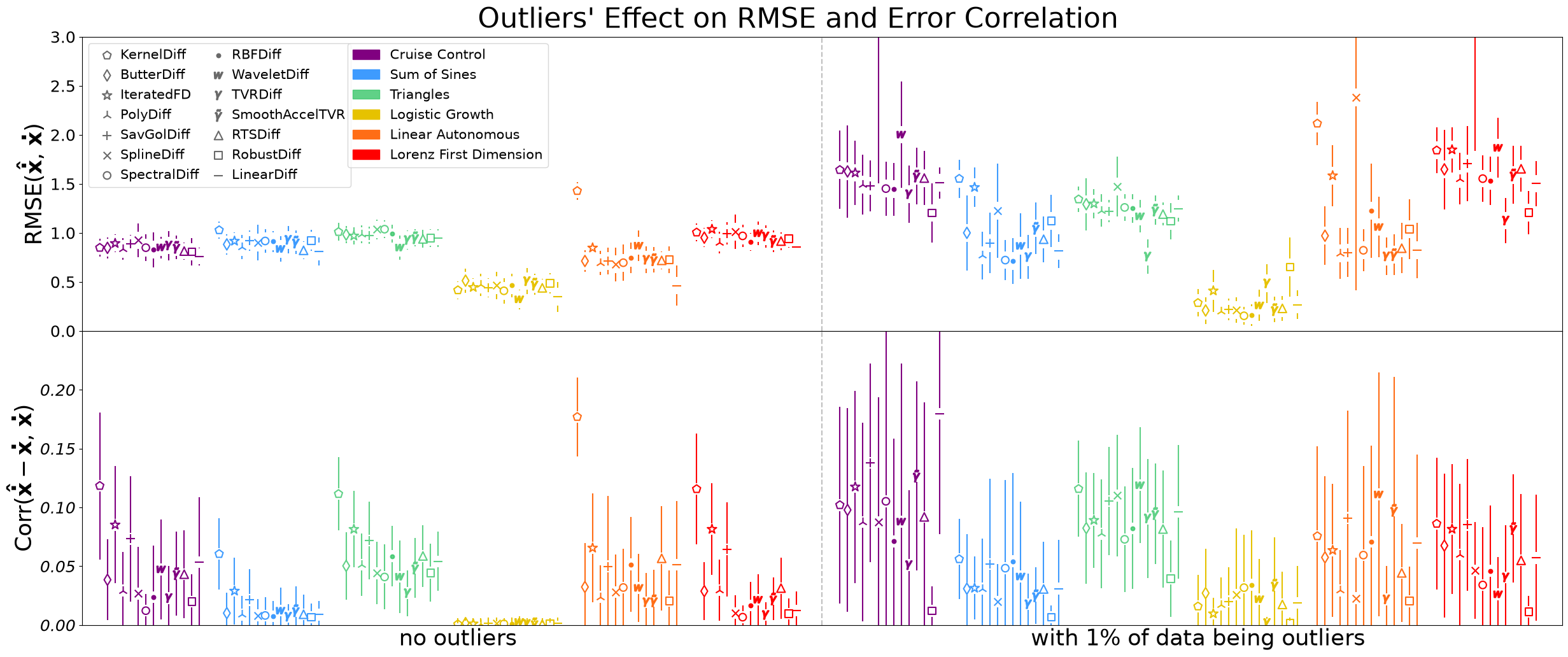}
  \vspace{-3mm}
  \caption{Outliers generally cause worse RMSE and Error Correlation.}
\end{figure}

\begin{figure}[!t]\label{fig:vary-noise-type}
  \centering
  \includegraphics[width=0.99\textwidth]{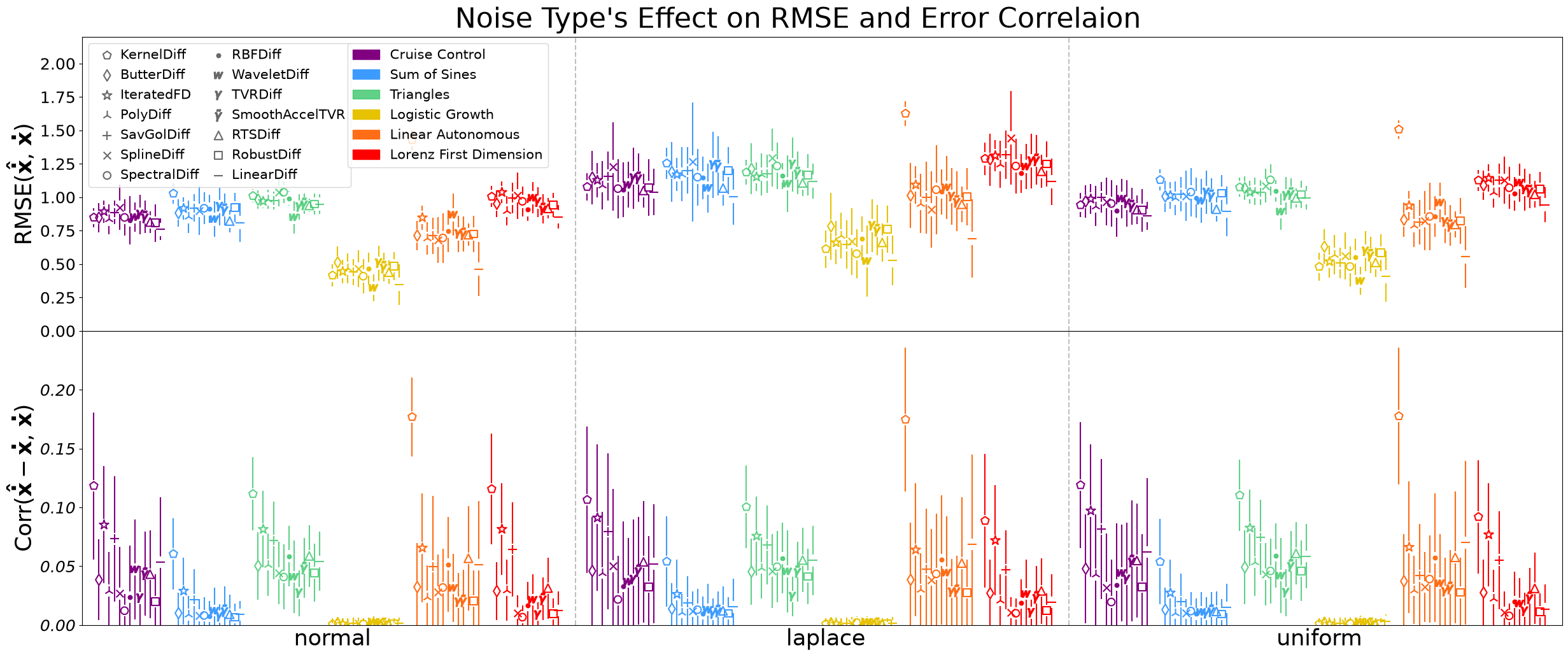}
  \vspace{-3mm}
  \caption{Noise type does not have a profound effect on RMSE or Error Correlation, although the Laplace distribution, with the heaviest tails and greatest noise energy, causes slightly worse performance.}
\end{figure}

\begin{figure}[!t]\label{fig:vary-noise-scale}
  \centering
  \includegraphics[width=0.99\textwidth]{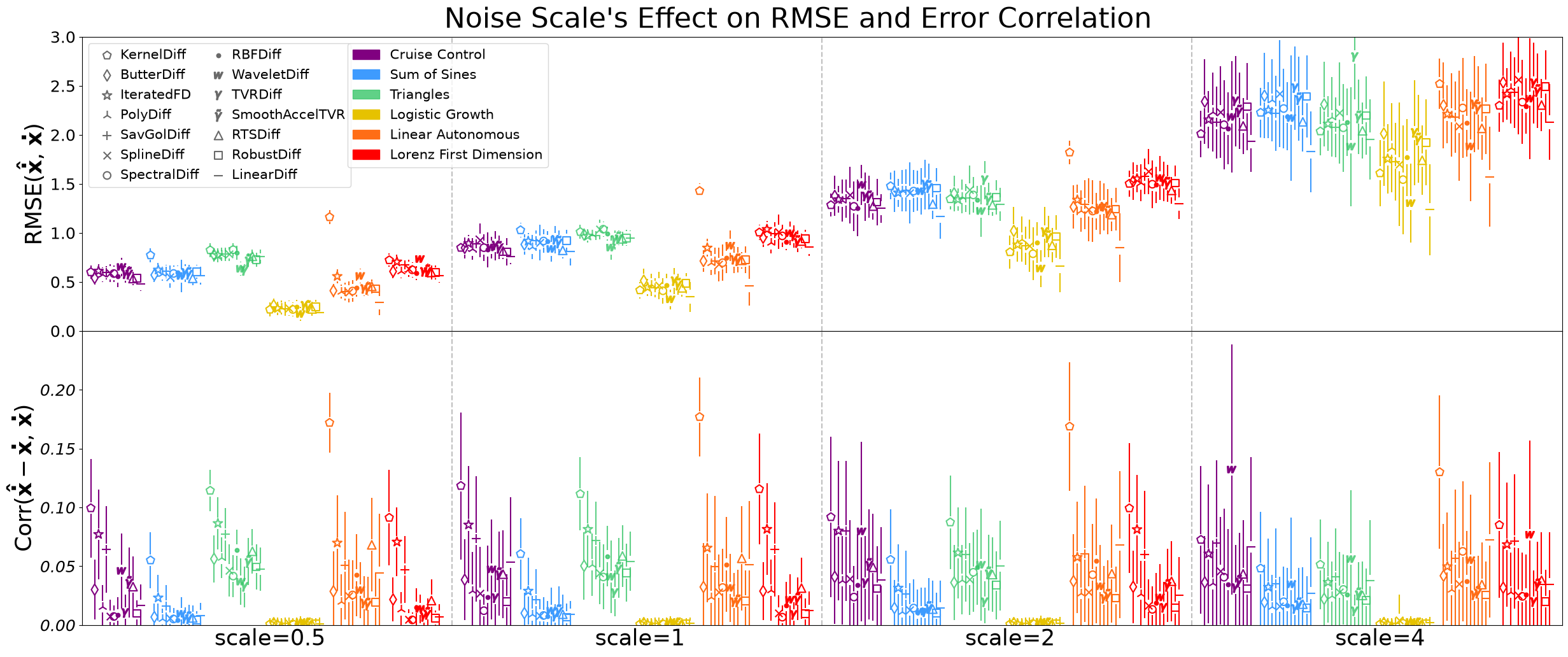}
  \vspace{-3mm}
  \caption{Increasing noise scale degrades average RMSE and Error Correlation simultaneously and increases variability for all methods. Scale=1 corresponds to $\mathcal{N}(0,0.1)$.}
\end{figure}

\begin{figure}[!t]\label{fig:vary-dt}
  \centering
  \includegraphics[width=0.99\textwidth]{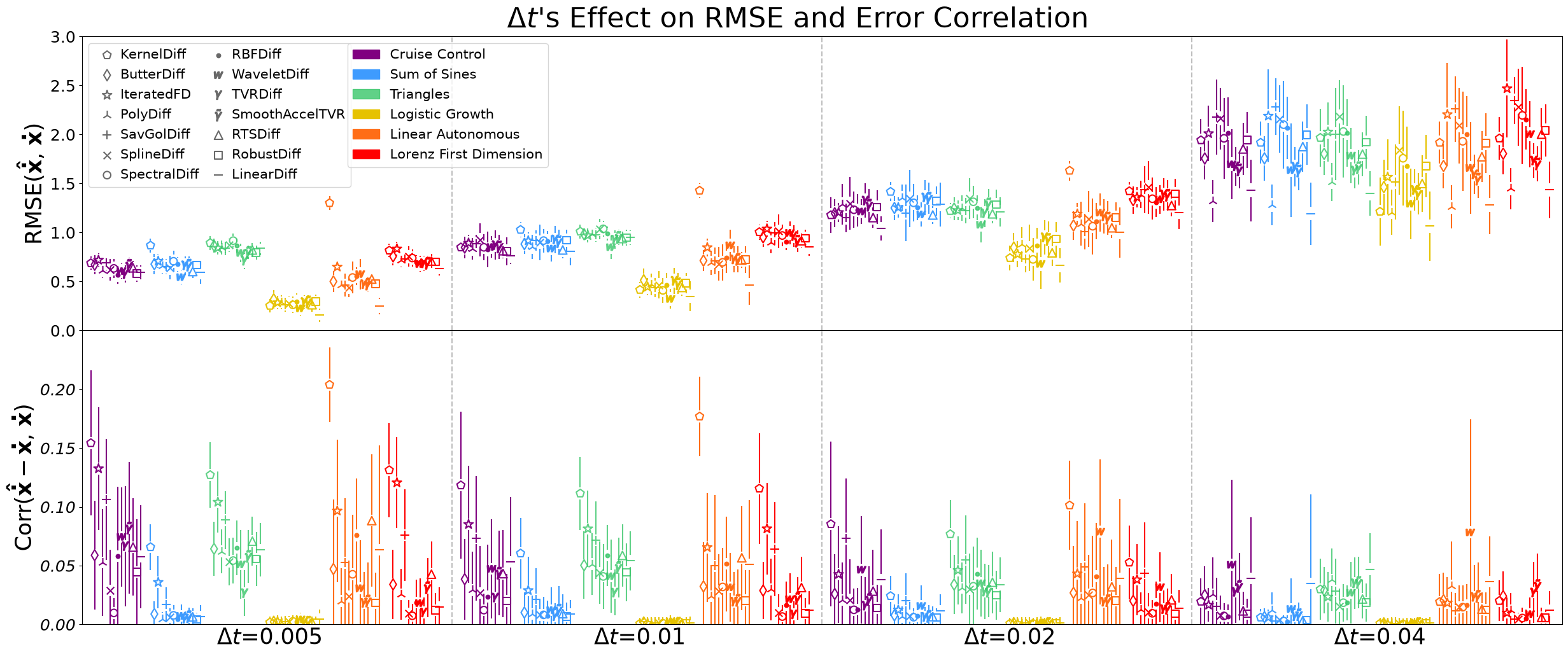}
  \vspace{-3mm}
  \caption{Larger $\Delta t$ are difficult to handle for all methods.}
\end{figure}

\begin{figure}[!t]\label{fig:vary-bandlimit}
  \centering
  \includegraphics[width=0.99\textwidth]{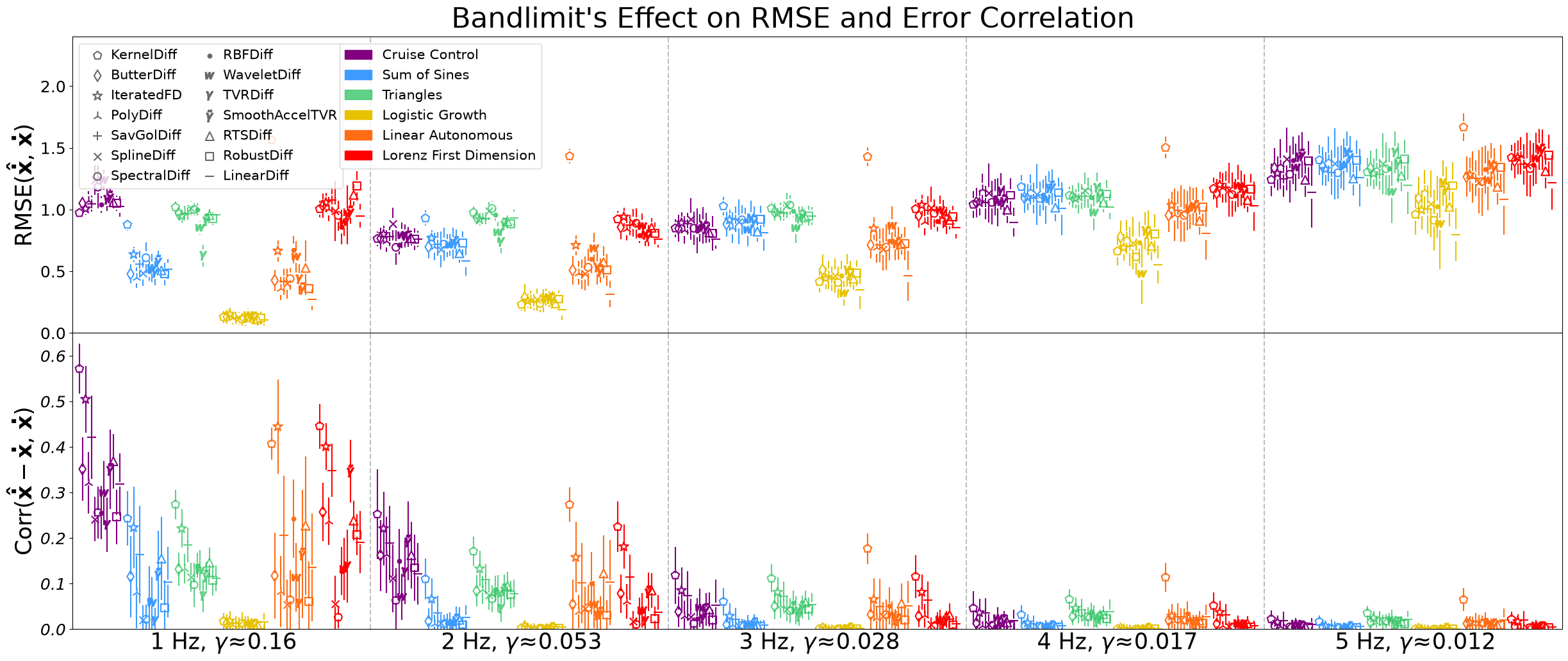}
  \vspace{-2mm}
  \caption{Simulations without sharp corners continue to achieve lower RMSE as the bandlimit is decreased, favoring more smoothness during optimization. But this comes at the cost of sharply increasing Error Correlation, unless the underlying signal truly has very little frequency content.}
\end{figure}

Figures \ref{fig:vary-outliers}--\ref{fig:vary-bandlimit} show RMSE and Error Correlation along each of these slices, with simulations marked by color, differentiation methods differentiated by symbol, and experimental configurations separated into boxes. A method's mean performance across the 52 randomly-seeded runs is marked by its symbol's location, with bars indicating sample standard deviation.\footnote{Confidence intervals on the mean can be calculated as a fraction of sample standard deviation, in this case about 0.278 for 95\%, per \texttt{scipy.stats.t.ppf(1-0.05/2, N-1)/np.sqrt(N)} with $N=52$.} For implementations of all methods, simulations, and the optimization loop used in this section, see \texttt{PyNumDiff v0.3}~\cite{pynumdiff}. Methods are plotted left to right in order of exposition in the above sections. For naming clarity: \texttt{KernelDiff} refers to applying a smoothing kernel followed by FD, and \texttt{SmoothAccelTVR} refers to second-order TVR followed by convolution with a Gaussian kernel to dull sharp corners.

\subsection{Recommendations} The most unmistakable pattern in Figures \ref{fig:vary-outliers}--\ref{fig:vary-bandlimit} is that all performances improve or degrade together with the difficulty of the problem, where difficulty is governed by:

\vspace{2mm}\begin{enumerate}
  \item the presence or absence of outliers
  \item the noise scale
  \item the size of $\Delta t$, because faster sampling raises the Nyquist frequency and allows better band separation of signal from noise
  \item the smoothness of the true signal, because a spectrum that falls off more sharply, e.g.~for Logistic Growth, allows more aggressive mitigation of high frequency noise without cutting into signal.
\end{enumerate}\vspace{2mm}

In our experiments, optimizing robust loss, \autoref{eqn:pavel-cost} rather than \autoref{eqn:floris-cost}, greatly improves all methods' performance in the presence of outliers. Only two methods, \texttt{RobustDiff} and \texttt{TVRDiff} with Huber loss, are specifically designed to handle outliers, by pairing an assumed model (\autoref{sec:kalman-constant-deriv}) with techniques from robust estimation (\autoref{sec:robust-estimation}). They achieve consistently lower bias than other methods and are most accurate for some but not all simulations (\autoref{fig:vary-outliers}). With a high number of hyperparameters and heavy convex solves, \texttt{RobustDiff} is unfortunately the most expensive to optimize. Demonstrating robustness of a different kind, the sliding-window methods, \texttt{PolyDiff} (\autoref{sec:sliding-poly}) and \texttt{LinearDiff} (\autoref{sec:local-linear}), degrade less than other methods as the distance between samples is increased (\autoref{fig:vary-dt}). \texttt{LinearDiff} is also most robust to a too-high choice of bandlimit (\autoref{fig:vary-bandlimit}).

At lower bandlimit, noise scale, and step size, and even in the presence of outliers, \texttt{TVRDiff} (\autoref{sec:tvr}) does noticeably best on the \texttt{Triangles} simulation, because its derivatives can match the piecewise character of the true derivative. \texttt{LinearDiff} also typically does best on the sum-of-sines and linear autonomous simulations, which are based on linear ODEs. These suggest that pairing the model to the problem can have a tangible benefit, but for most simulations the best performer is inconsistent across sweeps and hardly distinguished from the pack. This is good news, because it means data series do not tend to favor particular methods for unknown reasons, that sophisticated methods tend to do equally well for all but special cases with intuitive rationales.

Remaining distinctions are less obvious. No method distinguishes itself as better able to handle a specific type of noise (\autoref{fig:vary-noise-type}), and all methods are similarly challenged by increasing noise scale (\autoref{fig:vary-noise-scale}). The choice of 3 Hz bandlimit is vindicated, because increasing $f$ (seeking less smooth derivatives) causes accuracy to worsen, while decreasing $f$ (seeking smoother derivatives) causes substantially worsening bias for most simulations (\autoref{fig:vary-bandlimit}), although some smoother data series would be better served by lower bandlimits.

All methods are equally handicapped within a configuration, so results on each simulation within each frame are comparable. We can therefore analyze statistics about methods' accuracy and bias \textit{rankings}. \texttt{KernelDiff} and \texttt{IteratedFD} perform worst, averaging to $9.6^\text{th}$ and $9.2^\text{th}$ places on RMSE, respectively. They also rank worst on error correlation, at $11.7$ and $10.3$. \texttt{LinearDiff} and \texttt{RTSDiff} rank best on RMSE, in $3.9^\text{th}$ and $5.8^\text{th}$ places on average, with middling error correlation rankings both at $7.7$. Least biased are \texttt{RobustDiff} and \texttt{TVRDiff}, both at average rank just over $5$. \texttt{LinearDiff}, \texttt{WaveletDiff}, \texttt{SpectralDiff}, and \texttt{RTSDiff} achieve an RMSE among the lowest (best) three in 62\%, 39\%, 26\%, and 22\% of configurations, respectively, and on the other side \texttt{RTSDiff} has an RMSE in the highest (worst) three in only 2.7\%, with \texttt{LinearDiff} in second place at 5.2\% and \texttt{PolyDiff} third best at 13\%.

To further compare methods, auxiliary considerations are summarized in \autoref{fig:flexibility}. These factors emphasize the value of algorithm flexibility, which can be underappreciated in theoretical discussion but is often important in practice. Indeed, one of the reasons neural networks have become so dominant is their endless capacity for recombination into architectures tailored to multivariate, multidimensional, sequential, incomplete, and structured data.\footnote{``Time cannot wither nor custom stale her infinite variety." --Shakespeare on the perceptron} For its outstanding versatility and high base performance, we recommend \texttt{RTSDiff} (\autoref{sec:kalman-constant-deriv}) as a default general-purpose method, with \texttt{LinearDiff} an excellent runner-up.

\begin{figure}[!t]\label{fig:flexibility}
\centering
\begin{tikzpicture}[font=\scriptsize]
\def\nrows{13}
\def\ncols{8}
\def\rowheight{0.5cm}
\def\colwidth{1cm}
\def\colheaders{
    high\\[-1.2pt]accuracy/0.25,
    low\\[-1.1pt]bias/0.35,
    outlier\\[-1.2pt]robust/0.25,
    efficient\\[-1.2pt]to run/0.25,
    fast to\\[-1.2pt]optimize/0.25,
    no heavy\\[-1pt]dependency/0.25,
    variable\\[-1.2pt]step size/0.35,
    graceful with\\[-1pt]missing data/0.35}
\def\rowlabels{Kernel Filtering, Butterworth Filtering, Iterated FD, Sliding Polyfit, Savitzky-Golay, Spline Fit, Fourier on Extension, Radial Basis Fit, Wavelets, TVR, RTS Smoothing, Robust MAP, Local Linear Dyn.}
\foreach \i in {0,...,\ncols} { 
    \draw (\i*\colwidth,0) -- (\i*\colwidth,\nrows*\rowheight);
}
\foreach \j in {0,...,\nrows} { 
    \draw (0,\j*\rowheight) -- (\ncols*\colwidth,\j*\rowheight);
}
\foreach \head/\offset [count=\i from 0] in \colheaders {
    \node[anchor=south west, align=left, rotate=30] at ({(\i+\offset)*\colwidth},{(\nrows-0.2)*\rowheight}) {\head};
}
\foreach \head [count=\j from 0] in \rowlabels {
    \node[anchor=east] at (0,{\nrows*\rowheight - (\j+0.5)*\rowheight}) {\head};
}
\def\griddata{{2,0,0,1,1,1,0,0},
              {1,1,0,1,1,1,0,0},
              {2,2,0,1,1,1,2,0},
              {1,1,0,1,1,1,1,1},
              {1,1,0,1,1,1,0,0},
              {1,1,2,1,1,1,1,1},
              {1,1,0,1,1,1,2,0},
              {1,1,0,1,1,1,1,0},
              {1,1,0,1,1,1,0,0},
              {1,1,1,1,1,0,0,0},
              {1,1,0,1,1,1,1,1},
              {1,1,1,1,2,0,1,1},
              {1,1,0,1,2,0,1,1}};
\foreach \row [count=\y from 0] in \griddata {
  \foreach \cell [count=\x from 0] in \row {
    \ifnum\cell=1
      \node at ({(\x+0.5)*\colwidth},{\nrows*\rowheight-(\y+0.5)*\rowheight}) {\textcolor{green!60!black}{$\checkmark$}};
    \else
      \ifnum\cell=0
        \node at ({(\x+0.5)*\colwidth},{\nrows*\rowheight-(\y+0.5)*\rowheight}) {\textcolor{red}{\textsf{X}}};
      \else
        \node at ({(\x+0.5)*\colwidth},{\nrows*\rowheight-(\y+0.5)*\rowheight}) {$\boldsymbol{\sim}$};
      \fi
    \fi
  }
}
\end{tikzpicture}
\vspace{-1mm}
\caption{Ability of methods to handle auxiliary considerations. Black tildes mean neither yes nor no, either mediocre performance or that a method could be extended to handle a situation, but at the cost of computational efficiency or mathematical or code simplicity. Optimization speed is heavily dependent on the size of hyperparameter search space. ``Heavy dependency" refers to convex solvers, which, although easy to install and use in modern form~\cite{cvxpy}, are based on underlying compilers and iterative algorithms that can fail to solve efficiently or at all, a source of potential brittleness.}
\end{figure}

\section{Irregularly-Spaced Samples}
\label{sec:practical}

We have made reference to a consistent step size, $\Delta x$ or $\Delta t$, throughout this review, but in the context of real data, sample spacing is anything but guaranteed and can even be highly variable. This is a wrinkle worth addressing, lest differentiation be practically limited to a fraction of experimental data streams.

Because each differentiation method is built from distinct mathematics, some handle irregular spacing more naturally than others. Especially in cases where performance may not be all that decisive (\autoref{fig:flexibility}), this can have massive bearing on method choice. \autoref{ta:irregular-dt-workability} summarizes which methods discussed in this review can be made to handle irregular steps and why, and subsequent subsections take a closer look at how several methods natively handle or can be extended to handle this situation.

\begin{table}[!t]
\caption{\label{ta:irregular-dt-workability} How amenable differentiation methods are to variable step size, $\widetilde{\Delta t}$.}\vspace{-2mm}
\centering
\begin{tabular}{lc>{\raggedright\footnotesize\arraybackslash}m{9cm}}
  \hline \textbf{Method}\rule{0pt}{2.65ex} & $\mathbf{\widetilde{\Delta t}}$\textbf{?} & \normalsize\textbf{Notes}\\ \hline
  AutoDiff & N/A & Not for differentiating data samples, rather for sampling fixed derivatives. But can be sampled wherever.\\
  Finite Difference & $\boldsymbol\sim$ & Can work, but more expensive and takes more care, especially if applied iteratively.\\
  Fourier Spectral & $\boldsymbol\sim$ & There is the Non-Uniform FFT~\cite{nufft}, although it is not included in common packages.\\
  Chebyshev & \textsf{X} & No speed advantage if not cosine-spaced samples, and bad basis for separating noise.\\
  Finite Elements & $\checkmark$ & Works with irregular $\Delta t$ on meshes internally, but beware relies on analytic equations as input rather than samples.\\
  Kalman Methods & $\checkmark$ & Evaluating $e^M,\ M\in \mathbb{R}^{m \times m}$ costs $O(m^3)$~\cite{matrix-exp}, so same asymptotic complexity, though more involved math.\\
  Smoothing $\rightarrow$ FD & \textsf{X} & Filters do not behave consistently with irregular $\Delta t$.\\
  Sliding Polyfit & $\checkmark$ & Use sliding window in independent-variable space to capture groups of points across the domain.\\
  Savitzky-Golay & \textsf{X} & Ceases to be a consistent filter, so essentially back to sliding polyfit.\\
  Splines & $\checkmark$ & Can be solved as efficiently as uniform case, and placing knots is no more complicated than before.\\
  Radial Basis Fit & $\checkmark$ & Local basis functions can be plopped over data anywhere, although they interact less if further separated.\\
  Wavelets & \textsf{X} & The discrete wavelet transform requires uniform spacing.\\
  TVR & \textsf{X} & Complicated FD equations turn the cost function gnarly.\\
  Local Linear Dyn. & $\checkmark$ & Integrals and polynomials feeding the fits are indifferent to gap sizes.\\ \hline
\end{tabular}
\end{table}

\subsection{Finite Difference with Irregular Steps}
\label{sec:fd-irregular-dt}

There is nothing stopping us from solving the stencil Vandermonde inverse problem of \autoref{eqn:fd-generalized} with irregular stencils. We do this by multiplying both sides by $\Delta x$ to put stencil locations back in units of the independent variable and then using whichever fractional or decimal stencil locations, $[s_0, ..., s_{S-1}]$, match data locations.

For example, using a ``centered" stencil on three points with distances $d_{01}$ between the first two and $d_{12}$ between the last two gives:
\begin{align*}
y'_n \approx&\ c_0 y_{n-1} + c_1 y_n + c_2 y_{n+1}\\
=&\ c_0 [y_n - d_{01}\cdot y'_n + \frac{(-d_{01})^2}{2!} y''_n + O(d_{01}^3)] + c_1 y_n \\
&+ c_2 [y_n + d_{12}\cdot y'_n + \frac{d_{12}^2}{2!} y''_n + O(d_{12}^3)]\\
\rightarrow & \begin{bmatrix}
1 & 1 & 1\\
-d_{01} & 0 & d_{12} \\
\frac{d_{01}^2}{2} & 0 & \frac{d_{12}^2}{2}
\end{bmatrix}
\begin{bmatrix} c_0 \\ c_1 \\ c_2 \end{bmatrix}
= \begin{bmatrix} 0 \\ 1 \\ 0 \end{bmatrix}
\end{align*}

Solving bespoke linear inverse problems like this incurs an extra $O(\text{order}^3)$ cost per data point. We should also be mindful that error bounds degrade, because $\|\bc\|_2 \leq 1/\sigma_\text{min}$, and the smallest singular value of the matrix, $\sigma_\text{min}$, is maximized with equal distances. In addition, if iterating this process (\autoref{sec:iteratedfinitedifference}), we have to be careful to use correct trapezoid widths when integrating, and handle stencils adaptively at the endpoints to avoid coefficients with magnitude over 1, which amplify noise.

\subsection{Splines with Irregular Steps}

Continuous B-splines (\autoref{fig:bsplines}) and their analytic derivatives can be sampled anywhere in the domain, and knots can be added anywhere, so splines naturally handle irregularly-spaced data. Wherever samples are taken, there are always only $O(d)$ nonzero basis functions at each, where $d$ is the degree of the spline, so the matrix $\bB$ from \autoref{eqn:spline-opt-vec-form} is always sparse, $\bB^\mathrm{T}\bB$ remains banded, and the fit problem is just as easy to solve.

\subsection{Kalman Smoothing with Irregular Steps}
\label{sec:kalman-irregular-dt}

Our exposition of Kalman filtering in \autoref{sec:kalman-filter} uses a discrete-time model of system dynamics, which is convenient when thinking in terms of computation, but we can imagine this model is actually arising from the continuous evolution of a system over constant increments of time~\cite{axelssonandgustafsson, crassidis_junkins}. We can write the underlying continuous model as:

$$\frac{d\bx}{dt} = \mathcal{A}\bx(t) + \mathcal{B}\bu(t) + \bw(t)\footnote{Continuous-time white noise can be defined as the \textit{weak} derivative of a continuous Brownian motion, $\int_a^b \bw(t)dt = \int_a^b d\bbeta(t)$. Because paths of Brownian motion, $\bbeta$, have fractal structure, with slopes up and down of arbitrarily large magnitude over infinitesimal intervals, they are too rough to differentiate~\cite{lorig, sarkka-sdes}. It is therefore mathematically imprecise to write $\bw(t)$ outside an integral; it is a \textit{distribution}, not a function (see \autoref{fig:weak-derivatives}). However, it is common in physics and engineering to write stochastic differential equations heuristically rather than in differential form (multiplying through by $dt$) or integral form.}$$

\noindent where script $\mathcal{A}$ and $\mathcal{B}$ are the continuous-time state evolution and control matrices, respectively, and $\bx$, $\bu$, and $\bw$ now vary continuously. Note that sampled observations are still discrete and obey $\by(t_n) = \bC\bx(t_n) + \bv_n,\ \bv_n \sim \mathcal{N}(0, \bR)$, just as before.

We would like to know $\bx(t_n + \Delta t)$ for some time in the future, so we solve the above with the method of integrating factor. We make the simplifying assumption that $\bu$ is constant over this interval, so-called ``zero-order hold", which is valid if control inputs or synchronized data streams only get updates at sampling times, a reasonable assumption if they come from a computed control or sensing loop. We substitute $\tau = t - t_n$ because all terms will become integrands, and integrating from 0 simplifies the math.
$$\frac{d\bx}{d\tau} = \mathcal{A}\bx(\tau) + \mathcal{B}\bu + \bw(\tau)
\rightarrow \underbrace{e^{-\mathcal{A}\tau}\frac{d\bx}{d\tau} - e^{-\mathcal{A}\tau}\mathcal{A}\bx(\tau)}_{=\frac{d}{d\tau}\Big(e^{-\mathcal{A}\tau} \bx(\tau)\Big)} = e^{-\mathcal{A}\tau}\big(\mathcal{B}\bu + \bw(\tau)\big)$$
\begin{align*}
&\rightarrow \underbrace{\int\limits_0^{\Delta t} \frac{d}{d\tau}\Big(e^{-\mathcal{A}\tau} \bx(\tau)\Big) d\tau}_{\mathclap{e^{-\mathcal{A}\tau} \bx(\tau)\Big|_0^{\Delta t} = e^{-\mathcal{A}\Delta t} \bx(\Delta t) - \mathbb{I}\bx(0)}} = \int\limits_0^{\Delta t} e^{-\mathcal{A}\tau}\big(\mathcal{B}\bu + \bw(\tau)\big) d\tau\\
&\rightarrow e^{-\mathcal{A}\Delta t}\bx(\Delta t) = \bx(0) + \int\limits_0^{\Delta t} e^{-\mathcal{A}\tau}\mathcal{B} d\tau\ \bu + \int\limits_0^{\Delta t} e^{-\mathcal{A}\tau}\bw(\tau) d\tau\\
&\rightarrow \bx(\Delta t) = e^{\mathcal{A}\Delta t}\bx(0) + \int\limits_0^{\Delta t} e^{\mathcal{A}(\Delta t-\tau)}\mathcal{B}d\tau\ \bu + \underbrace{\int\limits_0^{\Delta t} e^{\mathcal{A}(\Delta t-\tau)}\bw(\tau)d\tau}_{\bw_d(\Delta t) \sim \mathcal{N}}\footnotemark
\end{align*}
\footnotetext{This latter integral can be written w.r.t.~a Brownian motion differential, $d\bbeta$. Here the integrand is deterministic, so we can use $\int_0^Tg(t)d\bbeta(t) \sim \mathcal{N}\big(0, \int_0^T g^2(t)dt\big)$~\cite{lorig}. More generally, integrands can depend on $\bbeta$, as $\int_{t_0}^{t_{N-1}} f\big(t, \bbeta(t)\big)d\bbeta(t)= \lim_{N \to \infty} \sum_{k=0}^{N-1} f\big(t_k^*, \bbeta(t_k^*)\big)\big(\bbeta(t_{k+1}) - \bbeta(t_k)\big)$ for $t_k^* \in [t_k, t_{k+1}]$. The value of this expression depends on the choice of $t_k^*$, because the nonzero quadratic variation of Brownian motion leads to a $(\Delta\bbeta)^2$ term in the Taylor expansion of stochastic $f$ around $t_k$, whereas for smooth functions quadratic variation vanishes. To make the differential equation well-defined (with a unique solution), a common choice is $t_k^* = t_k$, which eliminates the quadratic term and yields the It\^o integral~\cite{lorig, sarkka-sdes}.}

\noindent where $\bw_d$ is discrete-time noise. The middle term can be simplified a bit with a change of variables: Letting $\tilde\tau = \Delta t-\tau$, then when $\tau = 0,\ \tilde\tau = \Delta t$, and when $\tau= \Delta t,\ \tilde\tau = 0$, and $\frac{d\tilde\tau}{d\tau} = -1 \rightarrow d\tau = -d\tilde\tau$.
$$\int\limits_0^{\Delta t} e^{\mathcal{A}(\Delta t-\tau)}\mathcal{B}d\tau = \int
\limits_{\Delta t}^0 e^{\mathcal{A}\tilde\tau}\mathcal{B}(-d\tilde\tau) = \int\limits_0^{\Delta t} e^{\mathcal{A}\tilde\tau}\mathcal{B}d\tilde\tau = \bB(\Delta t)$$

\noindent where the bold symbol denotes the discrete-time matrix.

The noise term, $\bw_d$, is a zero-mean Gaussian random variable because integrating $\bw(\tau)$ against $e^{\mathcal{A}(\Delta t-\tau)}$, which is just another matrix, constitutes a linear functional, which preserves the mean and Gaussian character of the process, only stretching its covariance. We can find this new covariance as:
\begin{equation}\label{eqn:discrete-cov}
\begin{aligned}
&\mathbb{E}[\bw_d(\Delta t)\bw_d(\Delta t)^\mathrm{T}] - \overset{\text{\normalsize 0}}{\cancel{\mathbb{E}[\bw_d]\mathbb{E}[\bw_d^\mathrm{T}]}} = \mathbb{E}\Big[\int\limits_0^{\Delta t} e^{\mathcal{A}(\Delta t-\tau)} \bw(\tau)\bw(\tau)^\mathrm{T} e^{\mathcal{A}^\mathrm{T}(\Delta t-\tau)}d\tau\Big]\\
&\quad\quad\quad= \int\limits_0^{\Delta t} e^{\mathcal{A}(\Delta t-\tau)} \underbrace{\mathbb{E}[\bw(\tau)\bw(\tau)^\mathrm{T}]}_{\mathclap{\text{constant, continuous-time } \mathcal{Q}}} e^{\mathcal{A}^\mathrm{T}(\Delta t-\tau)}d\tau = \int\limits_0^{\Delta t} e^{\mathcal{A}\tilde{\tau}} \mathcal{Q}e^{\mathcal{A}^\mathrm{T}\tilde\tau}d\tilde\tau = \bQ(\Delta t)
\end{aligned}
\end{equation}

\noindent where $\mathcal{Q}$ is constant because it is the expected value of a stationary random matrix, which is deterministic, irrespective of $\tau$, and in the last step we have done the same $\tilde\tau = \Delta t - \tau$ trick as before.

We can now find the estimated state and error covariance for time offset $\Delta t$ by taking expected values:
\begin{equation}\label{eqn:variable-dt-x}
\begin{aligned}
\hbx(\Delta t) &= \mathbb{E}[\bx(\Delta t)] = \mathbb{E}[e^{\mathcal{A}\Delta t}\bx(0) + \bB(\Delta t)\bu + \bw_d(\Delta t)]\\
&= e^{\mathcal{A}\Delta t}\mathbb{E}[\bx(0)] + \mathbb{E}[\bB(\Delta t)\bu] + \mathbb{E}[\bw_d(\Delta t)] = e^{\mathcal{A}\Delta t}\hbx(0) + \bB(\Delta t)\bu
\end{aligned}
\end{equation}
\begin{equation}\label{eqn:variable-dt-P}
\begin{aligned}
\bP(\Delta t) &= \mathbb{E}\big[\big(\bx(\Delta t) - \mathbb{E}[\bx(\Delta t)]\big)\big(\bx(\Delta t) - \mathbb{E}[\bx(\Delta t)]\big)^\mathrm{T}\big]\\
&= \mathbb{E}\big[\big(e^{\mathcal{A}\Delta t}\bx(0) + \cancel{\bB(\Delta t)\bu} + \bw_d(\Delta t) - e^{\mathcal{A}\Delta t}\mathbb{E}[\bx(0)] -\cancel{\bB(\Delta t)\bu}- \mathbb{E}[\bw_d(\Delta t)]\big)\big(\text{ditto}\big)^\mathrm{T}\big]\\
&=\mathbb{E}\big[\big(e^{\mathcal{A}\Delta t}\big(\bx(0) -\mathbb{E}[\bx(0)]\big) + \big(\bw_d(\Delta t) - \mathbb{E}[\bw_d(\Delta t)]\big)\big)\big(\text{ditto}\big)^\mathrm{T}\big]\\
&=e^{\mathcal{A}\Delta t} \underbrace{\mathbb{E}\big[\big(\bx(0) -\mathbb{E}[\bx(0)]\big) \big(\bx(0) -\mathbb{E}[\bx(0)]\big)^\mathrm{T}\big]}_{P(0)}e^{\mathcal{A}^\mathrm{T}\Delta t} + \overset{\text{\normalsize 0}}{\cancel{\text{cross terms}}} + \bQ(\Delta t)
\end{aligned}
\end{equation}

\noindent where cross terms are 0 because noise is assumed to be uncorrelated with state estimation error.

By treating any $\bx_n, \bP_n$ as initial conditions through the time shift $\tau = t - t_n$, we can use Equations \ref{eqn:variable-dt-x} and \ref{eqn:variable-dt-P} to find $\bx_{n+1},\bP_{n+1}$ for arbitrary step size. All we have left to do is compute discrete-time $\bA(\Delta t) = e^{\mathcal{A}\Delta t}$, $\bB(\Delta t)$, and $\bQ(\Delta t)$ given continuous-time $\mathcal{A}$, $\mathcal{B}$, and $\mathcal{Q}$. The state evolution matrix exponential can be taken directly, but the discrete control matrix and process noise covariance involve integrals. Thankfully, there is a very elegant trick~\cite{integration-via-expM, axelssonandgustafsson, crassidis_junkins} to integrate by exponentiating a special matrix, which allows us to solve for all these quantities at once.

First consider the generic block triangular matrix:

$$\mathcal{M} = \begin{bmatrix} \mathcal{M}_1 & \mathcal{M}_2 \\ 0 & \mathcal{M}_3\end{bmatrix}$$

\noindent Taking a matrix exponential is equivalent to plugging the matrix into the Taylor expansion:

$$e^{\mathcal{M}t} = \mathbb{I} + \mathcal{M}t + \frac{\mathcal{M}^2t^2}{2!} + \frac{\mathcal{M}^3t^3}{3!} + \cdots$$

\noindent The sum of integer powers of a block triangular matrix remains block triangular, so the result has structure:
$$\bM(t) = \begin{bmatrix} \bM_1(t) & \bM_2(t) \\ 0 & \bM_3(t)\end{bmatrix}$$

We can thus set up a differential equation:
\begin{align*}
&\frac{d}{dt}e^{\mathcal{M}t} = \mathcal{M}e^{\mathcal{M}t} = \mathcal{M}\bM(t)\\
&\rightarrow \frac{d}{dt}\begin{bmatrix} \bM_1 & \bM_2 \\ 0 & \bM_3\end{bmatrix} = \begin{bmatrix} \mathbf{\dot{M}}_1 & \mathbf{\dot{M}}_2 \\ 0 & \mathbf{\dot{M}}_3\end{bmatrix} = \begin{bmatrix} \mathcal{M}_1 & \mathcal{M}_2 \\ 0 & \mathcal{M}_3\end{bmatrix}\begin{bmatrix} \bM_1(t) & \bM_2(t) \\ 0 & \bM_3(t)\end{bmatrix}
\end{align*}

\noindent where over-dots represent the derivative, which can be safely distributed to each block independently. By equating blocks, we can turn this into three equations:

$$\mathbf{\dot{M}}_1 = \mathcal{M}_1\bM_1,\quad \mathbf{\dot{M}}_2 = \mathcal{M}_1\bM_2 + \mathcal{M}_2\bM_3,\quad \mathbf{\dot{M}}_3 = \mathcal{M}_3\bM_3$$

\noindent These can be solved by integrating factor. Blocks 1 and 3 have the same governing equation, so we show the solution procedure only of block 1:
\begin{align*}
&\mathbf{\dot{M}}_1 = \mathcal{M}_1\bM_1 \rightarrow \mathbf{\dot{M}}_1 - \mathcal{M}_1\bM_1 = 0 \rightarrow e^{-\mathcal{M}_1\tau}\mathbf{\dot{M}}_1 - e^{-\mathcal{M}_1\tau}\mathcal{M}_1\bM_1 = 0\\
&\rightarrow \frac{d}{d\tau}\big(e^{-\mathcal{M}_1 \tau}\bM_1(\tau)\big) = 0 \rightarrow \int\limits_0^t\frac{d}{d\tau}\big(e^{-\mathcal{M}_1 \tau}\bM_1\big)d\tau = \int\limits_0^t 0d\tau\\
& \rightarrow e^{-\mathcal{M}_1\tau}\bM_1(\tau)\Big|_0^t = 0 \rightarrow e^{-\mathcal{M}_1t}\bM_1(t) - \mathbb{I} \underbrace{\bM_1(0)}_{\mathclap{\text{on-diagonal block of } e^{\mathcal{M}0} = \mathbb{I}}} = 0 \rightarrow \bM_1(t) = e^{\mathcal{M}_1 t}
\end{align*}

\noindent For block 2:
\begin{align*}
&\mathbf{\dot{M}}_2 - \mathcal{M}_1 \bM_2= \mathcal{M}_2 \bM_3 \rightarrow \int\limits_0^t \frac{d}{d\tau}(e^{-\mathcal{M}_1\tau} \bM_2)d\tau = \int\limits_0^t e^{-\mathcal{M}_1\tau}\mathcal{M}_2 \bM_3d\tau\\
&\rightarrow e^{-\mathcal{M}_1\tau}\bM_2(\tau)\Big|_0^t = e^{-\mathcal{M}_1 t}\bM_2(t) - \mathbb{I}\underbrace{\bM_2(0)}_{\mathclap{0 \text{ because off-diagonal in } e^{\mathcal{M}0}}}\\
&\rightarrow \bM_2(t) = \int\limits_0^t e^{\mathcal{M}_1(t - \tau)}\mathcal{M}_2 \bM_3(\tau)d\tau = \int\limits_0^t e^{\mathcal{M}_1\tilde\tau}\mathcal{M}_2 \bM_3(t - \tilde\tau)d\tilde\tau
\end{align*}

\noindent where in the last step we use the same variable transformation as before. Substituting the solution for $\bM_3(t)$, we get:
$$\bM_2(t) = \int\limits_0^t e^{\mathcal{M}_1\tilde \tau}\mathcal{M}_2 e^{\mathcal{M}_3(t - \tilde\tau)}d\tilde\tau$$

Notice that the solution for blocks 1 and 3 looks like the basic matrix exponentiation we have to do to find $\bA$, and the solution for block 2 looks suspiciously like the integral expression for $\bQ$ from \autoref{eqn:discrete-cov}. Let's choose $\mathcal{M}_1 = \mathcal{A}$, $\mathcal{M}_2 = \mathcal{Q}$, and $\mathcal{M}_3 = -\mathcal{A}^\mathrm{T}$. We then form the matrix:
$$\mathcal{M} = \begin{bmatrix} \mathcal{A} & \mathcal{Q} \\ 0 & -\mathcal{A}^\mathrm{T} \end{bmatrix}$$

\noindent and exponentiate to produce:
$$\bM(\Delta t) = \begin{bmatrix} e^{\mathcal{A}\Delta t} & \int\limits_0^{\Delta t} e^{\mathcal{A}\tilde\tau} \mathcal{Q} e^{\mathcal{A}^\mathrm{T}(\tilde\tau - \Delta t)} d\tilde\tau\\ 0 & e^{-\mathcal{A}^\mathrm{T}\Delta t}\end{bmatrix} = \begin{bmatrix} \bA(\Delta t) & \bQ(\Delta t)\bA\!^{-\mathrm{T}}(\Delta t) \\ 0 & \bA\!^{-\mathrm{T}}(\Delta t)\end{bmatrix}$$

\noindent where $^{-\mathrm{T}}$ denotes the inverse transpose. Similarly, choosing $\mathcal{M}_1 = \mathcal{A}$, $\mathcal{M}_2 = \mathcal{B}$, and $\mathcal{M}_3 = 0$ yields:
$$\bM(\Delta t) = \begin{bmatrix} \bA(\Delta t) & \bB(\Delta t) \\ 0 & \mathbb{I}\end{bmatrix}$$

\noindent So we can simply pick out a couple of the matrices we need, and find the other by multiplication since $\bA\!^{-\mathrm{T}}\bA^\mathrm{T} = e^{-\mathcal{A}^\mathrm{T} \Delta t}e^{\mathcal{A}^\mathrm{T}\Delta t} = \mathbb{I}$.

For all the math that got us here, this is very short to code and fairly cheap to run: For a model of dimension $\nu$ the matrix $\mathcal{M} \in \mathbb{R}^{2\nu \times 2\nu}$ costs $O(\nu^3)$ to exponentiate~\cite{matrix-exp}, but both the Kalman filter (\autoref{algo:kalman-filter-algo}) and RTS smoother (\autoref{eqn:rts}) are already $O(\nu^3)$ per iteration due to dense matrix multiplications or inversions, so variable time steps are slightly more complicated and expensive but do not actually worsen the asymptotic complexity.

RTS smoothing extends to irregular steps just as easily as the Kalman filter, because it simply takes the history of all \textit{a priori} and \textit{a posteriori} $\bx$ and $\bP$ estimates and combines distributions, agnostic to whether they come from Equations \ref{eqn:variable-dt-x} and \ref{eqn:variable-dt-P} or from the update equations in \autoref{algo:kalman-filter-algo}. The exponentiation causes covariances $\bQ$ and $\bP$ to grow continuously between steps, but the optimal rule for combining Gaussians remains the same regardless of their size.

\section{Conclusions}
\label{sec:conclusion}


Numerical differentiation is a fundamental operation that appears across all scientific and engineering domains, part of the essential computational infrastructure for extracting insight from measurements and simulations. Ultimately, success demands matching method assumptions to problem characteristics. To that end, this review has organized the landscape of differentiation scenarios into \textit{five families}: static, analytic relationships best handled by Automatic Differentiation; regular and irregular simulations where Spectral Methods and Finite Elements excel, respectively; and cases of noisy data that bifurcate into model-based approaches leveraging Kalman filtering frameworks and model-free smoothing techniques.

The key takeaways are that (1) more restrictive assumptions enable better performance, because more efficient, accurate methods demand more structure, and the corollary that (2) where structure is available, the choice of method is largely forced by algorithms' strengths and limitations. Periodic signals unlock the power of Fourier spectral methods; known system dynamics enable optimal Kalman smoothing; and even naive constant-derivative models provide remarkably effective regularization. Thus the upper paths in \autoref{fig:flowchart} lead to definite conclusions.

For practitioners facing noisy data without prior models, the loss function framework presented in Equations \ref{eqn:floris-cost}–\ref{eqn:pavel-cost}, combined with the bandlimit heuristic of \autoref{eqn:gamma-heuristic}, enables principled hyperparameter optimization to balance fidelity to data with a smoothness prior. Solutions consistently approach Pareto-optimal, reducing the burden from navigating dozens of method-specific parameters to selecting a single smoothness parameter informed by signal bandwidth. In most applications a balanced estimate is a good starting point, but different end uses can tip the balance---e.g.~control inputs favoring smoothness to prevent high control efforts~\cite{crassidis_junkins}, even at the cost of higher error, or data-driven modeling applications preferring low error correlation at the cost of RMSE to prevent overemphasis of low frequency dynamics~\cite{sindy}.

Our comprehensive experimental comparison of optimized methods on noisy data reveals that most achieve similar accuracy under ideal conditions, with auxiliary considerations often proving decisive (\autoref{fig:flexibility}). \texttt{RTSDiff} (\autoref{sec:kalman-constant-deriv}) emerges as a versatile general-purpose choice, offering a high performance floor alongside exceptional flexibility for variable step sizes and data interpolation. For particularly large data gaps, \texttt{LinearDiff} and \texttt{PolyDiff} (\autoref{sec:polyfit}) degrade most gracefully. \texttt{LinearDiff} is most resilient against poor choices of bandlimit, and \texttt{RobustDiff} provides superior outlier rejection at the cost of substantially increased optimization computation. Occasionally, prior knowledge suggests a particular choice, such as \texttt{TVRDiff} (\autoref{sec:tvr}) to identify piecewise derivatives, \texttt{LinearDiff} if the data was generated by a linear ODE, \texttt{SpectralDiff} where frequency plays a central role, or \texttt{WaveletDiff} for multi-scale features. One can even apply multiple methods to the same data to reveal challenging regions of a series where methods disagree or to corroborate hyperparameter choices.

Greater challenges remain broadly open, including nonstationary noise, adaptive priors, and uncertainty quantification. A parallel literature in control~\cite{moja-one-sided} addresses one-sided estimation, only incidentally touched here (\autoref{sec:noisy-with-knowledge}). But there is a broad pattern worth recognizing, which sets a fundamental limit on this field: Every method in this review works by imposing an assumption---a bandlimit, a smoothness prior, an interpolation form, a model, governing equations---and where noise shares the signal's character so no such assumption can separate the two, the problem is underdetermined rather than merely hard, and no method will be able to recover a true signal or derivative estimate. Knowing this limit can be valuable in itself, a sign marking where to quit.

\section*{Acknowledgments}  
We are especially grateful to Sasha Aravkin for extensive discussions on convex optimization techniques related to robust and nonsmooth formulations. We would also like to thank Jing Yu for lending her clarity and depth to guide us through some of the complex relationships among robust control methods. And finally we recognize Alexey Yermakov for thoroughly proofreading an early draft. The authors acknowledge support from the National Science Foundation AI Institute in Dynamic Systems (grant number 2112085).

\bibliographystyle{siamplain}
\bibliography{references.bib}

\end{document}